\documentclass[a4paper,11pt,british,oneside]{article}

\usepackage[utf8]{inputenc}
\usepackage[dvipsnames,svgnames,x11names,hyperref]{xcolor}
\usepackage{amsmath,amssymb,csquotes,enumerate,dsfont,babel,mathtools,textgreek,hyphenat,needspace,microtype}
\usepackage[amsmath,amsthm,thmmarks,hyperref,thref]{ntheorem}
\definecolor{linkcol}{rgb}{0,0.2,0.6}
\definecolor{citecol}{rgb}{0,0.4,0}
\definecolor{green}{rgb}{0.67,0.77,0}
\usepackage[square,sort,numbers]{natbib}
\setcitestyle{citesep={;}}
\usepackage[pdfpagelabels,bookmarks,colorlinks,linkcolor=linkcol,citecolor=citecol,urlcolor=black]{hyperref}
\newcounter{thcounter}
\numberwithin{thcounter}{section}
\numberwithin{equation}{section}

\theoremstyle{break}
\newtheorem{lemma}[thcounter]{Lemma}
\newtheorem{proposition}[thcounter]{Proposition}
\newtheorem{theorem}[thcounter]{Theorem}
\newtheorem{corollary}[thcounter]{Corollary}

\newtheorem{question}[thcounter]{Question}
\theorembodyfont{\upshape}
\newtheorem{definition}[thcounter]{Definition}
\newtheorem{remark}[thcounter]{Remark}
\newtheorem{example}[thcounter]{Example}

\usepackage{wrapfig,lscape,rotating}

\providecommand\numberthis{\addtocounter{equation}{1}\tag{\theequation}}

\providecommand{\NN}{\mathbb{N}}
\providecommand{\ZZ}{\mathbb{Z}}
\providecommand{\QQ}{\mathbb{Q}}
\providecommand{\CC}{\mathbb{C}}

\providecommand{\GG}{\mathbb{G}}

\mathchardef\mhyphen="2D

\newcommand{\Tr}{\operatorname{Tr}}

\newcommand{\Hom}{\operatorname{Hom}}
\newcommand{\Tor}{\operatorname{Tor}}
\newcommand{\Ext}{\operatorname{Ext}}
\newcommand{\Irr}{\operatorname{Irr}}
\newcommand{\mult}{\operatorname{mult}}
\newcommand{\onb}{\operatorname{onb}}
\newcommand{\op}{\mathrm{op}}

\newcommand{\id}{\mathrm{id}}

\newcommand{\Rep}{\operatorname{Rep}}
\providecommand{\ryosup}[3]{\mathord{\raisebox{0.8ex}[0ex][0ex]{\scriptsize $#1$}{#2}\hspace{-0.25ex}\raisebox{0.8ex}[0ex][0ex]{\scriptsize $#3$}}}
\providecommand{\lsup}[2]{\ryosup{#1}{#2}{}}

\newcommand{\I}{\mathds{1}}
\newcommand{\im}{\operatorname{im}}

\newcommand{\diff}{\mathop{\mathrm{\mathstrut d}}\!}

\renewcommand{\bar}{\overline}
\renewcommand{\tilde}{\widetilde}
\renewcommand{\hat}{\widehat}
\makeatletter
\let\@temp\phi
\let\phi\varphi
\let\varphi\@temp

\let\@temp\epsilon
\let\epsilon\varepsilon
\let\varepsilon\@temp
\makeatother

\newcommand{\Fun}{\operatorname{Fun}}
\newcommand{\Orb}{\operatorname{Orb}}
\newcommand{\BHom}{\overline{\mathrm{Hom}}}
\newcommand{\BEnd}{\overline{\mathrm{End}}}

\newcommand{\Tot}{\operatorname{Tot}}
\newcommand{\counit}{\varrho}
\newcommand{\cftimes}[1]{\mathbin{\overline{\otimes}_{#1}}}
\newcommand{\cfpow}[2]{\overline{\otimes}_{#1}^{\,#2}}
\newcommand{\htimes}{\mathbin{\overline{\otimes}}}
\newcommand{\ot}{\leftarrow}
\newcommand{\cone}{\operatorname{cone}}
\begin{document}

\def\sectionautorefname{section}
\def\subsectionautorefname{section}
\def\subsubsectionautorefname{section}

\begin{center}
    {\LARGE\boldmath\bf $L^2$-Betti numbers of $C^*$-tensor categories associated with totally disconnected groups}

    \bigskip

    {\sc Matthias Valvekens}
    
    \bigskip
    
    \textit{E-mail address:} \texttt{matthias.valvekens@kuleuven.be}\\[0.5ex]
    {\small\sc KU~Leuven, Department of Mathematics\\ Celestijnenlaan 200B, B-3001 Leuven, Belgium}
\end{center}
\begin{abstract}
\noindent We prove that the $L^2$-Betti numbers of a rigid $C^*$-tensor category vanish in the presence of an almost-normal subcategory with vanishing $L^2$-Betti numbers, generalising a result of \cite{bfs}.
We apply this criterion to show that the categories constructed from totally disconnected groups in \cite{arano-vaes} have vanishing $L^2$-Betti numbers.
Given an almost-normal inclusion of discrete groups $\Lambda<\Gamma$, with $\Gamma$ acting on a type $\mathrm{II}_1$ factor $P$ by outer automorphisms, we relate the cohomology theory of the quasi-regular inclusion $P\rtimes\Lambda\subset P\rtimes\Gamma$ to that of the Schlichting completion $G$ of $\Lambda<\Gamma$.
If $\Lambda<\Gamma$ is unimodular, this correspondence allows us to prove that the $L^2$-Betti numbers of $P\rtimes\Lambda\subset P\rtimes\Gamma$ are equal to those of $G$.
\end{abstract}
\section{Introduction}

The theory of $L^2$-invariants originates with Atiyah, who defined $L^2$-Betti numbers for $\Gamma$-coverings of compact manifolds \cite{atiyah} as a tool to study elliptic differential operators.
In this setting, one has a discrete group $\Gamma$ acting freely on a noncompact manifold $\bar{X}$ with compact quotient $X=\bar{X}/\Gamma$.
If the covering space $\bar{X}$ is contractible, these $L^2$-Betti numbers are invariants of $\Gamma$.
Later, Cheeger and Gromov introduced an $L^2$-cohomology theory for actions of discrete groups on $CW$-complexes \cite{cheeger-gromov}, which gives rise to $L^2$-Betti numbers $\beta^{(2)}_n(\Gamma)$ for arbitrary countable discrete groups $\Gamma$ via the action of $\Gamma$ on its classifying space.
Over the years, $L^2$-invariants have found application in many parts of mathematics, see e.g.\@ \cite{gaboriau,lueck-l2-invariants}.
Given these connections with other research areas, considerable effort has been made to define $L^2$-invariants for various kinds of ``quantum'' symmetries.
For example, given a compact quantum group $\mathbb{G}$ of Kac type, one can take the $L^2$-Betti numbers $\beta^{(2)}_n(\hat{\mathbb{G}})$ of its discrete dual in the sense of \cite{kyed-l2betti}.
Various familiar results about $L^2$-Betti numbers of discrete groups carry over to the quantum setting---see e.g.\@ \cite{kyed-coamenable,bkr16,krvv}---although concrete computations are in general much more difficult.

Similarly, $L^2$-Betti numbers have also been defined for \textit{rigid $C^*$-tensor categories}, with the introduction of a Hochschild-type (co)homology theory for such categories in \cite{psv-cohom}.
Taking (co)homology with coefficients in the regular representation and measuring the dimension of the resulting spaces in the appropriate way, one then obtains a definition of $L^2$-Betti numbers.
In fact, \cite{psv-cohom} formulates these concepts in the more general framework of \textit{quasi-regular inclusions} of type $\mathrm{II}_1$ factors.
The first main result of this paper is a vanishing criterion for the $L^2$-Betti numbers of a rigid $C^*$-tensor category, which we apply to the categories constructed from totally disconnected groups in \cite{arano-vaes}.
When the categories involved are actual discrete groups, our vanishing criterion reduces to a known result of \cite{bfs}.
In the second half, we discuss the cohomology theory and $L^2$-Betti numbers of quasi-regular inclusions coming from \textit{discrete Hecke pairs}, i.e.\@ inclusions of discrete groups $\Lambda<\Gamma$ such that $[\Lambda:\Lambda\cap g\Lambda g^{-1}]<\infty$ for all $g\in\Gamma$.
Concretely, given an outer action of $\Gamma$ on a type $\mathrm{II}_1$ factor $P$, the crossed product inclusion $P\rtimes \Lambda \subset P\rtimes\Gamma$ is quasi-regular.
The second main result of this paper essentially identifies the cohomology theory of $P\rtimes\Lambda\subset P\rtimes\Gamma$ with the continuous cohomology of the \textit{Schlichting completion} $G$ of $\Lambda<\Gamma$.
When $G$ is unimodular, this correspondence also recovers the $L^2$-Betti numbers $\beta^{(2)}_n(G)$, as defined in \cite{kpv13,hdp-thesis}.

Rigid $C^*$-tensor categories come about naturally in the study of quantum symmetries of various kinds, providing a unifying language to talk about discrete groups, representation categories of compact (quantum) groups and standard invariants of finite-index subfactors.
The study of rigid $C^*$-tensor categories is in many respects informed and inspired by progress made in the theory of discrete groups over the past decades.
Following this principle, the fact that approximation and rigidity properties like amenability, property~(T) and the Haagerup property have played a central role in understanding discrete groups has spurred efforts to formulate these properties in the framework of $C^*$-tensor categories.
In the context of subfactor standard invariants, this programme was first carried out through Popa's \textit{symmetric enveloping algebra} \cite{popa-symm-env-1,popa-symm-env-2,popa-annals-01,brothier-weak-am}, and later through the development of a (unitary) representation theory for rigid $C^*$-tensor categories \cite{pv-repr-subfactors,neshveyev-yamashita,ghosh-jones}; see also \cite{psv-cohom,adlw-1,adlw-2}. 
The recent developments concerning $L^2$-Betti numbers fit naturally into this line of research. 

Already in \cite{psv-cohom}, various computational properties of $L^2$-Betti numbers for quasi-regular inclusions were established, including a K\"unneth-type formula for tensor products, a free product formula, and a Cheeger--Gromov-type vanishing result for amenable quasi-regular inclusions.
All of these specialise to analogous statements in the context of rigid $C^*$-tensor categories.
Later, \cite{krvv} established that $L^2$-Betti numbers of rigid $C^*$-tensor categories scale along finite-index inclusions, and that the $L^2$-Betti numbers of the dual $\hat{\mathbb{G}}$ of a Kac-type compact quantum group $\mathbb{G}$ agree with those of its representation category $\Rep_f(\mathbb{G})$.
Moreover, since the $L^2$-Betti numbers of a rigid $C^*$-tensor category actually only depend on the Morita equivalence class of the category \cite[Proposition~7.4]{psv-cohom}, this allows for quite some flexibility in computing $L^2$-Betti numbers of discrete quantum groups.
Taken together, these results provide computations of the $L^2$-Betti numbers of several discrete quantum groups and representation categories of compact quantum groups, including the free unitary quantum groups of Van Daele and Wang \cite{vandaele-wang-uqg} and various wreath products \cite[Theorem~5.2]{krvv}; see also \cite{bkr16,kyed-raum-l2betti}.

Let $\mathcal{C}$ be a full $C^*$-tensor subcategory of a rigid $C^*$-tensor category $\mathcal{D}$.
The goal of the first part of this article is to show that if $\mathcal{C}$ has vanishing $L^2$-Betti numbers and $\mathcal{C}$ is (almost) normal in $\mathcal{D}$ (see \thref{def:normal-almost-normal} below), then the larger category $\mathcal{D}$ must also have vanishing $L^2$-Betti numbers.
This generalises a result of Bader, Furman and Sauer \cite[Corollary~1.4]{bfs}.

Given a totally disconnected locally compact group $G$ with a compact open subgroup $K$, Arano and Vaes \cite{arano-vaes} define a rigid $C^*$-tensor category $\mathcal{C}_f(K<G)$.
For all such $K<G$, the category $\Rep_f(K)$ of finite-dimensional unitary representations of $K$ sits inside $\mathcal{C}_f(K<G)$ as a normal subcategory.
Since $\Rep_f(K)$ has vanishing $L^2$-Betti numbers, our criterion yields that the $L^2$-Betti numbers of $\mathcal{C}_f(K<G)$ also vanish.
Our results also imply that certain bicrossed product discrete quantum groups (in particular those considered in \cite{vv-property-t}) have vanishing $L^2$-Betti numbers.

In the second part of the article, we study a related but fundamentally different problem.
Given an almost-normal inclusion $\Lambda<\Gamma$ of discrete groups with $\Gamma$ acting on a type $\mathrm{II}_1$ factor $P$ by outer automorphisms, the crossed product inclusion $P\rtimes\Lambda\subset P\rtimes\Gamma$ is quasi-regular and irreducible.
In this context, it commonly happens that analytic properties of $\Lambda<\Gamma$ (and by extension $P\rtimes\Lambda\subset P\rtimes\Gamma$) correspond to analytic properties of the Schlichting completion $G$ of $\Lambda<\Gamma$, see e.g.\@ \cite{ad-hal,larsen-palma,popa-annals-01}.
Since the Schlichting completion $G$ of $\Lambda<\Gamma$ is unimodular whenever the inclusion $P\rtimes\Lambda\subset P\rtimes\Gamma$ is, one can ask the same question about the $L^2$-Betti numbers $\beta^{(2)}_n(G)$ and $\beta^{(2)}_n(P\rtimes\Lambda\subset P\rtimes\Gamma)$.
We show that both sequences of $L^2$-Betti numbers agree, provided that one takes $\beta^{(2)}_n(G)$ with respect to an appropriate normalisation of the Haar measure on $G$.

\subsection*{Acknowledgements}
The author would like to express his gratitude to his advisor Stefaan~Vaes for suggesting this topic and for sharing his insights on the vanishing criterion for almost-normal subgroups on the one hand, but also for pointing out several technical mistakes in an earlier draft of this paper, and for providing extremely valuable feedback throughout the process.
He is also grateful to Adam~Skalski for the hospitality extended during the author's stay at IMPAN, where part of the work on this project was done.
This work was partially supported by European Research Council Consolidator Grant 614195 RIGIDITY, by long term structural funding~-- Methusalem grant of the Flemish Government, and by the FWO--PAS project VS02619N: von Neumann algebras arising from quantum symmetries.

\section{Preliminaries}
\subsection{Dimension theory}
\subsubsection{Elementary notions}
    We frequently make use of the L\"uck dimension for von Neumann algebras with a normal, semifinite and faithful (abbreviated n.s.f.) trace.
    This dimension function (due to L\"uck \cite{lueck-dimtheo-1} in the finite setting) generalises the Murray--von Neumann dimension to arbitrary modules over von Neumann algebras equipped with a faithful n.s.f.\@ trace.
    Throughout, all von Neumann algebras we deal with will have separable predual, and all Hilbert spaces will be separable.
    For an overview of the theory in the case of von Neumann algebras equipped with a normal faithful tracial \textit{state}, we refer to \cite[\S~6.1]{lueck-l2-invariants}.
    The definition in the finite case can be restated essentially verbatim for semifinite traces \cite[Definition~B.17]{hdp-thesis}.
    \begin{definition}[\cite{hdp-thesis}]
        Let $M$ be a von Neumann algebra with n.s.f.\@ trace $\tau$, and $H$ a right $M$-module.
        Then define
        \begin{align*}
            \dim_{-M}(H) = \sup\{(\Tr\otimes\tau)(P)\mid\, &P\in M_n(\CC)\otimes M, (\Tr\otimes\tau)(P)<\infty, \\
            &P(\CC^n\otimes M)\hookrightarrow H\text{ as }M\text{-modules} \}\,.
        \end{align*}
        The definition for left modules is analogous.
    \end{definition}
    We refer to \cite[\S~5.4]{kpv13} and \cite[Appendix~B]{hdp-thesis} for an in-depth discussion of this dimension function in the semifinite setting.
    If there is no possibility of confusion as to whether we take the dimension from the left or from the right, we often simply write $\dim_M$ instead of $\dim_{M-}$ or $\dim_{-M}$.
    As explained in \cite[Remark~3.8]{krvv}, the dimension function $\dim_{-M}$ can alternatively be defined as
    \[
        \dim_{-M}(H) = \sup\{ \tau(p)\dim_{-pMp}(Hp)\mid p\in M, \tau(p)<\infty\}\,,
    \]
    where we consider the corner $pMp$ equipped with the tracial state given by $\tau(p)^{-1}\tau(\cdot)$.
    This observation reduces many questions about dimension theory for semifinite traces to the finite case.

    This is especially true when working with so-called \textit{locally finite} modules, which will cover all modules of interest to us.
    \begin{definition}
        Let $M$ be a von Neumann algebra with n.s.f.\@ trace $\tau$.
        A right (resp.\@ left) $M$-module $H$ is \textit{locally finite} if for all $\xi\in H$ there exists a projection $p\in M$ such that $\tau(p)<\infty$ and $\xi p=\xi$ (resp.\@ $p\xi=\xi$).
    \end{definition}
    The category of locally finite (left or right) $M$-modules is closed under passage to subobjects, quotients and extensions.

    In homological computations, it is often useful to consider dimension isomorphisms instead of algebraic isomorphisms.
    \begin{definition}
        Let $M$ be a von Neumann algebra with n.s.f.\@ trace $\tau$, and let $H$, $K$ be (left or right) $M$-modules.
        Then a left $M$-module map $\phi:K\to H$ is a \textit{dimension monomorphism} if $\ker\phi$ is zero-dimensional, a \textit{dimension epimorphism} if $\operatorname{coker}\phi$ is zero-dimensional and a \textit{dimension isomorphism} if both conditions are satisfied.
        A sequence $H\stackrel{\phi}{\to} K\stackrel{\psi}{\to} L$ is \textit{dimension exact} if $\ker\psi/\im\phi$ is zero-dimensional.
    \end{definition}
    As the name implies, modules related by a dimension isomorphism in either direction must have the same dimension.

    The following criterion is easily deduced from its well-known counterpart for finite traces \cite[Theorem~2.4]{sauer-discr-groupoids}.
    \begin{proposition}
    \label{thm:zero-dim-criterion}
        Let $M$ be a von Neumann algebra with n.s.f.\@ trace $\tau$, and let $K$ be a locally finite left $M$-module.
        Then $K$ is zero-dimensional if and only if for every $\xi\in K$ and $\epsilon>0$ there exists a projection $p\in M$ such that $\tau(p)<\epsilon$ and $p\xi=\xi$.
    \end{proposition}
    \begin{remark}
        Let $H$ be a locally finite module.
        If $K\subset H$ is a submodule with the property that given $\epsilon>0$ and $\xi\in H$, one can always find a projection $p\in M$ with $\tau(p)<\epsilon$ and $(\I-p)\xi\in K$, we say that $K$ is \textit{rank dense} in $H$.
        By the proposition stated above, this is the same as saying that $H/K$ is zero-dimensional, or that the inclusion map of $K$ in $H$ is a dimension epimorphism.
        The terminology \textit{rank dense} comes from the \textit{rank topology} on $H$ 
        \cite[see e.g.][]{thom-rank-metric}, which is induced by the pseudometric defined by $\|\xi\|_{\mathrm{rank}} = \inf\{\tau(p)\mid p\in M \text{ projection s.t. } p\xi = \xi\}$.
        The statement that $\dim_M H=0$ is equivalent to $\|\xi\|_{\mathrm{rank}}=0$ for all $\xi\in H$.
    \end{remark} 
    The following variant of \cite[Lemma~A.5]{kpv13} remains true in the semifinite setting, with (mutatis mutandis) the same proof.
    \begin{lemma}
    \label{thm:dense-rank-dense-cyclic} 
        Let $M$ be a von Neumann algebra with n.s.f.\@ trace $\tau$, $p\in M$ a projection of finite trace and $K$ a dense $M$-submodule of $L^2(M,\tau)p$.
        Then $K$ is also rank dense in $L^2(M,\tau)p$.
    \end{lemma}
    \begin{definition}
    \label{def:tau-filtration}
        Let $M$ be a von Neumann algebra with n.s.f.\@ trace $\tau$.
        A $\tau$-\textit{filtration} is a family of projections $\mathcal{F}$ in $M$ satisfying
        \begin{itemize}
            \item $\tau(p)<\infty$ for all $p\in\mathcal{F}$;
            \item $p\vee q\in \mathcal{F}$ for all $p,q\in\mathcal{F}$;
            \item $\bigvee_{p\in\mathcal{F}} p = \I$.
        \end{itemize}
        Given a left $M$-submodule $K$ of the $M$-$M$-bimodule $L^2(M)$, we say that $\mathcal{F}$ is \textit{compatible} with $K$ if $Kp\subset K$ for all $p\in\mathcal{F}$.
        The space $K^{\mathcal{F}}=\{\xi\in K\mid \exists p\in\mathcal{F}: \xi p=\xi\}$ is then a left $M$-submodule of $K$.
    \end{definition}
    \begin{lemma}
    \label{thm:dense-rank-dense-filtration} 
        Let $M$ be a von Neumann algebra with n.s.f.\@ trace $\tau$ and $K$ a dense left $M$-submodule of $L^2(M)$.
        For any $\tau$-filtration $\mathcal{F}$ compatible with $K$, we have that $K^{\mathcal{F}}$ is rank dense in $L^2(M)^{\mathcal{F}}$.
    \begin{proof}
        Fix $\xi\in L^2(M)^{\mathcal{F}}$ and $\epsilon>0$.
        Then there exists $p\in\mathcal{F}$ such that $\xi p=\xi$.
        Apply \thref{thm:dense-rank-dense-cyclic} to $Kp\subset L^2(M)p$ to find a projection $q\in M$ such that $(\I-q)\xi\in Kp\subset K^{\mathcal{F}}$ and $\tau(q)<\epsilon$, as required.
    \end{proof}
    \end{lemma}

\subsubsection{Induction and scaling}
    In the sequel, we will frequently work with inclusions of von Neumann algebras equipped with semifinite traces rather than a single algebra.
    In this section, we review some of the techniques that frequently appear in this context.

    The method used to prove the scaling formula for $\lambda$-Markov inclusions discussed in \cite{krvv} appears as Theorem~3.17 in \cite{sauer-thesis}, albeit in the finite setting.
    For completeness, we state the semifinite formulation here.
    \begin{proposition}
    \label{thm:generic-scaling-formula}
    Let $M$ and $N$ be von Neumann algebras equipped with some fixed n.s.f.\@ trace.
    Consider an exact functor $F:\mathrm{Mod}_{-N}\to \mathrm{Mod}_{-M}$.
    Suppose that $F$ preserves filtered colimits and that there exists $\lambda\in (0,+\infty)$ such that
    \[
        \dim_M(F(H)) = \lambda \dim_N(H)
    \]
    for any finitely generated locally finite projective\footnote{Up to isomorphism, these are precisely the modules of the form $P(\CC^n\otimes N)$ for $P\in M_n(\CC)\otimes N$ a projection with finite trace.} right module $H$, then the same formula holds for general locally finite $H$.
    An analogous statement holds for left modules.
    \begin{proof}
        If $H$ is a finitely presented and locally finite $N$-module, there is an exact sequence $0\to K\to L\to H\to 0$ with $K$ finitely generated and $L$ finitely generated locally finite projective.
        Since $N$ is semihereditary, it follows that $K$ is projective as well.  
        Applying the exact functor $F$, we obtain an exact sequence $0\to F(K)\to F(L)\to F(H)\to 0$.
        Since $K,L$ are both finitely generated locally finite projective modules, the additivity property of $\dim_M$ and $\dim_N$ now yields that
        \[
            \dim_M(F(H)) = \dim_M(F(L)) - \dim_M(F(K)) = \lambda\dim_N(L)-\lambda \dim_N(K) = \lambda \dim_N(H)\,.
        \] 
        The passage from finitely presented locally finite modules to arbitrary locally finite modules using the cocontinuity property of $F$ is exactly the same as in \cite{sauer-thesis}.
    \end{proof}
    \end{proposition}

    The following result is a generalisation of \cite[Theorem~1.48]{sauer-thesis} and \cite[Theorem~6.29]{lueck-l2-invariants}, although we do not make any claims about the faithfulness of $-\otimes_N M$ here.
    We follow the proof given in \cite{lueck-l2-invariants}.
    \begin{proposition}
    \label{thm:vna-flatness}
    Let $M$ be a von Neumann algebra with a von Neumann subalgebra $N$.
    Then $M$ is flat as a left (resp.\@ right) $N$-module.
    \begin{proof} 
        Clearly the situation is symmetric, so we only prove the variant for $M$ as a left module.
        Flatness of $M$ is equivalent to the requirement that $\Tor_1^N(H,M)=0$ for all right $N$-modules $H$.
        We will verify this for progressively more general $H$.
        First, assume $H$ to be finitely presented.
        Then there is an exact sequence
        \begin{equation}
        \label{eqn:finitely-presented-exact-sequence}
            0\to K \to \CC^n\otimes N\to H\to 0
        \end{equation}
        with $K$ finitely generated.
        Since von Neumann algebras are semihereditary, $K$ is actually projective, so $K$ is a direct summand of $\CC^m\otimes N$ for some $m$, which means that there exists a (not necessarily self-adjoint) idempotent $E\in M_m(\CC)\otimes N$ such that $K\cong E(\CC^m\otimes N)$.
        Taking $P$ to be the right support projection of $E$, multiplication by $E$ gives an isomorphism $P(\CC^m\otimes N)\to K$, so without loss of generality, we may suppose that \eqref{eqn:finitely-presented-exact-sequence} is given in the form
        \[
            0\to P(\CC^m\otimes N) \stackrel{\theta}{\to} \CC^n\otimes N\to H\to 0\,.
        \]
        In particular, this allows us to consider the injection $\theta$ as an element of $M_{n\times m}(\CC)\otimes N$.
        The first five terms of the long exact sequence of $\Tor$ are given by 
        \[
            0\to \Tor_1^N(H,M)\to P(\CC^m\otimes M) \stackrel{\theta}{\to} \CC^n\otimes M\to H\otimes_N M\to 0
        \]
        because $\Tor_1^N(\CC^n\otimes N,M)=0$ on account of the fact that $\CC^n\otimes N$ is free.
        We keep the notation $\theta$ for the amplification of $\theta:P(\CC^m\otimes N)\to \CC^n\otimes N$ to $P(\CC^m\otimes M)\to \CC^n\otimes M$ since both are given by the same matrix in $M_{n\times m}(\CC)\otimes N\subset M_{n\times m}(\CC)\otimes M$.
        Hence, to show that $\Tor_1^N(H,M)=0$, it suffices to verify that $\theta$ remains injective after applying $-\otimes_N M$.
        However, the injectivity of $\theta: P(\CC^m\otimes N)\to \CC^n\otimes N$ is equivalent to the support projection of $\theta^*\theta$ in $M_m(\CC)\otimes N$ being $P$.
        Viewing $\theta^*\theta$ as an element of $M_m(\CC)\otimes M$ the support projection does not change, so we then find that $\theta:P(\CC^m\otimes M)\to \CC^n\otimes M$ remains injective.
        This concludes the proof for $H$ finitely presented.
        
        The remainder of the proof is exactly the same as \cite{lueck-l2-invariants}.
        Assume $H$ to be finitely generated. Then there exists an exact sequence
        \[
            0\to K\to \CC^n\otimes N\to H\to 0,
        \]
        where we no longer assume $K$ to be finitely generated.
        Let $K_0$ be an arbitrary finitely-generated submodule of $K$.
        Then $\CC^n\otimes N/K_0$ is finitely presented, so
        \[
            \Tor_1^N\left((\CC^n\otimes N/K_0), M\right)=0.
        \]
        Taking the direct limit as $K_0$ increases to $K$, we recover that $\Tor_1^N(H,M)=0$ because $\Tor$ commutes with direct limits.
        To pass from finitely generated $H$ to general $H$, one can similarly take the limit along the finitely generated submodules of $H$ to recover the desired result.
    \end{proof}
    \end{proposition}

    \begin{corollary}
    \label{thm:vna-dimension-preserving}
    Let $M$ be a von Neumann algebra with n.s.f.\@ trace $\tau$, and $N$ a von Neumann subalgebra such that $\tau\vert_N$ remains semifinite.
    On locally finite left (resp.\@ right) $N$-modules, the functor $M\otimes_N-$ (resp.\@ $-\otimes_N M$) has the following properties:
    \begin{itemize}
        \item it is dimension-preserving,
        \item it preserves dimension-exact sequences, and
        \item in particular, it preserves dimension mono-, epi- and isomorphisms.
    \end{itemize}
    \begin{proof}
        The first claim is immediate from \thref{thm:vna-flatness} and \thref{thm:generic-scaling-formula}.
        Note that $M\otimes_N-$ and $-\otimes_N M$ also map locally finite modules to locally finite modules.
        Since any such dimension-preserving exact functor preserves dimension-exact sequences of locally finite modules, the other claims also follow.
    \end{proof}
    \end{corollary}

    \begin{corollary}
    \label{thm:vna-hilbert-dim-iso}
    Let $M$ be a von Neumann algebra with n.s.f.\@ trace $\tau$, and $N$ a von Neumann subalgebra such that $\tau\vert_N$ remains semifinite.
    For any $\tau$-filtration $\mathcal{F}$ of projections in $N$, the multiplication map $\mu:M\otimes_N L^2(N,\tau)^{\mathcal{F}}\to L^2(M,\tau)^{\mathcal{F}}$ is a dimension isomorphism.
    \begin{proof}
        Define $\alpha: M^{\mathcal{F}}\to M\otimes_N N^{\mathcal{F}}$ by sending $m\in M^{\mathcal{F}}$ to $m\otimes_N p$ for $p\in\mathcal{F}$ such that $mp=m$.
        Then $\alpha$ inverts the natural multiplication map $M\otimes_N N^{\mathcal{F}}\to M^{\mathcal{F}}$ and is hence an isomorphism.
        Moreover, note that the embedding of $N^{\mathcal{F}}$ into $L^2(N,\tau)^{\mathcal{F}}$ is an injective dimension isomorphism by \thref{thm:dense-rank-dense-filtration}, so the amplification $\beta: M\otimes_N N^{\mathcal{F}} \to M\otimes_N L^2(N,\tau)^{\mathcal{F}}$ remains an injective dimension isomorphism by \thref{thm:vna-dimension-preserving}.
        Hence, $\beta\circ\alpha:M^{\mathcal{F}}\to M\otimes_N L^2(N,\tau)^{\mathcal{F}}$ is a dimension isomorphism.
        Composing this injection with $\mu$, we recover the embedding of $M^{\mathcal{F}}$ into $L^2(M,\tau)^{\mathcal{F}}$, which is of course also a dimension isomorphism.
        This already shows that the image of $\mu$ is rank dense.
        The fact that the kernel of $\mu$ is zero-dimensional is now straightforward to show using \thref{thm:zero-dim-criterion}.
        Indeed, for any $\xi\in\ker\mu$ and any $\epsilon>0$, there exists a projection $p\in M$ with $\tau(p)<\epsilon$ and $x\in M^{\mathcal{F}}$ such that $p\xi-\xi=(\beta\circ\alpha)(x)$.
        But then $x=(\mu\circ \beta\circ\alpha)(x)=\mu(p\xi-\xi)=0$, so $p\xi=\xi$, as required.
    \end{proof}
    \end{corollary}

\subsection{Rigid \texorpdfstring{$C^*$}{C*}-tensor categories}
\subsubsection{Basic notions and notation}
    A \textit{rigid $C^*$-tensor category} is a $C^*$-tensor category $\mathcal{C}$ where every object $\alpha\in\mathcal{C}$ admits a conjugate $\bar{\alpha}\in\mathcal{C}$.
    The tensor unit of $\mathcal{C}$---which we take to be irreducible---will be denoted by $\epsilon$.
    Rigid $C^*$-tensor categories are always semisimple, i.e.\@ every object decomposes into finitely many irreducible ones.
    The standard reference on rigid $C^*$-tensor categories is \cite[Chapter~2]{neshveyev-tuset}.

    The category $\Rep_f(\mathbb{G})$ of finite-dimensional unitary representations of a compact quantum group $\mathbb{G}$ is a rigid $C^*$-tensor category in a natural way, with the monoidal structure given by the tensor product of representations.
    Objects in $\Rep_f(\mathbb{G})$ are conjugate if and only if they are conjugate as unitary representations in the sense of \cite[Definition~1.4.5]{neshveyev-tuset}.

    Any discrete group $\Gamma$ can be viewed as a rigid $C^*$-tensor category by considering the category $\mathrm{Hilb}_f^\Gamma$ of finite-dimensional $\Gamma$-graded Hilbert spaces \cite[Examples~1.6.3, 2.1.2]{neshveyev-tuset}, with morphisms consisting of linear maps respecting the $\Gamma$-grading.
    Irreducible objects correspond to copies of $\CC$ labelled by elements of $\Gamma$, and the group multiplication induces a monoidal structure on $\mathrm{Hilb}_f^\Gamma$.
    The conjugation operation on $\mathrm{Hilb}_f^\Gamma$ corresponds to inversion in $\Gamma$.

    Throughout, we largely use the same notational conventions as \cite{psv-cohom} and \cite{krvv}.
    The collection of isomorphism classes of irreducible objects of a rigid $C^*$-tensor category $\mathcal{C}$ is denoted by $\Irr(\mathcal{C})$.
    We always assume $\Irr(\mathcal{C})$ to be a countable set.
    Additionally, we always implicitly work with a fixed choice of representatives for each isomorphism class, and we do not distinguish between isomorphism classes and their chosen representatives.

    Given $\alpha,\beta\in\mathcal{C}$, we denote the Banach space of morphisms $\alpha\to\beta$ by $(\beta,\alpha)$.
    The rigidity assumption forces $(\beta,\alpha)$ to be finite-dimensional. 
    In particular, all endomorphism spaces $(\alpha,\alpha)$ are finite-dimensional $C^*$-algebras.
    For all $\alpha\in\mathcal{C}$, the algebra $(\alpha,\alpha)$ comes with a natural trace $\Tr_\alpha: (\alpha,\alpha)\to \CC$ for $\alpha\in\mathcal{C}$ defined in terms of a \textit{standard solution to the conjugate equations} for $\alpha$. 
    These are maps $t_\alpha: \epsilon\to \bar{\alpha}\alpha$, $s_\alpha: \epsilon\to \alpha\bar{\alpha}$ satisfying the following three conditions
    \begin{align*}
        (s_\alpha^*\otimes\I)(\I\otimes t_\alpha)&=\I\,, \\
        (t_\alpha^*\otimes\I)(\I\otimes s_\alpha)&=\I\,, \\
        s_\alpha^*(T\otimes\I) s_\alpha &= t_\alpha^*(\I\otimes T) t_\alpha \qquad \forall T\in (\alpha,\alpha)\,. \numberthis\label{eqn:categorical-trace}
    \end{align*}
    Such standard solutions are unique up to unitary conjugacy.
    Moreover, the map $(\alpha,\alpha)\to (\epsilon,\epsilon)\cong \CC$ defined by the formula \eqref{eqn:categorical-trace} is always tracial and independent of the choice of standard solution.
    This is the \textit{categorical trace} $\Tr_\alpha$, which gives rise to an inner product on $(\beta,\alpha)$ by means of the formula
    \[
        \langle V,W\rangle = \Tr_\beta(VW^*) = \Tr_\alpha(W^*V)\,.
    \]
    Whenever we write $\onb(\beta,\alpha)$, we refer to an orthonormal basis of $(\beta,\alpha)$ with respect to this inner product.
    The positive real number
    \[
        d(\alpha)=\Tr_\alpha(\I)=s_\alpha^*s_\alpha=t_\alpha^*t_\alpha
    \]
    is called the \textit{quantum dimension} of $\alpha$.

    For our purposes, representations of rigid $C^*$-tensor categories are best understood in terms of the \textit{tube algebra}.
    In the case where $\Irr(\mathcal{C})$ is finite (i.e.\@ when $\mathcal{C}$ is a so-called \textit{fusion category}), this construction is due to Ocneanu~\cite{ocneanu-chirality}.
    In the general setting, the tube algebra point of view on the representation theory of rigid $C^*$-tensor categories was developed in \cite{ghosh-jones,psv-cohom}.
    Given a rigid $C^*$-tensor category $\mathcal{C}$, its tube algebra  $\mathcal{A}_{\mathcal{C}}$ has underlying vector space
    \[
         \mathcal{A}_{\mathcal{C}} = \bigoplus_{i,j,\alpha\in\Irr(\mathcal{C})} (i\alpha,\alpha j)\,.
    \]
    By abuse of notation, we usually identify $(i\alpha,\alpha j)$ with the summand indexed by $i,j,\alpha$ in the above direct sum.
    The product and star operation are defined by
    \begin{align*}
        V\cdot W &= \delta_{k,k'} \sum_{\substack{\gamma\in\Irr(\mathcal{C})\\ U\in \onb(\alpha\beta,\gamma)}} d(\gamma) (\I\otimes U^*)(V\otimes\I)(\I\otimes W)(U\otimes \I)\,\\
        V^\# &= (t_\alpha^*\otimes\I^{\otimes 2})(\I\otimes V^*\otimes\I) (\I\otimes\I\otimes s_\alpha)\,.
    \end{align*}
    where $V\in (i\alpha,\alpha k)$ and $W\in (k'\beta,\beta j)$.
    One checks that these definitions do not depend on any non-canonical choices, and that $\mathcal{A}_{\mathcal{C}}$ becomes a $\ast$-algebra with these operations.
    In order to eliminate any potential confusion with composition and adjoints of morphisms in $\mathcal{C}$, we explicitly denote the tube algebra operations by $\cdot$ and $\#$.

    For $i\in\Irr(\mathcal{C})$, let $p_i$ be the identity on $i$ considered as an element of $(i\epsilon,\epsilon i)\subset\mathcal{A}_{\mathcal{C}}$.
    Then $(p_i)_{i\in\Irr(\mathcal{C})}$ is a family of self-adjoint idempotents satisfying $p_i\cdot V \cdot p_j=V$ for all $V\in (i\alpha,\alpha j)$.
    Given a finite subset $F\subset\Irr(\mathcal{C})$, we denote the sum $\sum_{i\in F} p_i$ by $p_F$.
    Hence, while $\mathcal{A}_{\mathcal{C}}$ is typically not a unital algebra, every element $V\in\mathcal{A}_{\mathcal{C}}$ is supported under some $p_F$.

    The tube algebra comes with two natural linear functionals:
    \begin{align*}
        &\counit: \mathcal{A}_{\mathcal{C}} \to \CC: 
        &&V\in (i\alpha,\alpha j)\mapsto \begin{cases}
            \Tr_\alpha(V) & i=j=\epsilon\\
            0 & \text{otherwise}
        \end{cases}\,,\\
        &\tau: \mathcal{A}_{\mathcal{C}} \to \CC:
        &&V\in (i\alpha,\alpha j)\mapsto \begin{cases}
            \Tr_i(V) & i=j, \alpha=\epsilon\\
            0 & \text{otherwise}
        \end{cases}\,.
    \end{align*}
    The functional $\counit$ is in fact a $\ast$-homomorphism, referred to as the \textit{counit} or \textit{trivial representation} of $\mathcal{A}_{\mathcal{C}}$.
    This allows us to consider $\CC$ both as a left and right $\mathcal{A}_{\mathcal{C}}$-module.
    On the other hand, the functional $\tau$ is a faithful trace on $\mathcal{A}_{\mathcal{C}}$, and positive in the sense that $\tau(V^\#\cdot V)\geq 0$.
    The inner product it induces on $\mathcal{A}_{\mathcal{C}}$ satisfies
    \[
        \tau(V^\# \cdot W) = \begin{cases}
            d(\alpha)^{-1} \Tr_{\alpha j}(V^*W) & \alpha=\beta\\
            0 & \text{otherwise}
        \end{cases}
    \]
    where $V\in (i\alpha,\alpha j)$ and $W\in (i\beta,\beta j)$ \cite[see][Proposition~3.10]{psv-cohom}.
    Completing $\mathcal{A}_{\mathcal{C}}$ with respect to this inner product yields a GNS space $L^2(\mathcal{A}_{\mathcal{C}},\tau)$ on which $\mathcal{A}_{\mathcal{C}}$ acts faithfully by left (and right) multiplication \cite[Proposition~3.10]{psv-cohom}.
    If we now take $\mathcal{M}_{\mathcal{C}}$ to be the weak closure of $\mathcal{A}_{\mathcal{C}}$ embedded in $\mathcal{B}(L^2(\mathcal{A}_{\mathcal{C}},\tau))$ as left multiplication operators, 
    then $\tau$ extends uniquely to a n.s.f.\@ tracial weight on $\mathcal{M}_{\mathcal{C}}$.

    Following \cite{psv-cohom}, we put
    \[
        H_n(\mathcal{C}; V) = \Tor^{\mathcal{A}_{\mathcal{C}}}_n(V, \CC)
    \]
    for any nondegenerate right $\mathcal{A}_{\mathcal{C}}$-module $V$, where $\CC$ is considered as a left $\mathcal{A}_{\mathcal{C}}$-module.
    This generalises the definition of homology for actions of discrete groups on vector spaces.
    
    The Hilbert space $L^2(\mathcal{A}_{\mathcal{C}},\tau)$ is naturally an $\mathcal{M}_{\mathcal{C}}$-$\mathcal{M}_{\mathcal{C}}$-bimodule, but since we want to view it rather as an $\mathcal{A}_{\mathcal{C}}$-module from the right, it is more convenient to consider the subspace given by the following algebraic direct sum:
    \[
        L^2(\mathcal{A}_{\mathcal{C}},\tau)^0 = \bigoplus_{i\in\Irr(\mathcal{C})} L^2(\mathcal{A}_{\mathcal{C}},\tau)\cdot p_i\,.
    \] 
    Observe that this subspace remains an $\mathcal{M}_{\mathcal{C}}$-$\mathcal{A}_{\mathcal{C}}$-bimodule, and that $L^2(\mathcal{A}_{\mathcal{C}},\tau)^0$ is locally finite as a left $\mathcal{M}_{\mathcal{C}}$-module.
    We can analogously consider the $\mathcal{A}_{\mathcal{C}}$-$\mathcal{M}_{\mathcal{C}}$-bimodule $\lsup{0}{L^2(\mathcal{A}_{\mathcal{C}},\tau)}$ defined by the algebraic direct sum
    \[
        \lsup{0}{L^2(\mathcal{A}_{\mathcal{C}},\tau)} = \bigoplus_{i\in\Irr(\mathcal{C})} p_i\cdot L^2(\mathcal{A}_{\mathcal{C}},\tau)\,.
    \]
    Of course $L^2(\mathcal{A}_{\mathcal{C}},\tau)^0$ is precisely $L^2(\mathcal{A}_{\mathcal{C}},\tau)^{\mathcal{F}}$ as in \thref{def:tau-filtration}, where $\mathcal{F}$ is the filtration given by all projections $p_F$ for $F$ ranging over the finite subsets of $\Irr(\mathcal{C})$.  
    One can then introduce the $L^2$-Betti numbers of $\mathcal{C}$ via
    \[
        \beta^{(2)}_n(\mathcal{C}) = \dim_{\mathcal{M}_{\mathcal{C}}} H_n(\mathcal{C}; L^2(\mathcal{A}_{\mathcal{C}},\tau)^0)
    \]
    since the left $\mathcal{M}_{\mathcal{C}}$-module structure on $L^2(\mathcal{A}_{\mathcal{C}},\tau)^0$ naturally passes to the homology spaces.
    As explained in \cite[Proposition~6.4, Corollary~7.2]{psv-cohom}, the analogous cohomological definition using $\Ext_{\mathcal{A}_{\mathcal{C}}}^\bullet(\CC; -)$ is equivalent.

\subsubsection{Inclusions of rigid \texorpdfstring{$C^*$}{C*}-tensor categories}
\label{sec:cat-inclusions}
    Let $\mathcal{D}$ be a rigid $C^*$-tensor category and $\mathcal{C}$ a full $C^*$-tensor subcategory.
    The fullness guarantees that $\Irr(\mathcal{C})\subset\Irr(\mathcal{D})$.
    For $\alpha\in\Irr(\mathcal{D})$, put
    \[
        \alpha\mathcal{C} = \{\beta\in\Irr(\mathcal{D})\mid \exists\pi\in\mathcal{C}: \beta\hookrightarrow \alpha\pi\}\,,\qquad
        \mathcal{C}\alpha = \{\beta\in\Irr(\mathcal{D})\mid \exists\pi\in\mathcal{C}: \beta\hookrightarrow \pi\alpha\}\,.
    \]
    In this way, one obtains partitions of $\Irr(\mathcal{D})$ into left/right $\mathcal{C}$\textit{-cosets}.
    The set of left (resp.\@ right) cosets is denoted by $\mathcal{D}/\mathcal{C}$ (resp.\@ $\mathcal{D}\backslash\mathcal{C}$).

    For any object $\alpha\in\mathcal{D}$, the rigidity assumption guarantees that there exists a largest subobject of $\alpha$ that is contained in $\mathcal{C}$, which we denote by~$[\alpha]_{\mathcal{C}}$.

    Given some family of objects $S\subset\mathcal{D}$, denote 
    \[
        \tilde{\mathcal{A}}_{S} = \bigoplus_{i,j\in\Irr(\mathcal{D})} \bigoplus_{\alpha\in\Irr(\mathcal{D})\cap S} (i\alpha,\alpha j) \subset\mathcal{A}_{\mathcal{D}}\,.
    \]
    Note that the ambient category is suppressed in the notation.
    We denote by $e_S$ the orthogonal projection of $L^2(\mathcal{A}_{\mathcal{D}},\tau)$ onto the closure of $\tilde{\mathcal{A}}_{S}$.
    When $S=\{e\}$, $\tilde{\mathcal{A}}_S$ is the commutative $*$-algebra
    \begin{equation}
    \label{eqn:b-definition}
        \mathcal{B} = \mathrm{span}\{p_i\mid i\in\Irr(\mathcal{D})\}\,.
    \end{equation}

    In homological arguments involving the tube algebra, working with the inclusions $\mathcal{A}_{\mathcal{C}}\subset\mathcal{A}_{\mathcal{D}}$ and $\mathcal{M}_{\mathcal{C}}\subset\mathcal{M}_{\mathcal{D}}$ is not always practical, the main reason being that $\mathcal{A}_{\mathcal{D}}$ is a degenerate $\mathcal{A}_{\mathcal{C}}$-module.
    Instead, it is often better to consider the nondegenerate ring extension 
    \[
        \tilde{\mathcal{A}}_{\mathcal{C}} = \bigoplus_{i,j\in\Irr(\mathcal{D})} \bigoplus_{\alpha\in\Irr(\mathcal{C})} (i\alpha,\alpha j) \subset\mathcal{A}_{\mathcal{D}}\,.
    \]
    Taking the weak closure of the left regular representation of $\tilde{\mathcal{A}}_{\mathcal{C}}$ on $L^2(\tilde{\mathcal{A}}_{\mathcal{C}},\tau)$, we obtain a von Neumann algebra $\tilde{\mathcal{M}}_{\mathcal{C}}$.
    Using \cite[Lemma~3.11]{krvv}, one checks that $\tilde{\mathcal{M}}_{\mathcal{C}}$ is canonically isomorphic to the von Neumann algebra subalgebra of $\mathcal{M}_{\mathcal{D}}$ generated by $\tilde{\mathcal{A}}_{\mathcal{C}}$ acting on $L^2(\mathcal{A}_{\mathcal{D}},\tau)$ by multiplication, and that the projection $e_\mathcal{C}$ commutes with both left- and right representations of $\tilde{\mathcal{M}}_{\mathcal{C}}$ on $L^2(\mathcal{A}_{\mathcal{D}},\tau)$.
    This identification also respects the canonical tracial weights on $\tilde{\mathcal{M}}_{\mathcal{C}}$ and $\mathcal{M}_{\mathcal{D}}$.
    We therefore make no further distinction between $\tilde{\mathcal{M}}_{\mathcal{C}}$ and its image inside $\mathcal{M}_{\mathcal{D}}$.

    The restriction of $\tau$ to $\tilde{\mathcal{M}}_{\mathcal{C}}$ remains semifinite, and $e_{\mathcal{C}}$ induces a conditional expectation from $\mathcal{M}_{\mathcal{D}}$ onto $\tilde{\mathcal{M}}_{\mathcal{C}}$.
    As a consequence of \cite[Lemma~3.10]{krvv}, we have that
    \[
        \beta^{(2)}_n(\mathcal{C})=
        \dim_{\tilde{\mathcal{M}}_{\mathcal{C}}}\Tor^{\tilde{\mathcal{A}}_{\mathcal{C}}}_\bullet(L^2(\tilde{\mathcal{A}}_{\mathcal{C}},\tau)^0, \CC) \,.
    \]
    In other words, we can compute the $L^2$-Betti numbers of $\mathcal{C}$ using $\tilde{\mathcal{A}}_{\mathcal{C}}$ instead of $\mathcal{A}_{\mathcal{C}}$.
    Additionally, \cite[Proposition~3.12]{krvv} shows\footnote{The result in \cite{krvv} is stated for inclusions of finite index, but the projectivity claim in the statement applies to arbitrary inclusions with exactly the same proof.} that $\mathcal{A}_{\mathcal{D}}$ is projective over $\tilde{\mathcal{A}}_{\mathcal{C}}$.
    If $\mathcal{C}\subset\mathcal{D}$ is an inclusion of finite index (in the sense of \cite{krvv}, i.e.\@ $\#\mathcal{D}/\mathcal{C}<\infty$), the same holds for $\tilde{\mathcal{M}}_{\mathcal{C}}\subset\mathcal{M}_{\mathcal{D}}$, but this is no longer true in general.
    However, we still get that $\mathcal{M}_{\mathcal{D}}$ is a flat $\tilde{\mathcal{M}}_{\mathcal{C}}$-module as a consequence of \thref{thm:vna-flatness}.

    Put $\mathfrak{n}_\tau=\{x\in\mathcal{M}_{\mathcal{D}}\mid \tau(x^*x)<\infty\}$.
    For any family of objects $S\subset\mathcal{D}$, we define 
    \begin{equation}
    \label{eqn:z-module-def}
        \mathcal{Z}_S = e_S(\mathfrak{n}_\tau) \subset L^2(\mathcal{A}_{\mathcal{D}}, \tau)\,.
    \end{equation}
    When $S=\mathcal{C}$, $\mathcal{Z}_{\mathcal{C}}$ is equal to $\tilde{\mathcal{M}}_{\mathcal{C}} \cap \mathfrak{n}_\tau$.
    In general, it is not necessarily the case that $\mathcal{Z}_S\subset\mathcal{M}_{\mathcal{D}}$.
    However, when $S$ is a left or right coset, $\mathcal{Z}_S$ has the following properties.
    \begin{lemma}
    \label{thm:z-orbit-module}
    Let $\mathcal{D}$ be a rigid $C^*$-tensor category and $\mathcal{C}$ a full tensor subcategory.
    Fix $\alpha\in\Irr(\mathcal{D})$. 
    Then
    \begin{enumerate}[(i)]
    \item the projection $e_{\alpha\mathcal{C}}$ commutes with the right action of $\tilde{\mathcal{M}}_{\mathcal{C}}$ on $L^2(\mathcal{A}_{\mathcal{D}},\tau)$\,;
    \item $\mathcal{Z}_{\alpha\mathcal{C}}$ is a flat right $\tilde{\mathcal{M}}_{\mathcal{C}}$-submodule of $L^2(\mathcal{A}_{\mathcal{D}},\tau)$\,;
    \item for all finite subsets $F\subset\Irr(\mathcal{D})$, $p_F\cdot\mathcal{Z}_{\alpha\mathcal{C}}$ is a subset of $\mathcal{M}_{\mathcal{D}}$ and projective as a right $\tilde{\mathcal{M}}_{\mathcal{C}}$-module.
    \end{enumerate}
    Analogous results hold if one considers $\mathcal{C}\alpha$ instead of $\alpha\mathcal{C}$.
    \begin{proof}
        Recall from \cite[Lemma~3.11]{krvv} that 
        \begin{equation}
        \label{eqn:e-proj-expansion}
            p_i\cdot e_{\alpha\mathcal{C}}(\xi) = \frac{d(\alpha)}{d([\bar{\alpha}\alpha]_{\mathcal{C}})}\sum_{\substack{j\in\Irr(\mathcal{D})\\ W\in \onb(i\alpha,\alpha j)}} d(j) W\cdot e_{\mathcal{C}} (W^\#\cdot \xi)
        \end{equation}
        for all $i\in\Irr(\mathcal{D})$ and $\xi\in L^2(\mathcal{A}_{\mathcal{D}},\tau)$.
        Since the $(p_i)_{i\in\Irr(\mathcal{D})}$ sum up to the identity on $L^2(\mathcal{A}_{\mathcal{D}},\tau)$, this already implies (i).

        Choose $x\in\mathfrak{n}_\tau$.
        For every $i\in\Irr(\mathcal{D})$, substituting $\xi=x$ in \eqref{eqn:e-proj-expansion} produces an element in $\mathfrak{n}_\tau$, since the sum is finite and $e_{\mathcal{C}}(W^\#\cdot x) \in \tilde{\mathcal{M}}_{\mathcal{C}}\cap \mathfrak{n}_\tau$ for all $W\in\mathcal{A}_{\mathcal{D}}$.
        Hence, $p_F\cdot\mathcal{Z}_{\alpha\mathcal{C}}\subset\mathfrak{n}_\tau\subset\mathcal{M}_{\mathcal{D}}$ for all finite subsets $F\subset\Irr(\mathcal{D})$.

        The fact that $\mathcal{Z}_{\alpha\mathcal{C}}$ is a right $\tilde{\mathcal{M}}_{\mathcal{C}}$-module is an immediate consequence of~(i).
        For all $i\in\Irr(\mathcal{D})$, consider the right $\tilde{\mathcal{M}}_{\mathcal{C}}$-linear map
        \begin{equation}
        \label{eqn:z-projective-targ}
            \Phi_i: \mathcal{Z}_{\alpha\mathcal{C}} \to \bigoplus_{j\in\Irr(\mathcal{D})} (i\alpha,\alpha j)\otimes p_j\cdot \tilde{\mathcal{M}}_{\mathcal{C}}:
            \xi \mapsto \bigoplus_{j\in\Irr(\mathcal{D})} \left(\sum_{W\in (i\alpha,\alpha j)} d(j)W\otimes e_{\mathcal{C}}(W^\#\cdot \xi)\right)\,.
        \end{equation}
        Note that only finitely many terms in the direct sum are nonzero, so $\Phi_i$ is well-defined.
        Also note that the module appearing on the right side of \eqref{eqn:z-projective-targ} is projective.
        By \eqref{eqn:e-proj-expansion}, the map 
        \[
            \Psi_i: \bigoplus_{j\in\Irr(\mathcal{D})} (i\alpha,\alpha j)\otimes p_j\cdot \tilde{\mathcal{M}}_{\mathcal{C}}\to\mathcal{Z}_{\alpha\mathcal{C}}:
            \bigoplus_{j\in\Irr(\mathcal{D})} V_j\otimes W_j \mapsto \sum_{j\in\Irr(\mathcal{D})} V_j\cdot W_j
        \]
        has the property that $(\Psi_i\circ\Phi_i)(\xi)$ equals $p_i\cdot \xi$ (up to a positive scaling constant), so $\Phi_i$ restricts to an embedding of $p_i\cdot\mathcal{Z}_{\alpha\mathcal{C}}$ into a projective module as a direct summand.
        This implies that $p_i\cdot\mathcal{Z}_{\alpha\mathcal{C}}$ is itself projective.
        Taking direct sums, we get that $p_F\cdot\mathcal{Z}_{\alpha\mathcal{C}}$ is a projective $\tilde{\mathcal{M}}_{\mathcal{C}}$-submodule of $\mathcal{M}_{\mathcal{D}}$ for all $F\subset\Irr(\mathcal{D})$.

        The maps $\Phi_i$ additionally induce an injection
        \[
            \mathcal{Z}_{\alpha\mathcal{C}} \to \prod_{i\in\Irr(\mathcal{D})} \left(\bigoplus_{j\in\Irr(\mathcal{D})} (i\alpha,\alpha j)\otimes p_j\cdot \tilde{\mathcal{M}}_{\mathcal{C}}\right)\,.
        \]
        This realises $\mathcal{Z}_{\alpha\mathcal{C}}$ as a submodule of a direct product of flat $\tilde{\mathcal{M}}_{\mathcal{C}}$-modules.
        By \cite[Theorems~4.47, 4.67]{lam} and the fact that von Neumann algebras are semihereditary rings, this implies that $\mathcal{Z}_{\alpha\mathcal{C}}$ is flat, as claimed.
    \end{proof}
    \end{lemma}

\subsection{Quasi-regular inclusions}
\label{sec:quasireg-prelim}

The theory of $L^2$-Betti numbers for rigid $C^*$-tensor categories introduced in \cite{psv-cohom} can be understood in a more broad framework, where the objects of interest are \textit{quasi-regular} inclusions of type $\mathrm{II}_1$ factors.
\begin{definition}
    Let $T\subset S$ be an inclusion of type $\mathrm{II}_1$ factors.
    Define the \textit{quasi-normaliser} of $T$ inside $S$ by
    \begin{align*}
        \mathrm{QN}_S(T) &= \Big\{ x\in S\,\Big\vert\, \exists a_1,\ldots,a_n,b_1,\ldots,b_m\in S: xT\subset \sum_{i=1}^n T a_i, Tx\subset\sum_{j=1}^m b_j T\Big\}\,.
    \end{align*}
    We say that $T\subset S$ is \textit{quasi-regular} if $S=\mathrm{QN}_S(T)''$.
    In what follows, we always consider \textit{irreducible} quasi-regular inclusions, i.e.\@ $T'\cap S=\CC\I$.
\end{definition}
Proceeding as in \cite{psv-cohom}, one can define a tube algebra for any irreducible quasi-regular inclusion $T\subset S$ together with some choice of category of Hilbert $T$-$T$-bimodules $\mathcal{C}$.
In much the same way as for rigid $C^*$-tensor categories, this gives rise to a representation theory and a purely algebraic Hochschild-type (co)homology theory \cite[\S~7]{psv-cohom}.
If $T\subset S$ is additionally \textit{unimodular}---referring to the requirement that all $T$-$T$-subbimodules of $L^2(S)^{\cfpow{T}{n}}$, $n\geq 1$ have equal left- and right $T$-dimension---the tube algebra embeds densely into a von Neumann algebra with a canonical n.s.f.\@ trace, which allows for a definition of $L^2$-Betti numbers \cite[Definition~4.3]{psv-cohom}.

By realising a finitely generated rigid $C^*$-tensor category $\mathcal{C}$ as the category of $M$-$M$-bimodules appearing in the Jones tower of an extremal finite-index subfactor $N\subset M$, one recovers a unimodular quasi-regular inclusion $T\subset S$ by passing to the symmetric enveloping inclusion of $N\subset M$.
Defined with respect to an appropriate category of $T$-$T$-bimodules, the tube algebra of $T\subset S$ is strongly Morita equivalent to $\mathcal{A}_{\mathcal{C}}$ \cite[Proposition~3.12]{psv-cohom}.
In particular, this gives rise to an identification of the (co)homology theories of $\mathcal{C}$ and $T\subset S$.

While the tube algebra perspective is very useful to frame the representation theory of rigid $C^*$-tensor categories in this setting, we will take a slightly different (but ultimately equivalent) point of view.
For our purposes, it will be convenient to define cohomology spaces for $T\subset S$ with coefficients in any Hilbert $S$-$S$-bimodule, following \cite[\S~4]{psv-cohom}.
We recall the necessary definitions and notation below.
All computations at the level of quasi-regular inclusions will be cohomological in nature, so we omit the analogous homological definitions.
\begin{remark}
    Given a von Neumann algebra $T$, a right $T$-module $V$ and a left $T$-module $W$, the notation $V\otimes_T W$ will always refer to the algebraic relative tensor product of $V$ and $W$.
    If $T$ is a type $\mathrm{II}_1$ factor and $V,W$ are Hilbert $T$-modules, we denote the Connes tensor product of Hilbert $T$-modules by $V\cftimes{T} W$ to emphasise the difference.
    We use a similar notational convention for tensor powers.
    In the same vein, we use the notation $\BHom(-,-)$ and $\BEnd(-)$ to denote bounded morphisms and endomorphisms, as opposed to purely algebraic ones.
\end{remark}
Given a quasi-regular inclusion $T\subset S$, put $\mathcal{S}=\mathrm{QN}_S(T)$.
For any Hilbert $S$-$S$-bimodule $\mathcal{K}$ there is a Hochschild-type \cite{hochschild-rel-homalg} bar cochain complex with terms
\[
    C^n=\Hom_{T-T}(\mathcal{S}^{\otimes_T^n}, \mathcal{K})
\]
for $n\geq 0$, where we adopt the convention that $\mathcal{S}^{\otimes_T^0}=T$.
In other words, $C^0=\{\xi\in\mathcal{K}\mid \xi \text{ is } T\text{-central}\}$.
The differentials of $C^\bullet$ are given by
\begin{align*}
    &\partial:C^n\to C^{n+1}: \sum_{i=0}^{n+1} (-1)^i\partial_i\,,\\
    &(\partial_0 \phi)(s_0\otimes_T\cdots \otimes_T s_n) = s_0\phi(s_1\otimes_T\cdots \otimes_T s_n)\,,\\
    &(\partial_i \phi)(s_0\otimes_T\cdots \otimes_T s_n) = \phi(s_0\otimes_T\cdots \otimes_T s_{i-1}s_i\otimes_T\cdots\otimes_T s_n)\qquad i\in\{1,\ldots,n\}\,,\\
    &(\partial_{n+1} \phi)(s_0\otimes_T\cdots \otimes_T s_n) = \phi(s_0\otimes_T\cdots \otimes_T s_{n-1})s_n\,.
\end{align*}
We denote the cohomology of this complex by $H^\bullet(T\subset\mathcal{S}; \mathcal{K})$.
If $T\subset S$ is unimodular, the bounded $S$-$S$-bimodule endomorphisms of the \textit{regular} Hilbert $S$-$S$-bimodule given by
\begin{equation}
\label{eqn:reg-bimodule-def}
    \mathcal{H}_{\mathrm{reg}} = \bigoplus_{n\geq 2} L^2(S)^{\cfpow{T}{n}}
\end{equation}
form a von Neumann algebra $\mathcal{M}(T\subset S)$ that comes with a natural n.s.f.\@ trace \cite[see][\S~4]{psv-cohom}.
Hence, one gets $L^2$-Betti numbers by putting
\[
    \beta_n^{(2)}(T\subset\mathcal{S}) = \dim_{\mathcal{M}(T\subset S)} H^n(T\subset\mathcal{S}; \mathcal{H}_{\mathrm{reg}})\,.
\]
Given an extremal finite-index subfactor $N\subset M$, the SE-inclusion $T\subset S$ will be a unimodular quasi-regular inclusion of type $\mathrm{II}_1$ factors.
On the other hand, $N\subset M$ also comes with a natural rigid $C^*$-tensor category $\mathcal{C}(N\subset M)$, generated by the $M$-$M$-bimodules appearing in the Jones tower of $N\subset M$.
As explained in \cite[\S~7.2]{psv-cohom}, the $L^2$-Betti numbers of $T\subset S$ and $\mathcal{C}(N\subset M)$ agree.
In particular, one recovers the $L^2$-Betti numbers of a discrete group $\Gamma$ by choosing an outer action of $\Gamma$ on some type $\mathrm{II}_1$ factor $P$ and taking the $L^2$-Betti numbers of $P\subset P\rtimes\Gamma$.

\subsection{Hecke pairs and Schlichting completions}
\label{sec:hecke-pairs-prelim}
Throughout this section, $\Gamma$ will denote a discrete group and $\Lambda$ an almost-normal subgroup, i.e.\@ a subgroup of $\Gamma$ with the property that $[\Lambda:\Lambda\cap g\Lambda g^{-1}]<\infty$ for all $g\in\Gamma$.
All discrete groups we consider will be countable.
Inclusions $\Lambda<\Gamma$ of this form are referred to as \textit{Hecke pairs}.
Equivalently, one could require that the double coset $\Lambda g\Lambda$ be a finite subset of $\Gamma/\Lambda$ (or $\Lambda\backslash\Gamma$) for all $g\in\Gamma$.
In order to study Hecke pairs systematically, it is often useful to instead consider a (typically nondiscrete) completion of $\Gamma$, originally due to Schlichting \cite{schlichting-peri-stab} and further developed by Tzanev \cite{tzanev}.
Since we will make extensive use of this construction in the final part of this article, we recall some of the essentials here.
The reader may consult \cite{tzanev,klq} for more detailed exposition.

Letting $\Gamma$ act on $\Gamma/\Lambda$ and completing the resulting permutation group in the topology of pointwise convergence, one obtains a locally compact totally disconnected group $G$ acting on $\Gamma/\Lambda$ by permutations with compact stabilisers---for a more precise statement, see \cite[\S~3]{klq}.
In particular, the stabiliser $K$ of $\Lambda$ is a compact open subgroup, in which the image of $\Lambda$ is dense.
The pair $K<G$ is called the \textit{Schlichting completion} of $\Lambda <\Gamma$, and we have a natural identification of the coset spaces $\Gamma/\Lambda$ and $G/K$.
Note that the natural map $\Gamma\to G$ need not be injective in general, e.g.\@ when $\Lambda<\Gamma$ is a proper normal subgroup\footnote{In this degenerate case, the Schlichting completion is canonically isomorphic to the quotient group $\Gamma/\Lambda$.}.
In spite of this, we usually do not distinguish between elements of $\Gamma$ and their images in $G$.

The locally compact group $G$ can fail to be unimodular. In general, $G$ is unimodular precisely when $[\Lambda:\Lambda\cap x\Lambda x^{-1}]=[\Lambda:\Lambda\cap x^{-1}\Lambda x]$ for all $x\in \Gamma$.
In this case, we will call the Hecke pair $\Lambda<\Gamma$ unimodular as well.

\begin{example}
Classical examples of Hecke pairs include the following:
\begin{itemize}
    \item $\mathrm{SL}_n(\ZZ)\subset\mathrm{SL}_n(\ZZ[1/p])$ for any prime $p$, with Schlichting completion $\mathrm{PSL}_n(\ZZ_p)\subset\mathrm{PSL}_n(\QQ_p)$ \cite[Example~3.10]{shalom-willis}.
    \item 
    The Baumslag--Solitar group $BS(m, n) = \langle a, t\mid t^{-1} a^m t = a^n\rangle$ contains $\langle a\rangle\cong \ZZ$ as a natural almost-normal subgroup.
\end{itemize}
The first example is unimodular. 
The Hecke pairs arising from Baumslag--Solitar groups typically are not unimodular, unless $|m|=n$ (see e.g.\@ \cite[Lemma~9.1]{raum-lc-powers}).
\end{example}

We conclude this section by taking a look at the group von Neumann algebra $L(G)$ of $G$ when $G$ is unimodular.
Recall that $L(G)$ is the von Neumann algebra defined by completing $C_c(G)$ acting on $L^2(G)$ by convolution.
By unimodularity, the von Neumann algebra $L(G)$ comes equipped with a natural n.s.f.\@ Plancherel trace $\tau$---given $f\in C_c(G)$, $\tau$ assigns the value $f(e)$ to the element of $L(G)$ given by convolution with $f$.
Since the convolution representation of $C_c(G)$ depends on the normalisation of the Haar measure, so does the Plancherel trace.
For our purposes, however, the normalisation will always be fixed to assign measure $1$ to $K$, so we do not make the dependency explicit in the notation.

For $L\subset K$ a compact open subgroup, the convolution operator associated with the normalised characteristic function $\mu(L)^{-1}\chi_L\in C_c(G)$ defines a projection $p_L\in L(G)$, given by
\begin{equation}
\label{eqn:subgroup-invar-projection}
    (p_L\xi)(h) = \frac{1}{\mu(L)} \int_L \xi(g^{-1}h) \diff \mu(g)\,.
\end{equation}
As an operator on $L^2(G)$, $p_L$ is the projection onto the left-$L$-invariant vectors, so $p_L L^2(G)\cong \ell^2(L\backslash G)$.
In fact, $p_L L(G) p_L$ is exactly the commutant of the right $G$-action on $\ell^2(L\backslash G)$, and the restriction of the Plancherel trace to this corner of $L(G)$ is implemented by the vector state associated with the (non-normalised) characteristic function of $L$.
If $L=K$ we also have that $\ell^2(K\backslash G)\cong \ell^2(\Lambda\backslash\Gamma)$, which allows us to view $p_K L^2(G) p_K$ as acting on $\ell^2(\Lambda\backslash\Gamma)$.
The \textit{Hecke algebra} of $\Lambda<\Gamma$---i.e.\@ the convolution $*$-algebra spanned by the characteristic functions $\chi_{\Lambda x\Lambda}$---acts naturally on $\ell^2(\Lambda\backslash\Gamma)$ by convolution, and forms a weakly dense subalgebra of $p_KL(G)p_K$ in this way.

\subsection{Continuous cohomology and \texorpdfstring{$L^2$}{L2}-Betti numbers for locally compact totally disconnected groups}
\label{sec:locally-profinite-cohom}

In order to formulate a useful cohomology theory for locally compact groups, one needs to take the topology into account in various ways.
Consider a (second countable) locally compact group $G$.
Instead of arbitrary $G$-modules, the coefficient spaces are limited to (quasi-complete) locally convex topological vector spaces equipped with a continuous right $G$-action.
Given such a space $V$, one can then define $H^n_c(G;V)$ in terms of a bar resolution involving the spaces of continuous functions $C(G^n, V)$ for $n\geq 1$ equipped with a diagonal translation action (see e.g.\@ \cite[Proposition~I.1.2]{guichardet}).
In particular, there is a meaningful way to talk about cohomology with coefficients in a Hilbert space on which $G$ acts continuously from the right.
If $G$ is unimodular, this naturally leads to a definition of $L^2$-Betti numbers for $G$ \cite{hdp-thesis}.
We recall some of the exposition of \cite{guichardet} and \cite{hdp-thesis} here.

Take $G$ to be a unimodular locally compact group and fix a Haar measure on $G$. 
Then the left $L(G)$-module structure on $L^2(G)$ turns $H^n_c(G;L^2(G))$ into an $L(G)$-module, 
Following \cite[\S~3.1]{hdp-thesis}, this allows one to define
\begin{equation}
\label{eqn:locally-profinite-l2betti}
    \beta^{(2)}_n(G) = \dim_{L(G)} H^n_c(G; L^2(G))\,.
\end{equation}
Note that $L^2$-Betti numbers are only defined up to a constant multiplicative factor, since scaling the Haar measure also scales the $L^2$-Betti numbers.
If $G$ is discrete and we choose the Haar measure such that each group element is assigned measure 1, we recover the usual $L^2$-Betti numbers for discrete groups.
See \cite{kpv13,hdp-thesis} for a more in-depth treatment.

For general $G$, the right-hand side of \eqref{eqn:locally-profinite-l2betti} is tightly coupled with the topology of $G$ in various ways.
However, when one restricts attention to totally disconnected locally compact groups together with a choice of compact open subgroup, almost all topological subtleties disappear; this is also the point of view espoused in \cite[Chapter~5]{hdp-thesis}.
Since this algebraic description will be useful for us later, we sketch the construction here.
Our principal reference is \cite{guichardet}.

Consider a locally compact totally disconnected group $G$ with $K$ a compact open subgroup, and suppose that we are given some right $\CC[G]$-module $V$.
For all $n\geq 0$ put $C^n=\Fun((K\backslash G)^{n+1}, V)$, i.e.\@ the space of all functions from $G^{n+1}$ to $V$ that are left-$K$-invariant in each coordinate.
We let $G$ act on $C^n$ as follows:
\[
    (f\cdot g)(g_0,\ldots, g_n) = f(g_0g^{-1},\ldots,g_ng^{-1})\cdot g\,.
\]
Then there is a complex $0\to V\to C^0\to \cdots\to C^n\to\cdots$ of right $\CC[G]$-modules with differentials given by 
\begin{align*}
    (\partial f)(g_0, \ldots, g_{n+1}) = \sum_{i=0}^{n+1} (-1)^i f(g_0,\ldots, \hat{g}_i,\ldots,g_{n+1})
\end{align*}
and the augmentation $V\to C^0$ defined by sending a vector $v\in V$ to the constant function with value $v$.
A routine computation shows that this complex is exact.

If $V$ is a quasi-complete locally convex topological vector space equipped with a continuous right $G$-action, this a priori purely algebraic complex turns out to be an injective resolution of $V$ in the appropriate sense \cite[\S~I.2, Proposition~III.2.3]{guichardet}. 
This allows us to compute $H^\bullet_c(G; V)$ as the cohomology of the complex
\begin{equation}
\label{eqn:cont-cohom-algebr}
    0\to \Fun(K\backslash G, V)^G\to \cdots \to \Fun((K\backslash G)^{n+1}, V)^G\to \cdots
\end{equation}
where $W^G$ is the space of $G$-invariant vectors in the $\CC[G]$-module $W$.

In addition to the previous assumptions, we now take $G$ to be unimodular, and we choose the normalisation of the Haar measure $\mu$ on $G$ such that $\mu(K)=1$.
If we decrease the size of a compact open subgroup $L<K$ along a neighbourhood basis $\mathcal{B}$ of $\{e\}$ consisting of compact open subgroups of $K$, the projections $p_L$ defined in \eqref{eqn:subgroup-invar-projection} increase strongly to the identity.
Moreover, note that $\tau(p_L)=\mu(L)^{-1}=[K:L]$.
By \cite[Lemma~A.16]{kpv13}, this leaves us with
\begin{align*}
    \beta^{(2)}_n(G) &= \dim_{L(G)} H^n_c(G; L^2(G)) \\
    &= \sup_{L\in\mathcal{B}}\, \tau(p_L)\dim_{p_L L(G) p_L} p_L\left(H^n_c(G; L^2(G))\right)\\
    &= \sup_{L\in\mathcal{B}}\, [K:L]\dim_{p_L L(G) p_L} H^n_c(G; \ell^2(L\backslash G))\,. \numberthis\label{eqn:subgroup-shrinking-formula}
\end{align*}
In many respects, the corners $p_L L(G)p_L$ behave like the group von Neumann algebras of discrete groups.
Combining this with \eqref{eqn:cont-cohom-algebr}, one can then apply techniques from the discrete setting to understand the right-hand side of \eqref{eqn:subgroup-shrinking-formula}.

\begin{remark}
\label{rem:hecke-neighbourhood-basis}
It is a basic result in the theory of locally compact totally disconnected groups that the collection of \textit{all} compact open subgroups of $K$ always forms a neighbourhood basis of the identity.
However, sometimes more convenient choices are available.
If we start with a discrete Hecke pair $\Lambda<\Gamma$ with Schlichting completion $K<G$, the conjugates $\{gKg^{-1}\mid g\in \Gamma\}$ yield a neighbourhood subbasis of the identity in $G$ \cite[Definition~4.3, Proposition~4.5]{klq}.
Hence, we can take $\mathcal{B}$ to be the filter of subgroups of the form $K\cap g_1 Kg_1^{-1}\cap \cdots g_n K g_n^{-1}$ for $g_1,\ldots,g_n\in \Gamma$.
This choice of $\mathcal{B}$ turns out to be useful in \autoref{sec:hecke-l2betti-comp}.
\end{remark}

\subsection{Homological algebra}
\label{sec:homological-algebra}
\subsubsection{Bicomplexes in homological algebra}
On the algebraic side, the main ingredient in the proof of \thref{thm:vanishing-theorem} is a well-established tool from homological algebra.
To keep this paper as self-contained as possible, we recall the necessary preliminaries here---this also gives us the opportunity to justify why these techniques can be applied to rigorously compute homology ``up to dimension zero''.
\begin{definition}
    A \textit{bicomplex} in an abelian category $\mathcal{C}$ is a triple $((C_{p,q})_{p,q\in\ZZ},\partial^h, \partial^v)$ where $\partial^h: C_{p,q}\to C_{p,q-1}$ and $\partial^v: C_{p,q}\to C_{p-1,q}$ are referred to as the horizontal and vertical differentials, respectively.
    As usual, we suppress the sub-indices on the differentials, and we simply write $C_{\bullet,\bullet}$ to refer to the entire triple if there is no danger of confusion.
    The differentials should satisfy $\partial^h\partial^h=0$, $\partial^v\partial^v=0$ and $\partial^h\partial^v=-\partial^v\partial^h$, i.e.\@ all columns $C_{\bullet,q}$ and all rows $C_{p,\bullet}$ should be chain complexes in $\mathcal{C}$.

    A \textit{first quadrant bicomplex} is a bicomplex where $C_{p,q}=0$ whenever $p<0$ or $q<0$.
\end{definition}
Our principal reference for this section is \cite{weibel}.
\begin{remark}
    The abelian categories we work with will always be categories of modules, so we will only consider \textit{concrete} abelian categories in this section (i.e.\@ subcategories of the category of abelian groups).
\end{remark}
\begin{remark}
    As explained in \cite[\S~1.2.5]{weibel}, the category of bicomplexes in $\mathcal{C}$ is equivalent to the category of chain complexes that themselves consist of chain complexes in $\mathcal{C}$.
    One only needs to transform the anticommutativity condition $\partial^h\partial^v=-\partial^v\partial^h$ into $\partial^h\partial^v=\partial^v\partial^h$, which can be accomplished by a sign change in the differentials.
\end{remark}

The diagonals of a bicomplex form the so-called \textit{total complex} $\Tot(C_{\bullet,\bullet})$ with terms
\[
    \Tot(C_{\bullet,\bullet})_n = \bigoplus_{n=p+q} C_{p,q}
\]
and differentials $\partial^{\mathrm{Tot}} = \partial^h + \partial^v$.
The anticommutativity condition $\partial^h\partial^v=-\partial^v\partial^h$ ensures that $\partial^2=0$.

An important feature of bicomplexes is that they can serve as a way of combining and comparing complexes.
The tensor product bicomplex construction is one such mechanism.
\begin{definition}
    Let $R$ be a ring (with or without unit), and let $Q_\bullet$ and $P_\bullet$ be chain complexes of left $R$-modules and right $R$-modules, respectively.
    The \textit{tensor product bicomplex} $Q_\bullet \otimes_R P_\bullet$ has terms $Q_q\otimes_R P_p$, with vertical differentials given by $\I^q\otimes\partial^P$ and horizontal differentials given by $(-1)^p(\partial^Q\otimes \I^p)$.
\end{definition}
In many interesting cases, one can compare the homology of $Q_\bullet$ and $P_\bullet$ by means of the total complex $\Tot(Q_\bullet \otimes_R P_\bullet)$.
One typical such application is the standard proof of the theorem stating that one can compute $\Tor^R_n(A,B)$ by resolving either variable \cite[see][\S~2.7]{weibel}.

The true power of bicomplexes lies within the framework of spectral sequences.
However, for our purposes, the following more elementary result is already sufficient.
\begin{proposition}[Acyclic assembly lemma]
\label{thm:bicomplex-aal}
    Let $C_{\bullet,\bullet}$ be a first-quadrant bicomplex in some concrete abelian category and assume that $N\in\NN$ is such that every row $C_{k,\bullet}$ (resp.\@ column $C_{\bullet,k}$) is exact at position $C_{k,m}$ (resp.\@ $C_{m,k}$) for all $m\leq N$.
    Then $\Tot(C_{\bullet,\bullet})$ is exact at position $\Tot(C_{\bullet,\bullet})_m$ for all $m\leq N$.
\begin{proof}
    See e.g.\@ \cite[\S~2.7.3]{weibel} for a detailed proof in the case where one takes the rows (or columns) to be exact everywhere (i.e.\@ the case $N=\infty$).
    Since the argument will be useful for us later, we provide it here in the case where the rows $C_{k,\bullet}$ are taken to be exact at position $C_{k,m}$ for all $m\leq N$.
    The other case is equivalent after transposing the bicomplex.

    Fix $0\leq m\leq N$.
    A cycle $\xi\in\Tot(C_{\bullet,\bullet})_m$ consists of elements $\xi_n\in C_{m-n,n}$ satisfying the relation $\partial^h(\xi_{n+1})+\partial^v(\xi_n)=0$ for all $n\geq 0$.
    To show that $\xi=\partial^{\mathrm{Tot}}(\eta)$ for some $\eta\in\Tot(C_{\bullet,\bullet})_{m+1}$, we have to exhibit elements $\eta_n\in C_{m-n+1,n}$ satisfying
    \begin{equation}
    \label{eqn:aal-eta-relation-inductive}
        \xi_n = \partial^h(\eta_{n+1}) + \partial^v(\eta_n)
    \end{equation}
    for all $n\in\ZZ$.
    Note that we necessarily have to choose $\eta_n=0$ when $n>m+1$ or $n<0$, since $C_{m-n+1,n}=0$ in this case.
    Hence, we only have to provide $\eta_0,\ldots,\eta_{m+1}$, and show that the relation \eqref{eqn:aal-eta-relation-inductive} holds for $n\in\{0,\ldots,m\}$.
    For other values of $n$, the relation \eqref{eqn:aal-eta-relation-inductive} trivialises.

    We will construct the remaining $\eta_0,\ldots,\eta_{m+1}$ satisfying \eqref{eqn:aal-eta-relation-inductive} inductively, starting by putting $\eta_0=0$.
    In order for \eqref{eqn:aal-eta-relation-inductive} to hold when $n=0$, we have to choose $\eta_1\in C_{m,1}$ such that $\partial^h(\eta_1)=\xi_0$.
    This is possible because the row $C_{m,\bullet}$ is exact at position $C_{m,0}$ and $\partial^h(\xi_0)=0$ automatically.
    Suppose now that we have constructed $\eta_0,\ldots,\eta_k$ satisfying \eqref{eqn:aal-eta-relation-inductive} for $0\leq n\leq k-1$ where $1\leq k\leq m$.
    Then we compute
    \begin{equation}
    \label{eqn:aal-eta-relation-transfer}
    \begin{aligned}
        0 &= \partial^h(\xi_k) + \partial^v(\xi_{k-1}) 
        =\partial^h(\xi_k) + \partial^v(\xi_{k-1} - \partial^v(\eta_{k-1}))\\
        &=\partial^h(\xi_k) + \partial^v(\partial^h(\eta_k))
        =\partial^h(\xi_k-\partial^v(\eta_k))\,.
    \end{aligned}
    \end{equation}
    But the row $C_{m-k,\bullet}$ is exact at position $C_{m-k,k}$, since $k\leq m\leq N$. 
    This means that there exists $\eta_{k+1}\in C_{m-n,k+1}$ such that $\partial^h(\eta_{k+1})=\xi_k-\partial^v(\eta_k)$, exactly as required by \eqref{eqn:aal-eta-relation-inductive}.
\end{proof}
\end{proposition}
We will also need the following variant for dimension-exact sequences.
This result is probably known to the experts, but we provide a proof for completeness.
\begin{proposition}[Dimension-acyclic assembly lemma]
\label{thm:bicomplex-dimension-aal}
    Let $M$ be a von Neumann algebra with n.s.f.\@ trace $\tau$.
    Let $C_{\bullet,\bullet}$ be a first-quadrant bicomplex of locally finite (left or right) $M$-modules, and assume that $N\in\NN$ is such that every row $C_{k,\bullet}$ (resp.\@ column $C_{\bullet,k}$) is dimension exact at position $C_{k,m}$ (resp.\@ $C_{m,k}$) for all $m\leq N$.
    Then $\Tot(C_{\bullet,\bullet})$ is dimension exact at position $\Tot(C_{\bullet,\bullet})_m$ for all $m\leq N$.
\begin{proof}
    As in our proof of \thref{thm:bicomplex-aal}, we only consider the case where the rows are taken to be dimension exact up to the $N$th position.
    Fix $0\leq m\leq N$, $\epsilon>0$ and a cycle $\xi\in\Tot(C_{\bullet,\bullet})_m$, which we again decompose as elements $\xi_n\in C_{m-n,n}$ satisfying $\partial^h(\xi_{n+1})+\partial^v(\xi_n)=0$ for all $n\geq 0$.
    By \thref{thm:zero-dim-criterion}, it suffices to find a projection $p\in M$ and an element $\eta\in\Tot(C_{\bullet,\bullet})_m$ such that $\tau(p)<\epsilon$ and $(\I-p)\xi=\partial^{\mathrm{Tot}}(\eta)$.
    In the same spirit as the proof of \thref{thm:bicomplex-aal} we will exhibit projections $p_0\leq\cdots\leq p_m$ and elements $\eta_0,\eta_1,\ldots,\eta_{m+1}$ such that
    \begin{equation}
    \label{eqn:aal-approx-eta-relation-inductive}
        (\I-p_n)\xi_n = (\I-p_n)(\partial^h(\eta_{n+1}) + \partial^v(\eta_n))\,
        \qquad\text{and}\qquad \tau(p_n)<\epsilon 2^{n-m}
    \end{equation}
    for all $n\in\{0,\ldots,m\}$.
    Taking $p=p_m$, we can then assemble $(\I-p)\eta_0,\ldots,(\I-p)\eta_{m+1}$ into an element of $\Tot(C_{\bullet,\bullet})_{m+1}$ with the required property.
    Again, we bootstrap the induction by putting $\eta_0=0$.
    Applying \thref{thm:zero-dim-criterion} to the dimension-exactness of $C_{m,\bullet}$ at $C_{m,0}$, we find a projection $p_0\in M$ and $\eta_1\in C_{m,1}$ such that $\tau(p_0)<\epsilon 2^{-m}$ and $(\I-p_0)\xi=\partial_h(\eta_1)$.
    Then \eqref{eqn:aal-approx-eta-relation-inductive} is satisfied for $n=0$.
    For the induction step, suppose that $1\leq k\leq m$ and that we have constructed appropriate projections $p_0\leq\cdots\leq p_{k-1}$ and $\eta_0,\ldots,\eta_k$ satisfying \eqref{eqn:aal-approx-eta-relation-inductive} for $0\leq n\leq k-1$.
    Analogously to \eqref{eqn:aal-eta-relation-transfer}, one then finds that
    \[
        \partial^h\left[ (\I-p_{k-1})(\xi_k-\partial^v(\eta_k))\right] = 0\,.
    \]
    By dimension-exactness of $C_{m-k,\bullet}$ at position $C_{m-k,k}$, this implies the existence of a projection $q\in M$ and $\eta_{k+1}\in C_{m-k,k+1}$ such that $\tau(q)<\epsilon 2^{k-m-1}$ and 
    \[
        (\I-q)(\I-p_{k-1})(\xi_k-\partial^v(\eta_k))=\partial^h(\eta_{k+1})\,.
    \]
    If we now put $p_k=q\vee p_{k-1}$, the requirement \eqref{eqn:aal-approx-eta-relation-inductive} is satisfied for $n=k$, completing the inductive step.
\end{proof}
\end{proposition}
Alternatively, one could apply \thref{thm:bicomplex-aal} in the Serre quotient category of $M$-modules modulo zero-dimensional objects (cfr.\@ the discussion in \cite[\S~1.2]{thom-rank-metric}) to arrive at the same conclusion, but the above argument gives a more hands-on approach.

The homological algebra section of this paper can be reformulated somewhat more succinctly using the language of spectral sequences without substantially changing the general structure of the proofs, but in the interest of accessibility we have opted to make all homology computations as explicit as possible.

\subsubsection{Mapping cones}
Consider two chain complexes $C_\bullet$ and $D_\bullet$ in some concrete abelian category, and a chain map $\phi_\bullet: C_\bullet\to D_\bullet$.
Then $\phi$ induces natural maps $\tilde{\phi}_*: H_*(C_\bullet)\to H_*(D_\bullet)$ on homology.
It is then natural to ask when $\tilde{\phi}$ is an isomorphism.
Mapping cones \cite[\S~1.5]{weibel} provide an answer by associating a new complex to $\phi$.
\begin{definition}
    Given a chain map $\phi_\bullet:C_\bullet\to D_\bullet$, define $\cone(\phi)_n = C_{n-1}\oplus D_n$ with differentials
    \[
        \partial: \cone(\phi)_n\to \cone(\phi)_{n-1}: \xi\oplus\eta \mapsto -\partial(\xi) \oplus (\partial(\eta)-\phi_{n-1}(\xi))\,.
    \]
    The complex $\cone(\phi)_\bullet$ is called the \textit{mapping cone} of $\phi$.
\end{definition}

The main result concerning mapping cones is the following.
\begin{proposition}
\label{thm:mapping-cone-iso}
A chain map $\phi_\bullet:C_\bullet\to D_\bullet$ induces an isomorphism on homology if and only if the mapping cone of $\phi$ is exact.
\begin{proof}
    See \cite[Corollary~1.5.4]{weibel}.
\end{proof}
\end{proposition}
In this sense, the homology groups of $\cone(\phi)_\bullet$ can be interpreted as measuring the extent to which $\tilde{\phi}_*: H_*(C_\bullet)\to H_*(D_\bullet)$ fails to be an isomorphism.

Using mapping cones and \thref{thm:mapping-cone-iso}, one can prove the following fact about first quadrant bicomplexes.
\begin{lemma}
    \label{thm:augmented-bicomplex-homology} 
    Let $C_{\bullet,\bullet}$ be a first quadrant bicomplex in some concrete abelian category.
    Let $C^0_{\bullet,\bullet}$ be the complex obtained from $C_{\bullet,\bullet}$ by replacing the zeroth column by zeroes.
    If $\Tot(C_{\bullet,\bullet})$ is exact, then there are natural isomorphisms
    \begin{equation}
    \label{eqn:augmented-bicomplex-homology}
        H_n(C_{\bullet,0}) \cong H_{n+1}(\Tot(C^0_{\bullet,\bullet}))
    \end{equation}
    for all $n\geq 0$.
\begin{proof}
    This is essentially the argument used in the proof of \cite[Theorem~2.7.2]{weibel}.
    We denote the horizontal (resp.\@ vertical) differentials of $C_{\bullet,\bullet}$ by $\partial^h$ (resp.\@ $\partial^v$).
    Write $A_p=C_{p-1,0}$, where we view $A_\bullet$ as a chain complex with differentials given by $\partial^v$. 
    Of course we can then canonically identify $H_{n+1}(A_\bullet)\cong H_n(C_{\bullet,0})$.
    Consider the maps 
    \begin{align*}
        \phi_n: \Tot(C_{\bullet,\bullet}^0)_n \to A_n: \bigoplus_{i\in\NN} C_{n-i,i}^0 \ni (\xi_i)_{i\in\NN} \mapsto -\partial^h(\xi_1).
    \end{align*}
    A routine computation shows that this gives rise to a chain map $\phi$ from $\Tot(C_{\bullet,\bullet}^0)$ to $A_\bullet$.
    On the other hand, there are isomorphisms
    \begin{align*}
        \cone(\phi)_{n+1} &\to \Tot(C_{\bullet,\bullet})_n: \left[\left(\bigoplus_{i\in\NN} C_{n-i,i}^0\right)\oplus C_{n,0}\right]\ni (\xi_i)_{i\in\NN} \oplus \eta\mapsto (\zeta_i)_{i\in\NN}\in \left(\bigoplus_{i\in\NN} C_{n-i,i}\right)\\
        &\text{where}\qquad \zeta_i=\begin{cases}
            \eta & i = 0\,, \\
            (-1)^n\xi_i & i \neq 0\,,
        \end{cases}
    \end{align*}
    realising a chain isomorphism from $\cone(\phi)_\bullet$ to the total complex of $C_{\bullet,\bullet}$ shifted by one.
    Hence, if said total complex is exact, \thref{thm:mapping-cone-iso} tells us that $\phi$ induces the required isomorphisms \eqref{eqn:augmented-bicomplex-homology} on homology.
\end{proof}
\end{lemma}

\section{A vanishing theorem for almost-normal inclusions of rigid \texorpdfstring{$C^*$}{C*}-tensor categories}
\subsection{Normal and almost-normal inclusions of rigid \texorpdfstring{$C^*$}{C*}-tensor categories}
Consider a rigid $C^*$-tensor category $\mathcal{D}$ and a full $C^*$-tensor subcategory $\mathcal{C}$.
For $\alpha,\beta\in\Irr(\mathcal{D})$, write $\alpha\sim\beta$ whenever there exist $\pi,\pi'\in\Irr(\mathcal{C})$ such that $\alpha\hookrightarrow \pi\beta\pi'$.
The relation $\sim$ is an equivalence relation, and we denote the equivalence class of $\alpha$ by $\mathcal{C}\alpha\mathcal{C}$.
We call these the \textit{double cosets} of $\mathcal{C}\subset\mathcal{D}$.
Obviously, $\alpha\mathcal{C}=\beta\mathcal{C}$ also implies that $\mathcal{C}\alpha\mathcal{C}=\mathcal{C}\beta\mathcal{C}$, so each double coset is a disjoint union of one-sided cosets.
This leads us to the following natural definition.
\begin{definition}
\label{def:normal-almost-normal}
    A full $C^*$-tensor subcategory $\mathcal{C}$ of a rigid $C^*$-tensor category $\mathcal{D}$ is said to be \textit{almost normal} if $\mathcal{C}\alpha\mathcal{C}$ is a union of finitely many left (resp.\@ right\footnote{The left and right versions are equivalent.}) cosets for all $\alpha\in\Irr(\mathcal{D})$.
    If $\mathcal{C}\alpha=\alpha\mathcal{C}$ for all $\alpha\in\Irr(\mathcal{D})$, we call the inclusion \textit{normal}.
\end{definition}

\begin{remark}
    \label{rem:almost-normal-qr-version}
    Consider an inclusion of rigid $C^*$-tensor categories $\mathcal{C}\subset\mathcal{D}$ with $\mathcal{D}$ finitely generated.
    If $\mathcal{C}\subset\mathcal{D}$ is almost normal, this inclusion gives rise to a tower $T\subset S\subset R$ of factors of type $\mathrm{II}_1$ in which all inclusions are irreducible and quasi-regular.
    To see why this is the case, fix a self-dual generating object $\alpha\in\mathcal{D}$ and use Popa's reconstruction theorem \cite{popa-lambda-lattices} to choose an extremal finite-index inclusion of type $\mathrm{II}_1$ factors $N\subset M$ such that $L^2(M_1)\mapsto\alpha$ extends to a unitary monoidal equivalence between $\mathcal{D}$ and the category of bifinite $M$-$M$-bimodules appearing in the Jones tower $N\subset M\subset M_1\subset\cdots$ \cite[see also][Example~5.1]{neshveyev-yamashita}.
    Fix an $M$-$M$-bimodule realisation $\mathcal{H}_\alpha$ of every $\alpha\in\Irr(\mathcal{D})$. 
    For simplicity, we choose $\mathcal{H}_\epsilon=L^2(M)$.
    Put $T=M\otimes M^{\mathrm{op}}$ and define
    \[
        S_0 = \bigoplus_{\alpha\in\Irr(\mathcal{C})} \mathcal{H}_\alpha^0\otimes\overline{\mathcal{H}_{\alpha}^0}
        \qquad\subset\qquad
        R_0 = \bigoplus_{\alpha\in\Irr(\mathcal{D})} \mathcal{H}_\alpha^0\otimes\overline{\mathcal{H}_{\alpha}^0}\,.
    \]
    Here, the direct sums and tensor products are algebraic and the symbol $-^0$ denotes passage to bounded vectors.
    The Longo--Rehren approach to Popa's symmetric enveloping algebra \cite[see also \citenum{pv-repr-subfactors}, Remark~2.7]{longo-rehren,masuda-lr} then provides $*$-algebra structures and traces on $S_0$ and $R_0$ that turn the inclusion $S_0\subset R_0$ into an inclusion of tracial $*$-algebras.
    By passing to the GNS representation associated with this canonical trace, we can complete $S_0$ and $R_0$ to von Neumann algebras $S$ and $R$.
    In this way, we get a tower of von Neumann algebras $T\subset S\subset R$.

    Since the inclusion $T\subset R$ is irreducible, $T\subset S\subset R$ is a tower of irreducible inclusions of type $\mathrm{II}_1$ factors.
    The inclusions $T\subset S$ and $T\subset R$ are quasi-regular and unimodular, with the $L^2$-Betti numbers recovering those of $\mathcal{C}$ and $\mathcal{D}$, respectively \cite[see][Proposition~3.12]{psv-cohom}.
    One then checks that the almost-normality of $\mathcal{C}$ in $\mathcal{D}$ implies that $R_0\subset\mathrm{QN}_R(S)$, so $\mathrm{QN}_R(S)$ is weakly dense in $R$, as required.
\end{remark}

\begin{example}
\label{exa:bruguieres-natale-normal}
    Let $\mathcal{C}$, $\mathcal{D}$ be rigid $C^*$-tensor categories, and $F:\mathcal{C}\to \mathcal{D}$ a unitary tensor functor.
    Recall that an object in a rigid $C^*$-tensor category is called \textit{trivial} if it is isomorphic to a direct sum of copies of the tensor unit $\epsilon$.
    Suppose that $F$ is normal in the sense of \cite{bruguieres-natale}, i.e.\@ for any object $\alpha\in\mathcal{C}$, there is a subobject $\alpha_0$ of $\alpha$ such that $F(\alpha_0)$ is the largest trivial subobject of $F(\alpha)$.
    Then $\mathcal{C}_0=\{\alpha\in \mathcal{C}\mid F(\alpha)\text{ trivial}\}$ is a normal subcategory of $\mathcal{C}$ in the sense of \thref{def:normal-almost-normal}.

    Indeed, by hypothesis, there exist nonzero natural numbers $(n_\gamma)_{\gamma\in\Irr(\mathcal{C}_0)}$ such that $F(\gamma)=\epsilon^{\oplus n_\gamma}$ for all $\gamma\in\Irr(\mathcal{C}_0)$.
    We also know that $F(\gamma)$ contains no trivial summands when $\gamma\in\Irr(\mathcal{C})\setminus \Irr(\mathcal{C}_0)$.
    This implies that
    \[
        \dim_{\CC} (F(\alpha),F(\beta)) = \mult(\epsilon, F(\bar{\alpha}\beta)) = \sum_{\gamma\in\Irr(\mathcal{C}_0)} n_\gamma \mult(\gamma, \bar{\alpha}\beta) 
    \]
    for all $\alpha,\beta\in\Irr(\mathcal{C})$. The right-hand side of this expression is nonzero if and only if $\beta\in\alpha\mathcal{C}_0$.
    However, taking conjugates, we also get that
    \[
        \dim_{\CC} (F(\alpha),F(\beta)) = \dim_{\CC} (F(\bar{\alpha}),F(\bar{\beta})) =\sum_{\gamma\in\Irr(\mathcal{C}_0)} n_\gamma \mult(\gamma, \alpha\bar{\beta}) 
    \]
    which is nonzero if and only if $\beta\in\mathcal{C}_0\alpha$.
    It follows that $\mathcal{C}_0\alpha=\alpha\mathcal{C}_0$ for all $\alpha\in\Irr(\mathcal{C})$, as claimed.

    Conversely, it is not true that any normal subcategory in the sense of \thref{def:normal-almost-normal} is the kernel of a normal functor.
    For one, the restriction of a normal tensor functor to its kernel $\mathcal{C}_0$ induces a fibre functor from $\mathcal{C}_0$ to $\mathrm{Hilb}_f$, and not all $C^*$-tensor categories admit such functors \cite[Example~2.3.7]{neshveyev-tuset}.
    This obstruction is arguably of a purely technical nature, and hence somewhat unsatisfying.
    However, since the existence of a normal functor with kernel $\mathcal{C}_0$ depends on much more than the fusion rules of $\mathcal{C}$, we strongly believe that there are many natural counterexamples that \textit{do} admit fibre functors, although we do not have any concrete such example.
\end{example}

\begin{remark}
    \label{rem:normal-quantum-subgroup}
    Given a compact quantum group $\mathbb{G}$ and a normal compact quantum subgroup $\mathbb{H}$ in the sense of \cite{wang-free-products,wang-scqg}, the restriction functor $\Rep_f(\mathbb{G})\to \Rep_f(\mathbb{H})$ is normal in the sense of \cite{bruguieres-natale}, and the subcategory $\mathcal{C}_0\subset\Rep_f(\mathbb{G})$ in \thref{exa:bruguieres-natale-normal} consists precisely of all representations with matrix coefficients in the subalgebra $C(\mathbb{G}/\mathbb{H})\subset C(\mathbb{G})$ \cite[Proposition~2.1]{wang-scqg}.
    In other words, $\mathcal{C}_0$ can be identified with the representation category $\Rep_f(\mathbb{G}/\mathbb{H})$. 
    We conclude that the normal compact quantum subgroup $\mathbb{H}\subset\mathbb{G}$ gives rise to a normal inclusion of representation categories $\Rep_f(\mathbb{G}/\mathbb{H})\subset \Rep_f(\mathbb{G})$.
\end{remark}

The categories constructed in \cite{arano-vaes} from totally disconnected groups provide another source of interesting examples.
For the convenience of the reader, we briefly explain how to construct these categories here.
Let $G$ be a locally compact totally disconnected group with a compact open subgroup $K$.
An $\ell^\infty(G/K)$-$G$-$\ell^\infty(G/K)$-\textit{module} consists of a Hilbert $\ell^\infty(G/K)$-$\ell^\infty(G/K)$-bimodule $\mathcal{H}$, equipped with a continuous unitary representation of $G$ that is equivariant with respect to the bimodule structure in the obvious way, i.e.\@
\[
    g\cdot (\chi_{sK}\cdot \xi\cdot\chi_{tK})=\chi_{gsK}\cdot (g\cdot \xi)\cdot \chi_{gtK}\,
\]
for $\xi\in\mathcal{H}$ and $s,t,g\in G$.
Taking the morphisms between $\ell^\infty(G/K)$-$G$-$\ell^\infty(G/K)$-modules to be all bounded linear maps intertwining the relevant actions, one gets the category $\mathcal{C}(K<G)$.
The Connes tensor product over $\ell^\infty(G/K)$ turns this category into a (typically nonrigid) $C^*$-tensor category.

A $\ell^\infty(G/K)$-$G$-$\ell^\infty(G/K)$-module $\mathcal{H}$ is said to be of \textit{finite rank} if $\chi_K\cdot\mathcal{H}$ (or equivalently $\mathcal{H}\cdot\chi_K$) is finite-dimensional.
The subcategory of $\mathcal{C}(K<G)$ given by all finite-rank $\ell^\infty(G/K)$-$G$-$\ell^\infty(G/K)$-modules is a rigid $C^*$-tensor category, which we will denote by $\mathcal{C}_f(K<G)$.
As stated in \cite{arano-vaes}, the finite-dimensional unitary representations of $K$ form a natural subcategory of $\mathcal{C}_f(K<G)$.
It turns out that this subcategory is normal in $\mathcal{C}_f(K<G)$.
This is the content of the next result.  
\begin{proposition}
\label{thm:av-normality}
    Let $G$ be a second-countable locally compact totally disconnected group, and $K<G$ a compact open subgroup.
    Let $\mathcal{C}$ be the full subcategory of $\mathcal{C}_f(K<G)$ consisting of all objects $\mathcal{H}\in\mathcal{C}_f(K<G)$ with the property that $\chi_{gK}\cdot\mathcal{H}=\mathcal{H}\cdot\chi_{gK}$ for all $g\in G$.
    Then the following hold:
    \begin{enumerate}[(i)]
        \item $\mathcal{C}$ is unitarily monoidally equivalent to $\Rep_f(K)$ \cite{arano-vaes},
        \item $\mathcal{C}$ is normal in $\mathcal{C}_f(K<G)$, and
        \item if moreover $G$ can be decomposed as $G=K\Gamma$ with $\Gamma<G$ a discrete subgroup such that $\Gamma\cap K=\{e\}$, then $\mathcal{C}$ is normal in $\mathcal{C}_f(K<G)$ in the stronger sense discussed in \thref{exa:bruguieres-natale-normal}.
    \end{enumerate}
\end{proposition}
The extra hypothesis in (iii) asserts that $(K,\Gamma)$ is a \textit{matched pair} for $G$, in the sense considered in \cite{kac-matched-pairs,baaj-skandalis,vaes-vainerman,fmp-bicrossed-products}.
Many examples of matched pairs $G=K\Gamma$ with $G$ totally disconnected, $K$ infinite compact and $\Gamma$ infinite discrete can be found by considering simply transitive actions on Euclidean buildings of type $\tilde{A}_2$ \cite{cmsz1,vv-property-t}.
Given such a matched pair, the associated bicrossed product quantum group $\mathbb{H}$ has $\mathcal{C}_f(K<G)$ as its representation category \cite[Proposition~6.1]{vv-property-t}.
\begin{proof}[Proof of \thref{thm:av-normality}]
    The first claim is explained in the discussion preceding Theorem~3.1 in \cite{arano-vaes}.
    To prove (ii), we show that the right (and left) cosets of $\mathcal{C}$ in $\mathcal{C}_f(K<G)$ are given by
    \begin{align*}
        \mathcal{C}_{KxK}&=\{\mathcal{H}\in\Irr(\mathcal{C}_f(K<G))\mid \chi_{KxK}\cdot\mathcal{H}\cdot\chi_K=\mathcal{H}\cdot\chi_K\}\\
        &=\{\mathcal{H}\in\Irr(\mathcal{C}_f(K<G))\mid \chi_K\cdot\mathcal{H}\cdot \chi_{Kx^{-1}K}=\chi_K\cdot\mathcal{H}\}
    \end{align*}
    for $KxK\in K\backslash G/K$.
    We now proceed to describe a set of canonical representatives of $\mathcal{C}_{KxK}$.
    For $x\in G$, define 
    \begin{equation}
    \label{eqn:av-cat-coset-canon-rep}
        \mathcal{K}_x = \ell^2(G/K\cap x^{-1} K x)
    \end{equation}
    with $G$ acting by left translation, and the $\ell^\infty(G/K)$-bimodule structure given by (cfr.\@ the discussion in \cite[\S~3]{arano-vaes})
    \[
        (F_1\cdot \xi \cdot F_2)(g) = F_1(g) \xi(g) F_2(gx^{-1}).
    \]
    It is not difficult to see that $\mathcal{K}_x$ is an irreducible element of $\mathcal{C}_{KxK}$.
    Moreover, the isomorphism class of $\mathcal{K}_x$ only depends on $KxK$.

    We now claim that $\mathcal{C}\mathcal{K}_x=\mathcal{K}_x\mathcal{C}=\mathcal{C}_{KxK}$, and that these cover all cosets of $\mathcal{C}$ in $\mathcal{C}_f(K<G)$.  
    The inclusions $\mathcal{C}\mathcal{K}_x\subset\mathcal{C}_{KxK}$ and $\mathcal{K}_x\mathcal{C}\subset\mathcal{C}_{KxK}$ are clear.
    To prove the reverse inclusions, fix an arbitrary object $\mathcal{H}\in\mathcal{C}_{KxK}$, and  define
    \[
        \mathcal{L}_1=\ell^2(G/K) \cftimes{\ell^\infty(G/K)} \mathcal{H}
        \qquad\qquad
        \mathcal{L}_2=\mathcal{H} \cftimes{\ell^\infty(G/K)} \ell^2(G/K)
    \]
    where $G$ acts diagonally, and the left- and right $\ell^\infty(G/K)$-module structures are both defined by multiplication on the $\ell^2(G/K)$-leg.
    By construction $\mathcal{L}_1,\mathcal{L}_2\in\mathcal{C}$.
    Under the identification of $\mathcal{C}$ with $\Rep_f(K)$, $\mathcal{L}_1$ (resp.\@ $\mathcal{L}_2$) is the object corresponding to the representation of $K$ on $\chi_K\cdot\mathcal{H}$ (resp.\@ $\mathcal{H}\cdot\chi_K$).

    We now embed $\mathcal{H}$ inside both $\mathcal{L}_1\mathcal{K}_x$ and $\mathcal{K}_x\mathcal{L}_2$, via
    \begin{align*}
        &\phi_1:\mathcal{H}\to \mathcal{L}_1\mathcal{K}_x: \\
        &\qquad\qquad\xi\mapsto \sum_{g(K\cap x^{-1}Kx)\in G/(K\cap x^{-1}Kx)} \left(\chi_{gK}\otimes_{\ell^\infty(G/K)} \xi\cdot\chi_{gx^{-1}K}\right)\otimes_{\ell^\infty(G/K)} \chi_{g(K\cap x^{-1}Kx)}\,,\\
        &\phi_2:\mathcal{H}\to \mathcal{K}_x\mathcal{L}_2: \\
        &\qquad\qquad\xi\mapsto \sum_{g(K\cap xKx^{-1})\in G/(K\cap xKx^{-1})} \chi_{gx(K\cap x^{-1}Kx)} \otimes_{\ell^\infty(G/K)} \left(\chi_{gxK}\cdot\xi\otimes_{\ell^\infty(G/K)} \chi_{gK}\right)\,.
    \end{align*}
    To check that $\phi_1$ and $\phi_2$ are morphisms in $\mathcal{C}_f(K<G)$, one first observes that both maps intertwine the respective actions of $G$.
    It is then sufficient to verify that $\phi_1$ and $\phi_2$ commute with $\chi_K\cdot-$ and $-\cdot\chi_K$, which is straightforward.
    The injectivity of $\phi_1$ and $\phi_2$ follow from the assumption that $\mathcal{H}\in\mathcal{C}_{KxK}$.

    Finally, we still have to argue that every irreducible object in $\mathcal{C}_f(K<G)$ occurs in some $\mathcal{C}_{KxK}$.
    For any $\mathcal{H}\in\mathcal{C}_f(K<G)$ and any double coset $KxK$, there is a canonical projection operator $p_{KxK}$ on $\mathcal{H}$ defined by
    \[
        p_{KxK}(\xi) = \sum_{g(K\cap xKx^{-1})\in G/K\cap xKx^{-1}} \chi_{gxK} \cdot \xi\cdot \chi_{gK}
    \]
    for $\xi\in\mathcal{H}$.
    Note that $p_{KxK}$ is an endomorphism of $\mathcal{H}$ for all $x\in G$, and that
    \begin{equation}
    \label{eqn:pkxk-relation}
        p_{KxK}(\xi\cdot\chi_K) = \chi_{KxK}\cdot\xi\cdot\chi_K
    \end{equation}
    for all $x\in G$ and $\xi\in\mathcal{H}$.
    Additionally, one easily verifies that $(p_{KxK})_{KxK\in K\backslash G/K}$ constitutes a family of pairwise orthogonal projections summing to the identity.
    If $\mathcal{H}$ is irreducible, it follows that there is precisely one double coset $KxK$ such that $p_{KxK}$ is the identity on $\mathcal{H}$, so \eqref{eqn:pkxk-relation} tells us that $\mathcal{H}\in \mathcal{C}_{KxK}$.
    This concludes the proof of (ii).

    To prove (iii), we construct a normal tensor functor (in the sense of \cite{bruguieres-natale}) from $\mathcal{C}_f(K<G)$ to the category $\mathrm{Hilb}_f^\Gamma$ of finite-dimensional Hilbert spaces graded over $\Gamma$.
    Concretely, put
    \[
        \mathbf{F}:\mathcal{C}_f(K<G)\to\mathrm{Hilb}_f^\Gamma: \mathcal{H}\mapsto \bigoplus_{\gamma\in\Gamma} \chi_K\cdot \mathcal{H}\cdot \chi_{\gamma K}\,.
    \]
    By definition, every morphism $\phi:\mathcal{H}\to\mathcal{K}$ in $\mathcal{C}_f(K<G)$ restricts to a linear map $\chi_K\cdot\mathcal{H}\cdot\chi_{\gamma K}\to \chi_K\cdot\mathcal{K}\cdot\chi_{\gamma K}$ for all $\gamma\in\Gamma$, which turns $\mathbf{F}$ into a unitary functor.
    The unitary functor $\mathbf{F}$ admits a unitary monoidal structure given by the natural unitaries 
    \[
        U_{\mathcal{H},\mathcal{K}}: \mathbf{F}(\mathcal{H}\mathcal{K})\to \mathbf{F}(\mathcal{H})\mathbf{F}(\mathcal{K}): 
        \xi\mapsto \sum_{\gamma\in\Gamma} (\I\otimes\pi_{\mathcal{K}}(\gamma)^*\chi_{\gamma K})\xi\,.
    \]
    Indeed, a straightforward computation shows that $U_{\mathcal{H},\mathcal{K}}$ is a morphism of $\Gamma$-graded Hilbert spaces for all $\mathcal{H},\mathcal{K}\in\mathcal{C}_f(K<G)$.
    The fact that $U_{\mathcal{H},\mathcal{K}}$ yields a monoidal structure on $\mathbf{F}$ then follows by exactly the same computation as the one in \cite[Proposition~6.1]{vv-property-t}.

    Finally, observe that the largest subobject of $\mathcal{H}$ that is mapped to a trivial object in $\mathrm{Hilb}_f^\Gamma$ is precisely $p_K(\mathcal{H})$.
    It follows that $\mathbf{F}$ is normal with kernel $\mathcal{C}$, as claimed.
\end{proof}

\begin{remark}
\label{rem:bicrossed-product-exact-sequence}
    In more abstract terms, the above argument demonstrates that the fibre functor on $\mathcal{C}_f(K<G)$ constructed in \cite[Proposition~6.1]{vv-property-t} from the matched pair decomposition $G=K\Gamma$ factors through $\mathrm{Hilb}_f^\Gamma$.
    More precisely, it is the result of composing $\mathbf{F}:\mathcal{C}_f(K<G)\to\mathrm{Hilb}_f^\Gamma$ with the forgetful functor $\mathrm{Hilb}_f^\Gamma\to\mathrm{Hilb}_f$ that discards the grading.

    Since said fibre functor realises $\mathcal{C}_f(K<G)$ as the representation category of the bicrossed compact quantum group $\mathbb{H}$ associated with the matched pair $(K,\Gamma)$, this is essentially the $C^*$-tensor category version of the exact sequence $1\to \hat{K}\to \hat{\mathbb{H}}\to \Gamma\to 1$ of discrete quantum groups.
\end{remark}

\begin{remark}
    We believe that for general $K<G$, the canonical copy of $\Rep_f(K)$ in $\mathcal{C}_f(K<G)$ will not always be the kernel of any normal functor, but we do not have any examples where we can prove that this is the case.
    This is a question in the same spirit as the one asked in \cite[Remark~6.2]{vv-property-t} about the existence of fibre functors on $\mathcal{C}_f(K<G)$.
    In particular, one could ask whether $\Rep_f(K)$ being normal inside $\mathcal{C}_f(K<G)$ in this stronger sense forces $G$ to be unimodular, or if there are any other examples besides the ones coming from matched pairs.
\end{remark}

\subsection{Vanishing results}
This section is devoted to the proof of the following vanishing theorem for $L^2$-Betti numbers of rigid $C^*$-tensor categories.
\begin{theorem}
    \label{thm:vanishing-theorem}
    Consider a rigid $C^*$-tensor category $\mathcal{D}$ and an almost-normal subcategory $\mathcal{C}$.
    Given $N\in\NN$ such that $\beta_n^{(2)}(\mathcal{C})=0$ for all $0\leq n\leq N$, we also have that $\beta_n^{(2)}(\mathcal{D})=0$ for all $0\leq n\leq N$.
\end{theorem}
If $\mathcal{D}$ is a countable discrete group, we recover \cite[Corollary~1.4]{bfs}.
The statement for general almost-normal inclusions of discrete groups was proven in \cite{bfs}, but the case of normal subgroups had been covered much earlier \cite[see][Theorem~3.3(2)]{lueck-dimtheo-2}.
From a homological algebra perspective, the techniques we employ in our argument are conceptually rather similar to those used in \cite{bfs}.
The dimension-theoretic part of the argument is more subtle, in part due to the fact that we have no categorical counterpart to the orbit-stabiliser theorem for subgroups acting on cosets.
\begin{remark}
    While we do not have any example of an almost-normal inclusion of rigid $C^*$-tensor categories that is non-normal, infinite index and non-equivalent to an almost-normal inclusion of discrete groups, making the extra assumption of normality does not substantially simplify the proof of \thref{thm:vanishing-theorem} (see also \thref{rem:almost-normal-qr-version}).
    We have therefore opted to state the theorem in this somewhat more general language.
\end{remark}
\begin{remark}
    In the context of inclusions $\Lambda<\Gamma$ of discrete groups, it is worth noting that specialising to normal subgroups allows one to approach the problem from several other angles.
    For example, the proof of \cite[Theorem~3.3(2)]{lueck-dimtheo-2} applies the Leray--Serre spectral sequence to the fibration of classifying spaces $B(\Lambda)\to B(\Gamma)\to B(\Gamma/\Lambda)$.
    In addition to this, there is an analogue of \thref{thm:vanishing-theorem} for strongly normal inclusions of discrete measured groupoids, proved in \cite[Theorem~1.3]{sauer-thom} by making use of a Grothendieck spectral sequence analogous to the classical Hochschild--Serre spectral sequence.
    To the best of the author's knowledge, these techniques are not immediately applicable to the situation of \thref{thm:vanishing-theorem}, even if the subcategory $\mathcal{C}$ is normal in $\mathcal{D}$.
\end{remark}
\begin{remark} 
    If one is only concerned with the first $L^2$-Betti number, there are versions of \thref{thm:vanishing-theorem} for discrete groups that go back even further---see e.g.\@ \cite[Theorem~0.7]{lueck-hilbert-modules}, \cite[Th\'eor\`eme~6.8]{gaboriau} or \cite[Corollary~1]{bmv}.
    Theorem~5.12 in \cite{peterson-thom} gives a vanishing criterion for the first $L^2$-Betti number under much weaker conditions on the subgroup $\Lambda<\Gamma$.
\end{remark}
Before turning to the proof of \thref{thm:vanishing-theorem}, we first discuss a number of corollaries and examples.

\begin{corollary} 
    \label{thm:infinite-amenable-subcategory-vt}
    Let $\mathcal{C}$ be a rigid $C^*$-tensor category.
    If $\mathcal{C}$ admits an almost-normal amenable subcategory with infinitely many classes of irreducible objects, then all $L^2$-Betti numbers of $\mathcal{C}$ vanish.
    \begin{proof}
        This follows from \thref{thm:vanishing-theorem} combined with the fact that amenable rigid $C^*$-tensor categories with infinitely many irreducibles have vanishing $L^2$-Betti numbers \cite[see][Corollary~8.5]{psv-cohom}.
    \end{proof}
\end{corollary}

\begin{corollary}
\label{thm:dqg-vanishing-theorem}
    Let $\mathbb{G}$ be a compact quantum group of Kac type, and $\mathbb{H}$ a normal compact quantum subgroup.
    Given $N\in\NN$ such that $\beta_n^{(2)}(\widehat{\mathbb{G}/\mathbb{H}})=0$ for all $0\leq n\leq N$, we also have that $\beta_n^{(2)}(\hat{\mathbb{G}})=0$ for all $0\leq n\leq N$.
    In particular, if $\GG$ admits an infinite coamenable quotient, $\beta_n^{(2)}(\hat{\GG})=0$ for all $n\geq 0$.
    \begin{proof}
        The first statement follows from \thref{thm:vanishing-theorem} by \thref{rem:normal-quantum-subgroup} and \cite[Theorem~4.2]{krvv}.
        If $\mathbb{K}$ is an infinite coamenable quotient of $\mathbb{G}$, then \cite[Corollary~6.2]{kyed-coamenable} states that $\hat{\mathbb{K}}$ has vanishing $L^2$-Betti numbers, so the second statement follows from the first.
        Alternatively, one can also conclude the second statement from \thref{thm:infinite-amenable-subcategory-vt} and the fact that $\Rep_f(\mathbb{K})$ is an amenable rigid $C^*$-tensor category with infinitely many irreducibles. 
    \end{proof}
\end{corollary}

Using \thref{thm:infinite-amenable-subcategory-vt}, we can also prove that the rigid $C^*$-tensor categories introduced in \cite{arano-vaes} have vanishing $L^2$-Betti numbers.
\begin{corollary}
    \label{thm:arano-vaes-vanishing}
    Let $G$ be a second-countable locally compact totally disconnected non-discrete group, and $K<G$ a compact open subgroup.
    Then the category $\mathcal{C}_f(K<G)$ has vanishing $L^2$-Betti numbers.
\begin{proof}
    Recall from \thref{thm:av-normality} that $\mathcal{C}_f(K<G)$ contains a copy of $\Rep_f(K)$ as a normal subcategory.
    We also know from \cite[Proposition~2.7.12]{neshveyev-tuset} that $\Rep_f(K)$ is amenable.
    Moreover, since we required $G$ to be non-discrete, $K$ must be an infinite compact group, which means that $\Rep_f(K)$ has infinitely many distinct irreducibles.
    The result now follows from \thref{thm:infinite-amenable-subcategory-vt}.
\end{proof}
\end{corollary}

Finally, our results also imply that certain bicrossed products have vanishing $L^2$-Betti numbers.
\begin{corollary}
    Let $G$ be a second-countable locally compact totally disconnected non-discrete group, $K<G$ a compact open subgroup, and $\Gamma<G$ a discrete subgroup such that $G=K\Gamma$ and $K\cap\Gamma=\{e\}$.
    Let $\mathbb{H}$ be the bicrossed product compact quantum group associated with the matched pair $(K,\Gamma)$.
    Then $\hat{\mathbb{H}}$ has vanishing $L^2$-Betti numbers.
\begin{proof}
    Recall that the existence of such a matched pair decomposition forces $G$ to be unimodular, and the bicrossed product compact quantum group $\mathbb{H}$ will always be of Kac type \cite[Theorem~3.4]{fmp-bicrossed-products}.
    The vanishing of all $L^2$-Betti numbers of $\mathbb{H}$ can be viewed as a corollary of \thref{thm:vanishing-theorem} in two ways: either via \thref{thm:dqg-vanishing-theorem}, or by passing through \thref{thm:arano-vaes-vanishing}.
    The first argument proceeds by observing that $\mathbb{H}$ fits into an exact sequence $1\to\hat{\Gamma}\to \mathbb{H}\to K\to 1$ of compact quantum groups.
    Since $K$ is infinite and coamenable \cite[see][Proposition~2.7.12]{neshveyev-tuset}, the conditions of \thref{thm:dqg-vanishing-theorem} are then satisfied.
    Alternatively, one could appeal to \cite[Proposition~6.1]{vv-property-t} to get that $\Rep_f(\mathbb{H})$ is unitarily monoidally equivalent to $\mathcal{C}_f(K<G)$.
    The conclusion then follows after applying \thref{thm:arano-vaes-vanishing} and \cite[Theorem~4.2]{krvv}.
\end{proof}
\end{corollary}

\subsection{Proof of the vanishing theorem}
The proof of \thref{thm:vanishing-theorem} proceeds in two major steps, which we treat in separate subsections.
First, we apply the techniques from homological algebra recalled in \autoref{sec:homological-algebra} to prove a ``generic'' dimension-vanishing criterion for arbitrary inclusions $\mathcal{C}\subset\mathcal{D}$ of rigid $C^*$-tensor categories.
This part of the argument is entirely algebraic, and involves almost no dimension theory.
In the second half, we show that the conditions of this generic vanishing criterion are satisfied under the hypotheses of \thref{thm:vanishing-theorem}.
This part of the proof uses very little homological machinery; it relies mostly on a local version of the dimension scaling formula from \cite{krvv}.

\subsubsection{A homological bicomplex}
In this section, we associate a first quadrant bicomplex to any full inclusion of rigid $C^*$-tensor categories $\mathcal{C}\subset\mathcal{D}$ and any nondegenerate right $\mathcal{A}_{\mathcal{D}}$-module $V$.
This bicomplex is built up out of the usual bar resolution of $V$ combined with an $\tilde{\mathcal{A}}_{\mathcal{C}}$-relative version of the bar complex, which we introduce in the lemma below.
\begin{lemma}
\label{thm:relative-resolution}
Let $\mathcal{D}$ be a rigid $C^*$-tensor category and $\mathcal{C}$ a full tensor subcategory.
Put
\[
    P_{\mathcal{D}/\mathcal{C}}^q=\mathcal{A}_{\mathcal{D}}^{\otimes_{\tilde{\mathcal{A}}_{\mathcal{C}}}(q+2)}
\]
The system
\begin{equation}
    \label{eqn:relative-resolution-bare}
    0\ot \mathcal{A}_{\mathcal{D}} \stackrel{m}{\ot} P_{\mathcal{D}/\mathcal{C}}^0 \ot \cdots \ot P_{\mathcal{D}/\mathcal{C}}^\bullet
\end{equation}
with differentials defined by $\partial=\sum_{i=0}^q (-1)^i \partial^i$ is an exact sequence of projective left and right $\mathcal{A}_{\mathcal{D}}$-modules, where the face maps $\partial^i$ are given by
\[
    \partial^i(V_0\otimes_{\tilde{\mathcal{A}}_{\mathcal{C}}}\cdots\otimes_{\tilde{\mathcal{A}}_{\mathcal{C}}} V_q\otimes_{\tilde{\mathcal{A}}_{\mathcal{C}}} V_{q+1}) 
    = V_0\otimes_{\tilde{\mathcal{A}}_{\mathcal{C}}} \cdots\otimes_{\tilde{\mathcal{A}}_{\mathcal{C}}} V_i\cdot V_{i+1}\otimes_{\tilde{\mathcal{A}}_{\mathcal{C}}} \cdots \otimes_{\tilde{\mathcal{A}}_{\mathcal{C}}} V_q \otimes_{\tilde{\mathcal{A}}_{\mathcal{C}}} V_{q+1}\, .
\]
and $m:\mathcal{A}_{\mathcal{D}}\otimes_{\tilde{\mathcal{A}}_{\mathcal{C}}} \mathcal{A}_{\mathcal{D}}\to \mathcal{A}_{\mathcal{D}}$ is the multiplication map.
For any nondegenerate right $\mathcal{A}_{\mathcal{D}}$-module $V$ (resp.\@ left $\mathcal{A}_{\mathcal{D}}$-module $W$) the respective induced complexes
\begin{equation}
    \label{eqn:relative-resolution-tp}
    \begin{aligned}
    &0\ot V\ot V\otimes_{\tilde{\mathcal{A}}_{\mathcal{C}}} \mathcal{A}_{\mathcal{D}}
    \ot\cdots \ot
    V\otimes_{\tilde{\mathcal{A}}_{\mathcal{C}}} \mathcal{A}_{\mathcal{D}}^{\otimes_{\tilde{\mathcal{A}}_{\mathcal{C}}}\bullet}\\
    &0\ot W\ot \mathcal{A}_{\mathcal{D}}\otimes_{\tilde{\mathcal{A}}_{\mathcal{C}}}W
    \ot\cdots \ot
    \mathcal{A}_{\mathcal{D}}^{\otimes_{\tilde{\mathcal{A}}_{\mathcal{C}}}\bullet}\otimes_{\tilde{\mathcal{A}}_{\mathcal{C}}} W
    \end{aligned}
\end{equation}
are therefore exact.
\begin{proof}
    The projectivity follows from \cite[Proposition~3.12]{krvv}, as explained in \autoref{sec:cat-inclusions}.
    The standard argument shows that $\partial\partial=0$.
    To show exactness, we exhibit an explicit homotopy.
    For convenience, write $\mathcal{A}_{\mathcal{D}}=P^{-1}_{\mathcal{D}/\mathcal{C}}$ and put $m=\partial_{-1}$.
    For $q\geq -1$, we then set
    \[
        h(p_i\cdot V_0\otimes_{\tilde{\mathcal{A}}_{\mathcal{C}}}\cdots\otimes_{\tilde{\mathcal{A}}_{\mathcal{C}}} V_q\otimes_{\tilde{\mathcal{A}}_{\mathcal{C}}} V_{q+1})
        =p_i\otimes_{\mathcal{\tilde{A}}_{\mathcal{C}}}V_0\otimes_{\tilde{\mathcal{A}}_{\mathcal{C}}}\cdots\otimes_{\tilde{\mathcal{A}}_{\mathcal{C}}} V_q\otimes_{\tilde{\mathcal{A}}_{\mathcal{C}}} V_{q+1}\,.
    \]
    Clearly $\partial h + h\partial$ is the identity map on \eqref{eqn:relative-resolution-bare}, so the identity is null-homotopic.
    It follows that the complexes in \eqref{eqn:relative-resolution-tp} are also exact, e.g.\@ by interpreting \eqref{eqn:relative-resolution-bare} as a deleted projective resolution of the zero module.
\end{proof}
\end{lemma}
If $\mathcal{C}$ is the trivial subcategory generated by the unit object, we of course recover the usual bar complex of $\mathcal{D}$. In this situation, we simply write $P_{\mathcal{D}}^\bullet$ instead of $P_{\mathcal{D}/\mathcal{C}}^\bullet$.

\begin{sidewaysfigure}[p]
\includegraphics[width=0.85\paperheight]{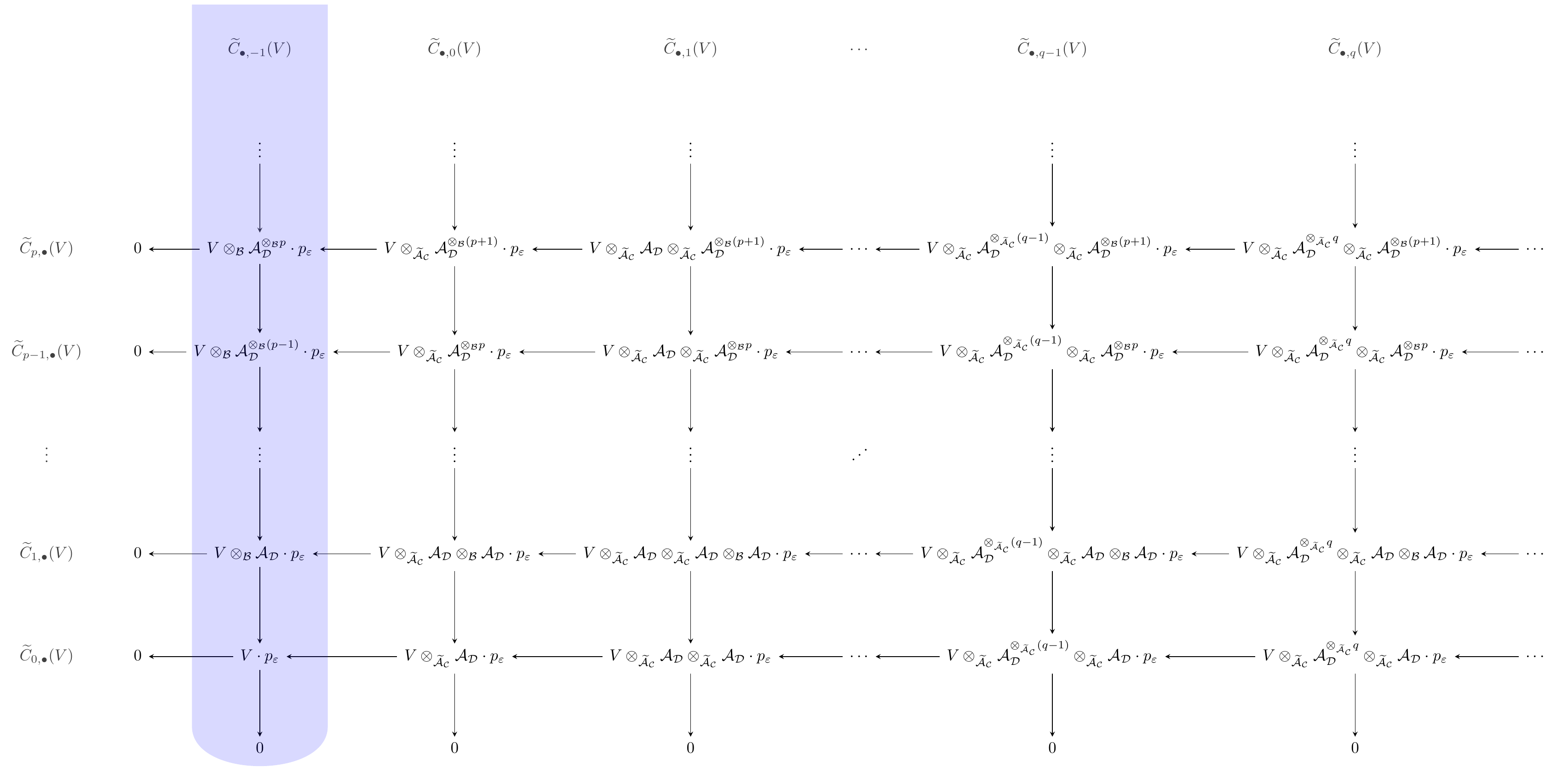}
\caption{%
    Diagram of the augmented bicomplex $\tilde{C}_{\bullet,\bullet}(V)$.
    To obtain $C_{\bullet,\bullet}(V)$ from $\tilde{C}_{\bullet,\bullet}(V)$, replace the shaded column $\tilde{C}_{\bullet,-1}(V)$ with zeroes.
}
\label{fig:augmented-bicomplex}
\end{sidewaysfigure}

\begin{definition}
Let $\mathcal{C}$ be a full tensor subcategory of a rigid $C^*$-tensor category $\mathcal{D}$ and $V$ a nondegenerate right $\mathcal{A}_{\mathcal{D}}$-module.
The first quadrant bicomplex $C_{\bullet,\bullet}(V)$ is defined as the $\mathcal{A}_{\mathcal{D}}$-tensor product complex of $V\otimes_{\mathcal{A}_{\mathcal{D}}}P_{\mathcal{D}/\mathcal{C}}^\bullet$ and $P_{\mathcal{D}}^\bullet\otimes_{\mathcal{A}_{\mathcal{D}}}\CC$, i.e.\@ 
\[
    C_{p,q}(V) \cong V\otimes_{\tilde{\mathcal{A}}_{\mathcal{C}}} \mathcal{A}_{\mathcal{D}}^{\otimes_{\tilde{\mathcal{A}}_{\mathcal{C}}} q} \otimes_{\tilde{\mathcal{A}}_{\mathcal{C}}} \mathcal{A}_{\mathcal{D}}^{\otimes_{\mathcal{B}} (p+1)}\cdot p_\epsilon
\] 
for $p,q\geq 0$. Here, $\mathcal{B}$ denotes the $*$-algebra defined in \eqref{eqn:b-definition}, i.e.\@ the $\ast$-subalgebra of $\mathcal{A}_{\mathcal{D}}$ spanned by the projections $p_i$, $i\in\Irr(\mathcal{D})$.

We also define an augmented version $\tilde{C}_{\bullet,\bullet}(V)$ of $C_{\bullet,\bullet}(V)$ by putting $\tilde{C}_{p,q}(V)=C_{p,q}(V)$ for $p,q\geq 0$ and $\tilde{C}_{p,-1}(V)=V\otimes_{\mathcal{B}} \mathcal{A}_{\mathcal{D}}^{\otimes_{\mathcal{B}} p}\cdot p_\epsilon$ for $p\geq 0$.
The extra differentials $\tilde{C}_{p,0}(V)\to \tilde{C}_{p,-1}(V)$ are induced by the module map $\sigma: V\otimes_{\tilde{\mathcal{A}}_{\mathcal{C}}} \mathcal{A}_{\mathcal{D}}\to V$.
The vertical differentials on $\tilde{C}_{p,-1}(V)$ are the usual ones for inhomogeneous chains.
See \autoref{fig:augmented-bicomplex} for a diagram of $\tilde{C}_{\bullet,\bullet}$.
\end{definition}

Before proceeding, we make a few elementary observations about the structure of $C_{\bullet,\bullet}(V)$ and $\tilde{C}_{\bullet,\bullet}(V)$.
First, note that the extra column $\tilde{C}_{\bullet,-1}(V)$ (drawn with shading in \autoref{fig:augmented-bicomplex}) is precisely the bar resolution of $V$ as a right $\tilde{\mathcal{A}}_{\mathcal{D}}$-module.
Therefore we have natural identifications of homology spaces
\begin{equation}
\label{eqn:augmented-column-homology}
    H_*(C_{\bullet,-1}(V))\cong \Tor_*^{\mathcal{A}_{\mathcal{D}}}(V,\CC) \cong H_*(\mathcal{D};V)\,.
\end{equation}
The homology of the other columns can be characterised as follows.
\begin{lemma}
\label{thm:column-homology}
Let $\mathcal{D}$ be a rigid $C^*$-tensor category, $\mathcal{C}$ a full tensor subcategory and $V$ a right $\mathcal{A}_{\mathcal{D}}$-module.
For all $n,q\geq 0$, we have that
\[
    H_n(C_{\bullet, q}(V)) \cong \Tor_n^{\tilde{\mathcal{A}}_{\mathcal{C}}}\left(V\otimes_{\tilde{\mathcal{A}}_{\mathcal{C}}}\mathcal{A}_{\mathcal{D}}^{\otimes_{\tilde{\mathcal{A}}_{\mathcal{C}}} q}, \CC\right).
\]
\begin{proof}
    Recall that there is a resolution of the counit $\counit: \mathcal{A}_{\mathcal{D}}\to \CC$ by projective left $\mathcal{A}_{\mathcal{D}}$-modules given by
    \begin{equation}
    \label{eqn:rho-resolution}
        0\ot \CC\stackrel{\counit}{\ot} \mathcal{A}_{\mathcal{D}}\cdot p_\epsilon\ot \cdots \ot \mathcal{A}_{\mathcal{D}}^{\otimes_{\mathcal{B}} \bullet}\cdot p_\epsilon\,.
    \end{equation}
    By \cite[Proposition~3.12]{krvv}, $\mathcal{A}_{\mathcal{D}}$ is projective as a left $\tilde{\mathcal{A}}_{\mathcal{C}}$-module. 
    It follows that \eqref{eqn:rho-resolution} is still a projective resolution of $\CC$ as a left $\tilde{\mathcal{A}}_{\mathcal{C}}$-module.
    Since $C_{\bullet,q}(V)$ is the result of applying $V\otimes_{\tilde{\mathcal{A}}_{\mathcal{C}}}\mathcal{A}_{\mathcal{D}}^{\otimes_{\tilde{\mathcal{A}}_{\mathcal{C}}} q}\otimes_{\tilde{\mathcal{A}}_{\mathcal{C}}}-$ to the deleted version of this resolution, the claim follows (see also \autoref{fig:augmented-bicomplex}).
\end{proof}
\end{lemma}

\begin{lemma}
\label{thm:augmented-row-homology}
Let $\mathcal{D}$ be a rigid $C^*$-tensor category, $\mathcal{C}$ a full tensor subcategory and $V$ a right $\mathcal{A}_{\mathcal{D}}$-module.
Then there are natural isomorphisms $H_*(\Tot(C_{\bullet,\bullet}(V))) \cong H_*(\mathcal{D}; V)$ in homology.
\begin{proof}
    We claim that the augmented complex $\tilde{C}_{\bullet,\bullet}(V)$ has exact rows.
    This is a consequence of the exactness of 
    \[
        0\ot V\ot V\otimes_{\tilde{\mathcal{A}}_{\mathcal{C}}} \mathcal{A}_{\mathcal{D}}
        \ot\cdots \ot
        V\otimes_{\tilde{\mathcal{A}}_{\mathcal{C}}} \mathcal{A}_{\mathcal{D}}^{\otimes_{\tilde{\mathcal{A}}_{\mathcal{C}}}\bullet}
    \]
    ensured by \thref{thm:relative-resolution}.
    Indeed, $\tilde{C}_{p,\bullet}(V)$ is the result of tensoring this complex with the projective (in particular flat) $\mathcal{A}_{\mathcal{D}}$-module $\mathcal{A}_{\mathcal{D}}^{\otimes_{\mathcal{B}} (p+1)}\cdot p_\epsilon$.
    Note that the augmented term $\tilde{C}_{p,-1}(V)$ is crucial to obtain exactness in $\tilde{C}_{p,0}(V)$.
    By \thref{thm:bicomplex-aal}, it follows that the total complex of $\tilde{C}_{\bullet,\bullet}(V)$ is exact.
    Applying \thref{thm:augmented-bicomplex-homology} to $\tilde{C}_{\bullet,\bullet}(V)$ (shifted by one horizontal degree), we get identifications
    \[
        H_*(\Tot(C_{\bullet,\bullet}(V))) \cong H_*(\tilde{C}_{\bullet,-1}(V)) \cong H_*(\mathcal{D}; V)
    \]
    so the claim is proven.
\end{proof}
\end{lemma}

Since $C_{\bullet,\bullet}(V)$ is functorial in $V$ and the isomorphisms in \thref{thm:column-homology} and \thref{thm:augmented-row-homology} are natural in $V$, any left $R$-module structure will carry over to the terms appearing in $C_{\bullet,\bullet}(V)$ in the obvious way, and factor through all homology computations.
We apply this to get the following generic vanishing criterion.  
\begin{proposition}
\label{thm:generic-vanishing-theorem}
Let $\mathcal{D}$ be a rigid $C^*$-tensor category, and $\mathcal{C}$ a full tensor subcategory.
Additionally, let $M$ be a von Neumann algebra with n.s.f.\@ trace $\tau$, and $V$ a $M$-$\mathcal{A}_{\mathcal{D}}$-bimodule such that $V$ is nondegenerate as an $\mathcal{A}_{\mathcal{D}}$-module and locally finite as an $M$-module.

Given $N\geq 0$ such that
\[
    \forall q\geq 0, n\leq N:\qquad \dim_M \Tor_n^{\tilde{\mathcal{A}}_{\mathcal{C}}}\left(V\otimes_{\tilde{\mathcal{A}}_{\mathcal{C}}}\mathcal{A}_{\mathcal{D}}^{\otimes_{\tilde{\mathcal{A}}_{\mathcal{C}}} q},\CC\right) = 0\,,
\]
we also have that $\dim_M H_n(\mathcal{D}; V) = 0$ for all $n\leq N$.
\begin{proof}
    By \thref{thm:column-homology} the columns $C_{\bullet,q}$ are dimension exact at $C_{n,q}$ for all $n\leq N$ and $q\geq 0$, so \thref{thm:bicomplex-dimension-aal} tells us that
    \[
        \dim_M H_n(\mathrm{Tot}(C_{\bullet,\bullet}(V)))=0.
    \]
    for all $n\leq N$.
    But we already knew that $H_n(\mathrm{Tot}(C_{\bullet,\bullet}(V)))\cong H_n(\mathcal{D};V)$ from \thref{thm:augmented-row-homology}, so the conclusion follows.
\end{proof}
\end{proposition}

\subsubsection{Inductively establishing dimension-exactness in the total complex}
In order to prove \thref{thm:vanishing-theorem}, we will show that the conditions of \thref{thm:generic-vanishing-theorem} are satisfied for $V=L^2(\mathcal{A}_{\mathcal{D}})^0$ and $M=\mathcal{M}_{\mathcal{D}}$ whenever $\mathcal{C}$ is almost normal in $\mathcal{D}$, in addition to having vanishing $L^2$-Betti numbers up to a certain degree.
Before doing so, we first discuss the relationship between the various algebras associated with $\mathcal{C}\subset\mathcal{D}$ when $\mathcal{C}$ is almost normal in $\mathcal{D}$.

For all double cosets $\mathcal{C}\alpha\mathcal{C}$, the projection $e_{\mathcal{C}\alpha\mathcal{C}}$ onto the closed linear span of $\tilde{\mathcal{A}}_{\mathcal{C}\alpha\mathcal{C}}\subset L^2(\mathcal{A}_{\mathcal{D}},\tau)$ satisfies
\[
    e_{\mathcal{C}\alpha\mathcal{C}} 
    = \sum_{\mathcal{C}\beta\subset\mathcal{C}\alpha\mathcal{C}} e_{\mathcal{C}\beta} 
    = \sum_{\gamma\mathcal{C}\subset\mathcal{C}\alpha\mathcal{C}} e_{\gamma\mathcal{C}}\,.
\]
Observe that the sums are finite, precisely because $\mathcal{C}\subset\mathcal{D}$ is almost normal.
It follows that the space $\mathcal{Z}_{\mathcal{C}\alpha\mathcal{C}}$ defined in \eqref{eqn:z-module-def} admits the following finite direct sum decompositions:
\[
    \mathcal{Z}_{\mathcal{C}\alpha\mathcal{C}} 
    = \bigoplus_{\mathcal{C}\beta\subset\mathcal{C}\alpha\mathcal{C}} \mathcal{Z}_{\mathcal{C}\beta} 
    = \bigoplus_{\gamma\mathcal{C}\subset\mathcal{C}\alpha\mathcal{C}} \mathcal{Z}_{\gamma\mathcal{C}} \,.
\]
By \thref{thm:z-orbit-module}, we get that $\mathcal{Z}_{\mathcal{C}\alpha\mathcal{C}}$ is a flat $\tilde{\mathcal{M}}_{\mathcal{C}}$-$\tilde{\mathcal{M}}_{\mathcal{C}}$-subbimodule of $L^2(\mathcal{A}_{\mathcal{D}},\tau)$, and that $p_F\cdot\mathcal{Z}_{\mathcal{C}\alpha\mathcal{C}}$ (resp.\@ $\mathcal{Z}_{\mathcal{C}\alpha\mathcal{C}}\cdot p_F$) is projective as a right (resp.\@ left) $\tilde{\mathcal{M}}_{\mathcal{C}}$-module for all finite $F\subset\Irr(\mathcal{D})$.
\begin{lemma}
    \label{thm:coset-tensor-iso}
    Let $\mathcal{C}\subset\mathcal{D}$ be an inclusion of rigid C*-tensor categories, and suppose that $\mathcal{C}$ is almost normal in $\mathcal{D}$.
    With the notation introduced in \autoref{sec:cat-inclusions}, the natural multiplication maps
    \begin{align*}
        \tilde{\mathcal{A}}_{\mathcal{C}\alpha\mathcal{C}}\otimes_{\tilde{\mathcal{A}}_{\mathcal{C}}} \lsup{0}{L^2(\tilde{\mathcal{A}}_{\mathcal{C}},\tau)} &\to e_{\mathcal{C}\alpha\mathcal{C}}\lsup{0}{L^2(\mathcal{A}_{\mathcal{D}},\tau)}\\
        L^2(\tilde{\mathcal{A}}_{\mathcal{C}},\tau)^0 \otimes_{\tilde{\mathcal{A}}_{\mathcal{C}}} \tilde{\mathcal{A}}_{\mathcal{C}\alpha\mathcal{C}}&\to e_{\mathcal{C}\alpha\mathcal{C}}L^2(\mathcal{A}_{\mathcal{D}},\tau)^0\\
        p_F\cdot\mathcal{Z}_{\mathcal{C}\alpha\mathcal{C}}\otimes_{\tilde{\mathcal{M}}_{\mathcal{C}}} L^2(\tilde{\mathcal{A}}_{\mathcal{C}},\tau) &\to p_F\cdot e_{\mathcal{C}\alpha\mathcal{C}} L^2(\mathcal{A}_{\mathcal{D}},\tau)\\ 
        L^2(\tilde{\mathcal{A}}_{\mathcal{C}},\tau)\otimes_{\tilde{\mathcal{M}}_{\mathcal{C}}} \mathcal{Z}_{\mathcal{C}\alpha\mathcal{C}}\cdot p_F &\to e_{\mathcal{C}\alpha\mathcal{C}} L^2(\mathcal{A}_{\mathcal{D}},\tau)\cdot p_F
    \end{align*}
    are isomorphisms for all $\alpha\in\Irr(\mathcal{D})$ and all finite $F\subset\Irr(\mathcal{D})$.
\begin{proof}
    Note that $p_F\cdot \mathcal{Z}_{\mathcal{C}\alpha\mathcal{C}}$ and $\mathcal{Z}_{\mathcal{C}\alpha\mathcal{C}}\cdot p_F$ are subsets of $\mathcal{M}_{\mathcal{D}}$ for all finite $F\subset\Irr(\mathcal{D})$, so the third and fourth multiplication maps are well-defined.
    In exactly the same way as the proof of \thref{thm:z-orbit-module}, one can write down explicit inverses for all these multiplication maps.
\end{proof}
\end{lemma}

\begin{lemma}
\label{thm:local-scaling-formula}
    Let $\mathcal{C}\subset\mathcal{D}$ be an almost-normal inclusion.
    Fix $\alpha\in\Irr(\mathcal{D})$, a locally finite left $\tilde{\mathcal{M}}_{\mathcal{C}}$-module $H$ and a locally finite right $\tilde{\mathcal{M}}_{\mathcal{C}}$-module $K$.
    Then
    \begin{align*}
        \dim_{\tilde{\mathcal{M}}_{\mathcal{C}}-}\left(\mathcal{Z}_{\mathcal{C}\alpha\mathcal{C}}\otimes_{\tilde{\mathcal{M}}_{\mathcal{C}}} H\right)
        &=\sum_{\mathcal{C}\beta\subset\mathcal{C}\alpha\mathcal{C}} \frac{d(\beta)^2}{d([\beta\bar{\beta}]_{\mathcal{C}})} \dim_{\tilde{\mathcal{M}}_{\mathcal{C}}-}(H)\\
        \dim_{-\tilde{\mathcal{M}}_{\mathcal{C}}}\left(K\otimes_{\tilde{\mathcal{M}}_{\mathcal{C}}}\mathcal{Z}_{\mathcal{C}\alpha\mathcal{C}}\right)
        &=\sum_{\beta\mathcal{C}\subset\mathcal{C}\alpha\mathcal{C}} \frac{d(\beta)^2}{d([\bar{\beta}\beta]_{\mathcal{C}})} \dim_{-\tilde{\mathcal{M}}_{\mathcal{C}}}(K)
    \end{align*}
\begin{proof}
    We will prove only the version for right modules for notational convenience, using the same methods as \cite{krvv}.
    Fix $\alpha\in\Irr(\mathcal{D})$ and representatives $\alpha_1,\ldots,\alpha_\kappa$ of orbits $\alpha_1\mathcal{C},\ldots,\alpha_\kappa\mathcal{C}$ covering $\mathcal{C}\alpha\mathcal{C}$.
    By \thref{thm:generic-scaling-formula} and the fact that $\mathcal{Z}_{\mathcal{C}\alpha\mathcal{C}}$ is flat over $\tilde{\mathcal{M}}_{\mathcal{C}}$, we only have to check the statement for locally finite finitely generated projective modules.
    Define a Hilbert space 
    \[
        \mathcal{H}=\bigoplus_{k=1}^\kappa\bigoplus_{i,j\in\Irr(\mathcal{D})} (i\alpha_k,\alpha_k j)
    \]
    and a right $\tilde{\mathcal{M}}_{\mathcal{C}}$-linear map into the Hilbert space tensor product
    \begin{align*}
        &U:e_{\mathcal{C}\alpha\mathcal{C}}L^2(\mathcal{A}_{\mathcal{D}},\tau)\to \mathcal{H}\htimes L^2(\tilde{\mathcal{A}}_{\mathcal{C}},\tau): \\
        &\qquad \xi\mapsto \sum_{k=1}^\kappa\sum_{\substack{i,j\in\Irr(\mathcal{D})\\ W\in\onb(i\alpha_k,\alpha_k j)}} \sqrt{\frac{d(j)d(\alpha_k)}{d([\bar{\alpha_k}\alpha_k]_{\mathcal{C}})}}W\otimes e_{\mathcal{C}}(W^\#\cdot\xi)\,.
    \end{align*}
    By \cite[Lemma~3.11]{krvv}, $U$ is isometric. Put $q=UU^*$, and observe that
    \[
        U^*: \mathcal{H}\otimes L^2(\tilde{\mathcal{A}}_{\mathcal{C}},\tau) \to e_{\mathcal{C}\alpha\mathcal{C}}L^2(\mathcal{A}_{\mathcal{D}},\tau): 
        (i\alpha_k,\alpha_k j) \otimes L^2(\tilde{\mathcal{A}}_{\mathcal{C}},\tau) \ni W\otimes \xi\mapsto \sqrt{\frac{d(j)d(\alpha_k)}{d([\bar{\alpha_k}\alpha_k]_{\mathcal{C}})}} W \cdot \xi\,.
    \] 
    The composition $UxU^*$ yields a well-defined operator on $\mathcal{H}\otimes L^2(\tilde{\mathcal{A}}_{\mathcal{C}},\tau)$ for all $x\in\tilde{\mathcal{M}}_{\mathcal{C}}$.
    Since this is a composition of right $\tilde{\mathcal{M}}_{\mathcal{C}}$-linear operators, we find that $UxU^*$ lies in the commutant of the right diagonal action of $\tilde{\mathcal{M}}_{\mathcal{C}}$.
    Hence, this gives rise to a unital $\ast$-homomorphism
    \[
        \Xi:\tilde{\mathcal{M}}_{\mathcal{C}}\to q(\mathcal{B}(\mathcal{H})\,\overline{\otimes}\,\tilde{\mathcal{M}}_{\mathcal{C}})q: x\mapsto UxU^*\,.
    \]
    Moreover, for $x\in \left(\tilde{\mathcal{M}}_{\mathcal{C}}\right)_+$
    \begin{align*}
        (\Tr_{\mathcal{B}(\mathcal{H})}\otimes\tau)(UxU^*)
        &=\sum_{k=1}^\kappa\sum_{\substack{i,j\in\Irr(\mathcal{D})\\ W\in\onb(i\alpha_k,\alpha_k j)}} \frac{d(j)d(\alpha_k)}{d([\bar{\alpha_k}\alpha_k]_{\mathcal{C}})}\tau(W^\#\cdot x\cdot W)\\
        &=\sum_{k=1}^\kappa\frac{d(\alpha_k)^2}{d([\bar{\alpha_k}\alpha_k]_{\mathcal{C}})}\tau(x)\,.
        \numberthis\label{eqn:local-scaling-formula-u-conj-trace}
    \end{align*}
    Let $P\in M_n(\CC)\otimes \tilde{\mathcal{M}}_{\mathcal{C}}$ be a projection of finite trace.
    We embed $P(\CC^n\otimes \mathcal{Z}_{\mathcal{C}\alpha\mathcal{C}})\hookrightarrow (\id\otimes \Xi)(P)(\CC^n\otimes \mathcal{H}\htimes L^2(\tilde{\mathcal{A}}_{\mathcal{C}},\tau))$
    by sending $\xi\otimes \eta$ to $\xi\otimes U(\eta)$.
    This is well-defined because the coefficients of $P$ lie in $\tilde{\mathcal{M}}_{\mathcal{C}}$.
    It follows that
    \begin{equation}
    \label{eqn:local-scaling-formula-upper}
        \dim_{-\tilde{\mathcal{M}}_{\mathcal{C}}} P(\CC^n\otimes \mathcal{Z}_{\mathcal{C}\alpha\mathcal{C}}) 
        \leq (\Tr_n\otimes\Tr_{\mathcal{B}(\mathcal{H})}\otimes\tau)((\id\otimes\Xi)(P))\,.
    \end{equation}
    Conversely, given $V\in (i\alpha_k,\alpha_k j)$, $\eta\in\CC^n$, $x\in \tilde{\mathcal{M}}_{\mathcal{C}}$, we also know that
    \begin{equation}
    \label{eqn:local-scaling-formula-u-star-ampl}
    \begin{aligned}
        (\I\otimes U^*)&\left((\id\otimes\Xi)(P)(\eta\otimes V\otimes x)\right) \\
        &=P(\eta\otimes U^*(V\otimes x))\\
        &= \sqrt{\frac{d(j)d(\alpha_k)}{d([\bar{\alpha_k}\alpha_k]_{\mathcal{C}})}} P(\eta\otimes V\cdot x)\qquad \in P(\CC^n\otimes\mathcal{Z}_{\mathcal{C}\alpha\mathcal{C}})\,.
    \end{aligned}
    \end{equation}
    Define
    \[
        \mathcal{H}_F=\bigoplus_{k=1}^\kappa\bigoplus_{\substack{i\in F\\ j\in\Irr(\mathcal{D})}} (i\alpha_k,\alpha_k j)\subset\mathcal{H}
    \]
    for finite $F\subset\Irr(\mathcal{D})$.
    Denote the projection of $\mathcal{H}$ onto $\mathcal{H}_F$ by $q_F$.
    By \eqref{eqn:local-scaling-formula-u-star-ampl}, $(\I\otimes U^*)$ restricts to an injective right-linear map
    \begin{equation}
    \label{eqn:local-scaling-formula-module-embedding}
        (\id\otimes\Xi)(P)(\CC^n\otimes\mathcal{H}_F\otimes \tilde{\mathcal{M}}_{\mathcal{C}})
        \to P(\CC^n\otimes \mathcal{Z}_{\mathcal{C}\alpha\mathcal{C}})\,.
    \end{equation}
    Since $F$ is finite, the left-hand side is a finitely generated locally finite projective module.
    Indeed, it is isomorphic to $P_0(\CC^n\otimes\mathcal{H}_F\otimes \tilde{\mathcal{M}}_{\mathcal{C}})$, where $P_0$ is the right support projection of $(\id\otimes\Xi)(P)(\I\otimes q_F\otimes\I)$.
    The trace of $P_0$ is equal to that of the left support of $(\id\otimes\Xi)(P)(\I\otimes q_F\otimes\I)$, which increases strongly to $(\id\otimes\Xi)(P)$ as $F\nearrow\Irr(\mathcal{D})$.
    It follows that the trace of $P_0$ converges to the trace of $(\id\otimes\Xi)(P)$.
    In summary, the dimension of the module on the left of \eqref{eqn:local-scaling-formula-module-embedding} increases towards
    \[
        (\Tr_n\otimes\Tr_{\mathcal{B}(\mathcal{H})}\otimes\tau)\left((\id\otimes\Xi)(P)\right) 
        = \sum_{k=1}^\kappa \frac{d(\alpha_k)^2}{d([\bar{\alpha_k}\alpha_k]_{\mathcal{C}})} (\Tr\otimes\tau)(P)
    \]
    via another application of \eqref{eqn:local-scaling-formula-u-conj-trace}.
    We conclude that we can embed finitely generated locally finite projective modules of dimension arbitrarily close to the right-hand side of \eqref{eqn:local-scaling-formula-upper} into $P(\CC^n\otimes \mathcal{Z}_{\mathcal{C}\alpha\mathcal{C}})$, yielding the reverse inequality.

    An appeal to \eqref{eqn:local-scaling-formula-u-conj-trace} now yields
    \begin{align*}
        \dim_{-\tilde{\mathcal{M}}_{\mathcal{C}}} P(\CC^n\otimes \mathcal{Z}_{\mathcal{C}\alpha\mathcal{C}}) 
        &= (\Tr_n\otimes\Tr_{\mathcal{B}(\mathcal{H})}\otimes\tau)((\id\otimes\Xi)(P))\\
        &=\sum_{i=1}^n (\Tr_{\mathcal{B}(\mathcal{H})}\otimes\tau)(UP_{ii}U^*)\\
        &=\sum_{k=1}^\kappa\frac{d(\alpha_k)^2}{d([\bar{\alpha_k}\alpha_k]_{\mathcal{C}})} \sum_{i=1}^n\tau(P_{ii})\\
        &=\sum_{k=1}^\kappa\frac{d(\alpha_k)^2}{d([\bar{\alpha_k}\alpha_k]_{\mathcal{C}})} (\Tr\otimes\tau)(P)\\
        &=\sum_{k=1}^\kappa\frac{d(\alpha_k)^2}{d([\bar{\alpha_k}\alpha_k]_{\mathcal{C}})}  \dim_{-\tilde{\mathcal{M}}_{\mathcal{C}}} P(\CC^n\otimes \tilde{\mathcal{M}}_{\mathcal{C}})\,.
    \end{align*}
\end{proof}
\end{lemma}
\begin{remark}
    Note that \cite[Lemma~3.2]{krvv} already guarantees that the left- and right scaling constants appearing in \thref{thm:local-scaling-formula} do not depend on the choice of coset representatives, as one would expect.
    However, the left- and right scaling constants need not coincide in general, even when $\mathcal{C}\subset\mathcal{D}$ is normal.
    Given a locally compact totally disconnected group $G$ with a compact open subgroup $K$, the normal inclusion $\Rep_f(K)\subset\mathcal{C}_f(K<G)$ provides an interesting class of examples.
    In this situation, the left- and right scaling constants agree for all cosets if and only if $G$ is unimodular.
    Indeed, given $\mathcal{H}\in\mathcal{C}_f(K<G)$, 
    \begin{align*}
        \dim([\bar{\mathcal{H}}\mathcal{H}]_{\Rep_f(K)}) &= \sum_{gK\in G/K} \dim(\chi_{gK}\cdot\mathcal{H}\cdot \chi_K)^2\,,\quad\text{and}\\
        \dim([\mathcal{H}\bar{\mathcal{H}}]_{\Rep_f(K)}) &= \sum_{gK\in G/K} \dim(\chi_{K}\cdot\mathcal{H}\cdot \chi_{gK})^2\,.
    \end{align*}
    Plugging in $\mathcal{H}=\mathcal{K}_x$ as defined in \eqref{eqn:av-cat-coset-canon-rep}, we recover 
    \[
        \dim([\bar{\mathcal{K}_x}\mathcal{K}_x]_{\Rep_f(K)})=[K:K\cap xKx^{-1}]
        \qquad\text{and}\qquad
        \dim([\mathcal{K}_x\bar{\mathcal{K}_x}]_{\Rep_f(K)})=[K:K\cap x^{-1}Kx]\,.
    \]
    These quantities are equal for all $x\in G$ if and only if $G$ is unimodular.
\end{remark}

\begin{remark}
    \label{rem:markov-bimodules}
    Observe that the structure of the proof of \thref{thm:local-scaling-formula} is globally the same as the argument for Markov inclusions used in the proof of \cite[Proposition~3.5]{krvv}.
    In fact, both statements can be unified by working with Markov \textit{bimodules} instead.
    Given two von Neumann algebras $M$, $N$ equipped with n.s.f.\@ traces and a Hilbert $M$-$N$-bimodule $\mathcal{H}$, we have access to two tracial weights on $M$: the original trace $\tau_M$, and the trace $\Tr^r_{\mathcal{H}}$ on $(N^\op)'\subset\mathcal{B}(\mathcal{H})$.
    We say that $\mathcal{H}$ is right $\lambda$-Markov\footnote{Observe that this is automatic if $M$ and $N$ are $\mathrm{II}_1$ factors and $\mathcal{H}$ is a bifinite $M$-$N$-bimodule.} if $\tau_M=\lambda\Tr_{\mathcal{H}}^r$ for some $\lambda\in (0,+\infty)$.
    Left $\lambda$-Markov bimodules are defined analogously.
    In the terminology of \cite{krvv}, an inclusion $N\subset M$ is $\lambda$-Markov if and only if $L^2(M,\tau_M)$ is right $\lambda$-Markov as an $M$-$N$-bimodule.

    Let $K$ be an (algebraic) $M$-$N$-subbimodule of $\mathcal{H}$ that contains a Pimsner--Popa basis of $\mathcal{H}$ as an $N$-module.
    If the tensor product functor $-\otimes_M K$ is exact on locally finite modules, then the general scaling formula $\dim_{-N} H\otimes_M K = \lambda^{-1}\dim_{-M} H$ holds for all locally finite right $M$-modules $H$.
    The proof proceeds along the same lines as \thref{thm:local-scaling-formula} and \cite[Proposition~3.5]{krvv}.
\end{remark}

We are now ready to prove \thref{thm:vanishing-theorem}.
\begin{proof}[Proof of \thref{thm:vanishing-theorem}]
    By \thref{thm:generic-vanishing-theorem}, it suffices to verify that
    \[
        \dim_{\mathcal{M}_{\mathcal{D}}} \Tor_n^{\tilde{\mathcal{A}}_{\mathcal{C}}}\left(L^2(\mathcal{A}_{\mathcal{D}},\tau)^0\otimes_{\tilde{\mathcal{A}}_{\mathcal{C}}}\mathcal{A}_{\mathcal{D}}^{\otimes_{\tilde{\mathcal{A}}_{\mathcal{C}}} q},\CC\right) = 0
    \]
    for all $q\geq 0$, $n\leq N$.
    Note that the above is equivalent to
    \begin{equation}
    \label{eqn:vanishing-theorem-goal}
        \dim_{\tilde{\mathcal{M}}_{\mathcal{C}}} \Tor_n^{\tilde{\mathcal{A}}_{\mathcal{C}}}\left(L^2(\tilde{\mathcal{A}}_{\mathcal{C}},\tau)^0\otimes_{\tilde{\mathcal{A}}_{\mathcal{C}}}\mathcal{A}_{\mathcal{D}}^{\otimes_{\tilde{\mathcal{A}}_{\mathcal{C}}} q},\CC\right) = 0
    \end{equation}
    for all $q\geq 0$, $n\leq N$. 
    Indeed, this follows by observing that $\mathcal{M}_{\mathcal{D}}\otimes_{\tilde{\mathcal{M}}_{\mathcal{C}}} -$ is exact and dimension-preserving in combination with the fact that $\mathcal{M}_{\mathcal{D}} \otimes_{\tilde{\mathcal{M}}_{\mathcal{C}}} L^2(\tilde{\mathcal{A}}_{\mathcal{C}},\tau)^0$ is dimension-isomorphic to $L^2(\mathcal{A}_{\mathcal{D}},\tau)^0$ (cfr.\@ \thref{thm:vna-hilbert-dim-iso}).

    We will prove \eqref{eqn:vanishing-theorem-goal} by induction on $q$.
    The base case $q=0$ is the hypothesis of the theorem stating that $\beta^{(2)}_n(\mathcal{C})=0$ for all $n\leq N$.
    For the induction step, observe that \thref{thm:coset-tensor-iso} provides us with isomorphisms of $p_F\cdot\tilde{\mathcal{M}}_{\mathcal{C}}\cdot p_F$-$\mathcal{A}_{\mathcal{D}}$-bimodules
    \begin{align*}
        p_F\cdot L^2(\tilde{\mathcal{A}}_{\mathcal{\mathcal{C}}},\tau)^0\otimes_{\tilde{\mathcal{A}}_{\mathcal{C}}}\mathcal{A}_{\mathcal{D}}^{\otimes_{\tilde{\mathcal{A}}_{\mathcal{C}}} (q+1)}
        &\cong \bigoplus_{\mathcal{C}\alpha\mathcal{C}\in\mathcal{C}\backslash \mathcal{D}/\mathcal{C}} p_F\cdot L^2(\tilde{\mathcal{A}}_{\mathcal{C}},\tau)^0\otimes_{\tilde{\mathcal{A}}_{\mathcal{C}}}\tilde{\mathcal{A}}_{\mathcal{C}\alpha\mathcal{C}}\otimes_{\tilde{\mathcal{A}}_{\mathcal{C}}}\mathcal{A}_{\mathcal{D}}^{\otimes_{\tilde{\mathcal{A}}_{\mathcal{C}}} q}\\
        &\cong \bigoplus_{\mathcal{C}\alpha\mathcal{C}\in\mathcal{C}\backslash \mathcal{D}/\mathcal{C}} p_F\cdot e_{\mathcal{C}\alpha\mathcal{C}}L^2(\mathcal{A}_{\mathcal{\mathcal{D}}},\tau)^0\otimes_{\tilde{\mathcal{A}}_{\mathcal{C}}}\mathcal{A}_{\mathcal{D}}^{\otimes_{\tilde{\mathcal{A}}_{\mathcal{C}}} q}\\
        &\cong \bigoplus_{\mathcal{C}\alpha\mathcal{C}\in\mathcal{C}\backslash \mathcal{D}/\mathcal{C}} p_F\cdot \mathcal{Z}_{\mathcal{C}\alpha\mathcal{C}} \otimes_{\tilde{\mathcal{M}}_{\mathcal{C}}} L^2(\tilde{\mathcal{A}}_{\mathcal{\mathcal{C}}},\tau)^0\otimes_{\tilde{\mathcal{A}}_{\mathcal{C}}}\mathcal{A}_{\mathcal{D}}^{\otimes_{\tilde{\mathcal{A}}_{\mathcal{C}}} q}
    \end{align*}
    for any finite subset $F\subset\Irr(\mathcal{D})$.
    Applying $H_n(\mathcal{C};-)$ and using the fact that $p_F\cdot\mathcal{Z}_{\mathcal{C}\alpha\mathcal{C}}$ is projective (and hence flat) as a right $\tilde{\mathcal{M}}_{\mathcal{C}}$-module, we recover isomorphisms of left $p_F\cdot\tilde{\mathcal{M}}_{\mathcal{C}}\cdot p_F$-modules
    \begin{align*}
        p_F\cdot H_n&\left(\mathcal{C}; L^2(\tilde{\mathcal{A}}_{\mathcal{\mathcal{C}}})^0\otimes_{\tilde{\mathcal{A}}_{\mathcal{C}}}\mathcal{A}_{\mathcal{D}}^{\otimes_{\tilde{\mathcal{A}}_{\mathcal{C}}} (q+1)}\right)\\
        &\cong \bigoplus_{\mathcal{C}\alpha\mathcal{C}\in\mathcal{C}\backslash \mathcal{D}/\mathcal{C}} H_n\left(\mathcal{C};p_F\cdot\mathcal{Z}_{\mathcal{C}\alpha\mathcal{C}} \otimes_{\tilde{\mathcal{M}}_{\mathcal{C}}} L^2(\tilde{\mathcal{A}}_{\mathcal{\mathcal{C}}})^0\otimes_{\tilde{\mathcal{A}}_{\mathcal{C}}}\mathcal{A}_{\mathcal{D}}^{\otimes_{\tilde{\mathcal{A}}_{\mathcal{C}}} q}\right) \\
        &\cong \bigoplus_{\mathcal{C}\alpha\mathcal{C}\in\mathcal{C}\backslash \mathcal{D}/\mathcal{C}} p_F\cdot \mathcal{Z}_{\mathcal{C}\alpha\mathcal{C}} \otimes_{\tilde{\mathcal{M}}_{\mathcal{C}}} H_n\left(\mathcal{C};L^2(\tilde{\mathcal{A}}_{\mathcal{C}})^0\otimes_{\tilde{\mathcal{A}}_{\mathcal{C}}} \mathcal{A}_{\mathcal{D}}^{\otimes_{\tilde{\mathcal{A}}_{\mathcal{C}}} q}\right)\,.
    \end{align*}
    Letting $F$ increase to $\Irr(\mathcal{D})$ and applying \cite[Lemma~A.16]{kpv13} leaves us with
    \begin{align*}
        \dim_{\tilde{\mathcal{M}}_{\mathcal{C}}} H_n&\left(\mathcal{C}; L^2(\tilde{\mathcal{A}}_{\mathcal{\mathcal{C}}})^0\otimes_{\tilde{\mathcal{A}}_{\mathcal{C}}}\mathcal{A}_{\mathcal{D}}^{\otimes_{\tilde{\mathcal{A}}_{\mathcal{C}}} (q+1)}\right)\\
        &=\sum_{\mathcal{C}\alpha\mathcal{C}\in \mathcal{C}\backslash \mathcal{D}/\mathcal{C}} \dim_{\tilde{\mathcal{M}}_{\mathcal{C}}}\left[\mathcal{Z}_{\mathcal{C}\alpha\mathcal{C}} \otimes_{\tilde{\mathcal{M}}_{\mathcal{C}}} H_n\left(\mathcal{C};L^2(\tilde{\mathcal{A}}_{\mathcal{C}})^0\otimes_{\tilde{\mathcal{A}}_{\mathcal{C}}} \mathcal{A}_{\mathcal{D}}^{\otimes_{\tilde{\mathcal{A}}_{\mathcal{C}}} q}\right)\right]\,.
    \end{align*}
    By the induction hypothesis and \thref{thm:local-scaling-formula}, all terms in the sum on the right-hand side are zero, so the claim follows.
\end{proof}

\section{\texorpdfstring{$L^2$}{L2}-Betti numbers associated with Hecke pairs}
Suppose that $\Gamma$ is a countable discrete group acting on a type $\mathrm{II}_1$ factor $P$ by outer automorphisms, and let $\Lambda<\Gamma$ be an almost-normal subgroup.
Then $P\rtimes\Lambda\subset P\rtimes\Gamma$ is a quasi-regular inclusion, which is unimodular in the sense of \cite[Definition~3.3]{psv-cohom} precisely when $\Lambda<\Gamma$ is.
On the other hand, the Schlichting completion $G$ is also unimodular in this case, so we have two ways to associate $L^2$-Betti numbers to the Hecke pair $\Lambda<\Gamma$.
These turn out to yield the same values, which is the content of the theorem below.
\begin{theorem}
\label{thm:hecke-inclusion-computation}
Let $\Lambda<\Gamma$ be a unimodular Hecke pair and $P$ a type $\mathrm{II}_1$ factor.
Denote the Schlichting completion of $\Lambda<\Gamma$ by $K<G$, and equip $G$ with the Haar measure normalised such that $\mu(K)=1$.
Given an outer action $\alpha$ of $\Gamma$ on $P$ by trace-preserving automorphisms, put $T=P\rtimes\Lambda$, $S=P\rtimes\Gamma$ and $\mathcal{S}=\mathrm{QN}_S(T)$.
Then $\beta^{(2)}_n(T\subset S) = \beta^{(2)}_n(G)$ for all $n\geq 0$.
\end{theorem}

\begin{example}
    \label{exa:l2-betti-contrast}
    Contrasted with \thref{thm:arano-vaes-vanishing}, \thref{thm:hecke-inclusion-computation} shows that the behaviour of $P\rtimes\Lambda\subset P\rtimes\Gamma$ is potentially different from that of $\mathcal{C}_f(K<G)$.
    To confirm that this is actually the case, we have to exhibit a Hecke pair $\Lambda<\Gamma$ such that the Schlichting completion is nondiscrete and has nonvanishing $L^2$-Betti numbers.
    Such examples exist in the literature: one can take $\Lambda=\mathrm{SL}_n(\ZZ)$ and $\mathrm{SL}_n(\ZZ[1/p])$ for $n\in\NN$ and $p$ a prime larger than $n$.
    Then \cite[Theorem~5.30]{hdp-thesis} shows that the $n$th $L^2$-Betti number of $G=\mathrm{PSL}_n(\QQ_p)$ does \textit{not} vanish.

    One can also use \thref{thm:hecke-inclusion-computation} together with the free product formula of \cite[Proposition~9.4]{psv-cohom} to produce such examples.
    Indeed, suppose that an infinite group $\Lambda$ embeds properly inside two discrete groups $\Gamma_1$ and $\Gamma_2$ such that $\Lambda<\Gamma_i$ is a unimodular reduced Hecke pair.
    Fix an outer action of the amalgamated free product $\Gamma=\Gamma_1 *_\Lambda \Gamma_2$ on a $\mathrm{II}_1$ factor $P$. 
    In particular, this induces outer actions of $\Lambda$, $\Gamma_1$ and $\Gamma_2$ on $P$.
    Note that the canonical copy of $\Lambda$ in $\Gamma$ is still a unimodular reduced Hecke pair.
    If we now put $T=P\rtimes\Lambda$, $S=P\rtimes\Gamma$ and $S_i=P\rtimes\Gamma_i$, we can identify
    \[
        (T\subset S) \cong (T\subset S_1 *_T S_2)\,,
    \]
    where the amalgamated free product is taken w.r.t.\@ the trace-preserving conditional expectations;
    see e.g.\@ \cite[Proposition~2.5]{ueda-amalg}.
    Let $K<G$ be the Schlichting completion of $\Lambda<\Gamma$.
    Since $\Lambda$ is infinite and $\Lambda<\Gamma$ is reduced, $G$ is nondiscrete.
    Moreover, by \cite[Proposition~9.4]{psv-cohom} and \thref{thm:hecke-inclusion-computation}, the first $L^2$-Betti number of $G$ is at least $1$.
    In particular, $\beta^{(2)}_1(\mathcal{C}_f(K<G))$ vanishes while $\beta^{(2)}_1(G)$ does not.
\end{example}

Before we can proceed with the proof of \thref{thm:hecke-inclusion-computation}, we need a number of technical auxiliary results, explained in the three subsections below.
In the process, we also obtain a cohomological dictionary between $\Lambda<\Gamma$ and $T\subset\mathcal{S}$, of which we discuss a few other applications.

\subsection{Bimodules and unitary representations}
\label{sec:bimod-unitary-reps}
Consider a discrete group $\Gamma$ acting on a type $\mathrm{II}_1$ factor $P$ by outer automorphisms and put $S=P\rtimes\Gamma$.
Let $\mathcal{H}$ be a Hilbert space equipped with a unitary right action of $\Gamma$.
We turn $L^2(S)\htimes\mathcal{H}$ into a Hilbert $S$-$S$-bimodule via
\begin{equation}
\label{eqn:g-rep-bimodule}
    x(\xi\otimes \eta) p u_g = x\xi pu_g\otimes \eta \cdot g\,,
\end{equation}
where $x\in S$, $p\in P$, $\xi\in L^2(S)$, $\eta\in\mathcal{H}$ and $g\in\Gamma$ (compare \cite[\S~2.1]{pv-superrigidity-amalgamated} and \cite[Proposition~2.7]{arano-vaes}).
The fact that \eqref{eqn:g-rep-bimodule} yields a normal right representation of $S$ on \eqref{eqn:g-rep-bimodule} is a consequence of the Fell absorption principle.
Often, $\mathcal{H}$ will be the permutation representation of $\Gamma$ on $\ell^2(I)$ coming from some right action of $\Gamma$ on a set $I$.
In this case, we abbreviate the $S$-$S$-bimodule $L^2(S)\htimes \ell^2(I)$ as $\mathcal{K}(\Gamma; I)$.
\begin{lemma}
\label{thm:i-bimod-intertw}
Let $\Gamma$ be a discrete group acting on a type $\mathrm{II}_1$ factor $P$ by outer automorphisms.
For any unitary antirepresentation $\Gamma\curvearrowright\mathcal{H}$, view $L^2(S)\htimes\mathcal{H}$ as an $S$-$S$-bimodule via \eqref{eqn:g-rep-bimodule}.
The assignment $F:\mathcal{H}\mapsto L^2(S)\htimes \mathcal{H}$ is a fully faithful functor from the category $\Rep(\Gamma^{\op})$ of unitary antirepresentations of $\Gamma$ to the category $\mathrm{HilbBimod}_{S-S}$ of Hilbert $S$-$S$-bimodules.
Here, the morphisms in $\Rep(\Gamma^{\op})$ and $\mathrm{HilbBimod}_{S-S}$ are bounded maps intertwining the relevant actions, and the action of $F$ on morphisms is defined by $F(V)=\I\otimes V$.

In particular, all Hilbert $S$-$S$-subbimodules of $L^2(S)\htimes\mathcal{H}$ are of the form $L^2(S)\htimes\mathcal{K}$ for some closed $\Gamma$-subspace $\mathcal{K}\subset\mathcal{H}$.
\begin{proof}
    We have to show that the map
    \[
        \Phi:\BHom_{-\Gamma}(\mathcal{H},\mathcal{K})\to \BHom_{S-S}(L^2(S)\htimes\mathcal{H},L^2(S)\htimes\mathcal{K}): V\mapsto \I\otimes V
    \]
    is an isomorphism for all $\mathcal{H},\mathcal{K}\in\Rep(\Gamma^{\op})$
    Clearly, $\Phi$ is well-defined and injective.
    To see that $\Phi$ is also surjective, fix a bounded $S$-$S$-bimodular map $W:L^2(S)\htimes\mathcal{H}\to L^2(S)\htimes\mathcal{K}$.
    Then, for all $\xi\in\mathcal{H}$, $W(\I\otimes\xi)$ is a $P$-central vector in $L^2(S)\htimes\mathcal{K}$.
    Expanding this in terms of an orthonormal basis of $\mathcal{K}$ and using the fact that $P\subset P\rtimes\Gamma$ is irreducible, we find that $W(\I\otimes\xi)=\I\otimes\eta$ for some unique $\eta\in\mathcal{K}$.
    If we now define $V(\xi)=\eta$, easy computation shows that $V:\mathcal{H}\to\mathcal{K}$ is a bounded $\Gamma$-intertwiner, and that $W=\I\otimes V$.
\end{proof}
\end{lemma}

The structure of $\mathcal{K}(\Gamma; I)$ is fairly easy to describe in terms of the action of $\Gamma$ on $I$.
\begin{lemma}
\label{thm:i-bimod-structure}
Let $\Gamma$ be a discrete group acting on a type $\mathrm{II}_1$ factor $P$ by outer automorphisms, and $I$ a set equipped with a right action of $\Gamma$.
Denote the set of $\Gamma$-orbits in $I$ by $\Orb(I)$.
Given $\alpha\in I$, denote the stabiliser of $\alpha$ in $\Gamma$ by $\Gamma_\alpha$.
Then the Hilbert $S$-$S$-bimodule $\mathcal{K}(\Gamma; I)$ decomposes as follows:
\begin{equation}
\label{eqn:i-bimod-decomp}
    \mathcal{K}(\Gamma; I) \cong \bigoplus_{\alpha\cdot\Gamma \in \Orb(I)} \mathcal{K}(\Gamma; \Gamma_\alpha\backslash \Gamma)\,.
\end{equation}

The left and right tracial weights \cite[\S~8.4]{adp} on $\BEnd_{S-S}(\mathcal{K}(\Gamma; I))$ agree.
If $\Gamma$ acts on $I$ with finite orbits, the categorical trace is given by
\begin{equation}
\label{eqn:i-bimod-cat-trace-fin-orb}
    \tau_{S-S}(V)= \sum_{\alpha\cdot\Gamma \in \Orb(I)} [\Gamma:\Gamma_\alpha] \langle V(\I\otimes \delta_{\alpha}), \I\otimes \delta_{\alpha}\rangle = \sum_{\alpha\cdot\Gamma \in \Orb(I)} [\Gamma:\Gamma_\alpha] \langle V_0(\delta_{\alpha}), \delta_{\alpha}\rangle
\end{equation}
for $V=\I\otimes V_0\in \BEnd_{S-S}(\mathcal{K}(\Gamma; I))_+$.
\begin{proof}
    The decomposition \eqref{eqn:i-bimod-decomp} is immediate from the orbit-stabiliser theorem.
    Since the vectors $(\I\otimes\partial_i)_{i\in I}$ are a basis for $\mathcal{K}(\Gamma; I)$ as a right $S$-module and as a left $S$-module, \eqref{eqn:i-bimod-cat-trace-fin-orb} follows.
\end{proof}
\end{lemma}

\subsection{The regular bimodule}
Throughout, we fix a countable discrete group $\Gamma$ with an outer action $\alpha:\Gamma\curvearrowright P$ on a type $\mathrm{II}_1$ factor $P$ once and for all. 
The crossed product $P\rtimes \Gamma$ will be denoted by $S$.
Given an almost-normal subgroup $\Lambda<\Gamma$, the regular $S$-$S$-bimodule \eqref{eqn:reg-bimodule-def} associated with $P\rtimes\Lambda\subset S$ can be written in a way that makes the connection with the Schlichting completion more apparent.
To this end, we apply the construction from the previous subsection.

For now, let $\Lambda<\Gamma$ be any subgroup of $\Gamma$ and put $T=P\rtimes\Lambda$.
Recall that a Hilbert $T$-$T$-bimodule $\mathcal{H}$ gives rise to a Hilbert $S$-$S$-bimodule $L^2(S)\cftimes{T} \mathcal{H} \cftimes{T} L^2(S)$ by induction.
On the bimodules discussed in the previous section, this induction operation has a particularly convenient description.
\begin{proposition}
\label{thm:i-bimod-induction}
Consider a set $I$ with a right action of $\Gamma$, and let $\Lambda_0<\Lambda$ be an arbitrary subgroup.
Then there are canonical unitary isomorphisms of $S$-$S$-bimodules
\begin{align*}
     \mathcal{K}(\Gamma; \Lambda_0\backslash\Gamma)&\cong L^2(S)\cftimes{T} \mathcal{K}(\Lambda; \Lambda_0\backslash\Lambda)\cftimes{T} L^2(S)\,,\tag{I}\\
     \mathcal{K}(\Gamma; \Lambda\backslash\Gamma \times I) &\cong L^2(S)\cftimes{T} \mathcal{K}(\Gamma; I)\,,\tag{II}\\
     \mathcal{K}(\Gamma; I\times\Lambda\backslash\Gamma ) &\cong \mathcal{K}(\Gamma; I) \cftimes{T} L^2(S)\,,\tag{III}\\
     &\cong L^2(S)\cftimes{T} \mathcal{K}(\Lambda; I)\cftimes{T} L^2(S)\,.\tag{IV}
\end{align*}
\begin{proof}
    \newcommand{\WidestEntry}{$xu_{g^{-1}}\otimes_T (\I\otimes \delta_{ig^{-1}})\otimes_T u_g$}%
    \newcommand{\SetToWidest}[1]{\makebox[\widthof{\WidestEntry}]{$#1$}}%
    The isomorphisms and their inverses are, in order,
    \begingroup\allowdisplaybreaks
    \begin{align*}
        \mathrm{(I)} 
        &\left\{\begin{array}{ccc}
            \SetToWidest{x \otimes \delta_{\Lambda_0 g}} &\mapsto& \SetToWidest{xu_{g^{-1}}\otimes_T (\I\otimes \delta_{\Lambda_0})\otimes_T u_g}\\[1ex]
            \SetToWidest{x \otimes_T (t\otimes\delta_{\Lambda_0 h}) \otimes_T p u_g} &\mapsto& \SetToWidest{xtpu_g\otimes \delta_{\Lambda_0 hg}}
        \end{array}\right.\\[2ex]
        \mathrm{(II)}
        &\left\{\begin{array}{ccc}
            \SetToWidest{x \otimes \delta_{\Lambda g,i}} &\mapsto& \SetToWidest{x u_{g^{-1}} \otimes_T (u_g\otimes \delta_i)}\\[1ex]
            \SetToWidest{x\otimes_T (pu_g\otimes \delta_i)} &\mapsto& \SetToWidest{xpu_g\otimes \delta_{\Lambda g, i}}
        \end{array}\right.\\[2ex]
        \mathrm{(III)}
        &\left\{\begin{array}{ccc}
            \SetToWidest{x \otimes \delta_{i,\Lambda g}} &\mapsto& \SetToWidest{(x u_{g^{-1}}\otimes \delta_{ig^{-1}})\otimes_T u_g}\\[1ex]
            \SetToWidest{(x\otimes \delta_i)\otimes_T pu_g} &\mapsto& \SetToWidest{xpu_g \otimes\delta_{ig, \Lambda g}}
        \end{array}\right.\\[2ex]
        \mathrm{(IV)}
        &\left\{\begin{array}{ccc}
            \SetToWidest{x \otimes \delta_{i,\Lambda g}} &\mapsto& \SetToWidest{xu_{g^{-1}}\otimes_T (\I\otimes \delta_{ig^{-1}})\otimes_T u_g}\,\\[1ex]
            \SetToWidest{x \otimes_T (t\otimes\delta_i) \otimes_T p u_g} &\mapsto& \SetToWidest{xtpu_g\otimes \delta_{ig,\Lambda g}}
        \end{array}\right.
    \end{align*}
    \endgroup
    where $x\in S$, $t\in T$, $p\in P$, $i\in I$ and $g,h\in\Gamma$.
\end{proof}
\end{proposition}

Inductive application of this result immediately yields the following.
\begin{corollary}
\label{thm:k-lambda-rel-tp}
    There are canonical isometric $S$-$S$-bimodule isomorphisms
    \[
        L^2(S)^{\cfpow{T}{n+1}} \cong \mathcal{K}(\Gamma; (\Lambda\backslash\Gamma)^n)
    \]
    for all $n\geq 0$.
    Here, we take $(\Lambda\backslash\Gamma)^0$ to be a one-point set equipped with the trivial action of $\Gamma$.
\end{corollary}
Interpreted in this way, the second claim of \thref{thm:i-bimod-induction} is simply a reformulation of the statement 
\[
    L^2(S)\cftimes{P\rtimes \Lambda_0} L^2(S) \cong L^2(S) \cftimes{T} \left(L^2(T)\cftimes{P\rtimes\Lambda_0} L^2(T)\right) \cftimes{T} L^2(S)\,.
\]

From now on, we assume $\Lambda$ to be an almost-normal subgroup of $\Gamma$ with $\Lambda<\Gamma$ unimodular.
Recall from \autoref{sec:quasireg-prelim} that the \textit{regular} Hilbert $S$-$S$-bimodule $\mathcal{H}_{\mathrm{reg}}$ associated with the quasi-regular inclusion $T\subset S$ is defined by
\[
    \mathcal{H}_{\mathrm{reg}}^0 = \bigoplus_{n\in\NN} L^2(S)^{\cfpow{T}{n}}
    \qquad\textrm{and}\qquad 
    \mathcal{H}_{\mathrm{reg}} = L^2(S)\cftimes{T}\mathcal{H}_{\mathrm{reg}}^0\cftimes{T} L^2(S)\,.
\]

As in \cite{psv-cohom}, $\mathcal{M}(T\subset S)$ refers to the von Neumann algebra $\BEnd_{S-S}(\mathcal{H}_{\mathrm{reg}})$. 
This algebra comes equipped with a canonical tracial weight $\Tr$ characterised by the following condition: if $\mathcal{H}\subset\mathcal{H}^0_{\mathrm{reg}}$ is a bifinite $T$-$T$-subbimodule and $p:\mathcal{H}^0_{\mathrm{reg}}\to \mathcal{H}$ is the associated orthogonal projection, we have that
\begin{equation}
\label{eqn:m-ts-tracial-weight}
    \Tr((\I\otimes p\otimes\I) x(\I\otimes p\otimes\I))=\Tr_{\mathcal{H}}(p\iota^* x\iota p)
\end{equation}
for all $x\in\mathcal{M}(T\subset S)_+$, where $\Tr_{\mathcal{H}}(\cdot)$ denotes the categorical trace\footnote{Since we assume $T\subset S$ to be unimodular, the categorical trace is equal to both the right and the left trace.} on $\BEnd_{T-T}(\mathcal{H})$ and $\iota:\mathcal{H}^0_{\mathrm{reg}}\to \mathcal{H}_{\mathrm{reg}}$ is the canonical induction map.

In light of \thref{thm:i-bimod-induction} and \thref{thm:k-lambda-rel-tp}, we will identify $\mathcal{H}_{\mathrm{reg}}^0$ and $\mathcal{H}_{\mathrm{reg}}$ with
\begin{equation}
\label{eqn:regular-bimodule-model}
    \mathcal{H}_{\mathrm{reg}}^0 = L^2(T)\oplus \bigoplus_{n\geq 0} \mathcal{K}(\Gamma; (\Lambda\backslash\Gamma)^n)
    \qquad \textrm{and} \qquad
    \mathcal{H}_{\mathrm{reg}} = \bigoplus_{n\geq 1} \mathcal{K}(\Gamma; (\Lambda\backslash\Gamma)^n)\,.
\end{equation}
A key technical tool in the $L^2$-Betti number computation is the observation that the regular bimodule $\mathcal{H}_{\mathrm{reg}}$ contains copies of $\mathcal{K}(\Gamma;\Lambda_0\backslash\Gamma)$ when $\Lambda_0$ ranges over a family of ``nice'' subgroups of $\Lambda$.
More concretely, let $\mathcal{F}_{\Lambda<\Gamma}$ be the set of subgroups of $\Lambda$ of the form $\Lambda_0=\bigcap_{i=0}^n x_i\Lambda x_i^{-1}$ for some $x_0=e,x_1\ldots, x_n\in\Gamma$.
Given any $\Lambda_0\in\mathcal{F}_{\Lambda<\Gamma}$, $\Lambda_0<\Gamma$ is of course still a Hecke pair, giving rise to the same Schlichting completion $G$ as $\Lambda<\Gamma$.
The following result then describes how to embed $\mathcal{K}(\Gamma;\Lambda_0\backslash\Gamma)$ inside $\mathcal{H}_{\mathrm{reg}}$ in a convenient position.
\begin{proposition}
\label{thm:lambda0-central-gen-vectors}
    For all $\Lambda_0\in\mathcal{F}_{\Lambda<\Gamma}$, there exists an $S$-$S$-bimodular isometry $W:\mathcal{K}(\Gamma; \Lambda_0\backslash\Gamma)\to \mathcal{H}_{\mathrm{reg}}$ that implements a $*$-isomorphism 
    \begin{equation}
    \label{eqn:lambda0-corner-endo-isos}
        \BEnd_{S-S}(\mathcal{K}(\Gamma; \Lambda_0\backslash\Gamma))
        \cong (WW^*)\mathcal{M}(T\subset S)(WW^*)
    \end{equation}
    satisfying
    \begin{equation}
    \label{eqn:lambda0-corner-traces}
        \Tr\left((WW^*)y(WW^*)\right)=[\Lambda:\Lambda_0]\langle (W^*yW) (\I\otimes \delta_{\Lambda_0}), (\I\otimes \delta_{\Lambda_0})\rangle
    \end{equation}
    for all $y\in\mathcal{M}(T\subset S)_+$.
    In addition, these isometries have the following properties:
    \begin{enumerate}[(i)]
        \item denoting the Schlichting completion of $\Lambda<\Gamma$ by $G$, and the respective closures of $\Lambda$ and $\Lambda_0$ by $K$ and $K_0$, there is a natural trace-preserving $*$-isomorphism
        \begin{equation} 
            (WW^*)\mathcal{M}(T\subset S)(WW^*)\cong p_{K_0} L(G)p_{K_0} \,,
        \end{equation}
        where $p_{K_0}$ is the projection in $L(G)$ defined by \eqref{eqn:subgroup-invar-projection};
        \item if we choose such a $W_{\Lambda_0}$ for all $\Lambda_0\in\mathcal{F}_{\Lambda<\Gamma}$ and denote the central support of $W_{\Lambda_0}W_{\Lambda_0}^*$ in $\mathcal{M}(T\subset S)$ by $z_{\Lambda_0}$, then
        \begin{equation}
        \label{eqn:lambda0-central-gen-dense}
            \mathcal{H}_{\mathrm{reg}} = \overline{\bigcup_{\substack{\Lambda_0\in\mathcal{F}_{\Lambda<\Gamma}}} z_{\Lambda_0}\mathcal{H}_{\mathrm{reg}}}\,.
        \end{equation}
    \end{enumerate}
    \begin{proof}
        Fix $\{x_1,\ldots,x_n\}\subset \Gamma$ and put $\Lambda_0=\Lambda\cap \bigcap_{i=1}^n x_i^{-1}\Lambda x_i$.
        Recall from \thref{thm:i-bimod-structure} that there is a $S$-$S$-bimodular embedding
        \[
            \mathcal{K}(\Gamma;\Lambda_0\backslash\Gamma)\to \mathcal{K}\left(\Gamma,(\Lambda\backslash\Gamma)^n\right)\subset\mathcal{H}_{\mathrm{reg}}^0:x\otimes\delta_{\Lambda_0 g} \mapsto x\otimes\delta_{\Lambda x_1g,\ldots,\Lambda x_ng}\,.
        \]
        By restriction, we get a $T$-$T$-bimodular embedding
        \[
            W_0: \mathcal{K}(\Lambda;\Lambda_0\backslash\Lambda)\to \mathcal{K}\left(\Gamma,(\Lambda\backslash\Gamma)^n\right)\subset\mathcal{H}^0_{\mathrm{reg}}\,.
        \]
        Through the isomorphisms in \thref{thm:i-bimod-induction}, this map induces an $S$-$S$-bimodular embedding
        \[
            W: \mathcal{K}(\Gamma;\Lambda_0\backslash\Gamma)\to \mathcal{K}\left(\Gamma,(\Lambda\backslash\Gamma)^{n+2}\right)\subset\mathcal{H}_{\mathrm{reg}}:
            x\otimes\delta_{\Lambda_0 g}\mapsto x\otimes\delta_{\Lambda g,\Lambda x_1g,\ldots,\Lambda x_ng,\Lambda g}\,.
        \]
        Again under the identifications of \thref{thm:i-bimod-induction}, the projections $\I\otimes_T W_0W_0^*\otimes_T\I$ and $WW^*$ represent the same endomorphism of $\mathcal{H}_{\mathrm{reg}}$.
        Using \thref{thm:i-bimod-structure} to compute the categorical trace on $\BEnd_{T-T}(\mathcal{K}(\Lambda;\Lambda_0\backslash\Lambda))$, it follows that $W$ satisfies \eqref{eqn:lambda0-corner-traces}.

        To prove (i), recall from \thref{thm:i-bimod-intertw} that $\BEnd_{S-S}(\mathcal{K}(\Gamma; \Lambda_0\backslash\Gamma))$ is given by the commutant of the right action of $\Gamma$ on $\ell^2(\Lambda_0\backslash\Gamma)$.
        In this case $\ell^2(\Lambda_0\backslash\Gamma)\cong \ell^2(K_0\backslash G)$, and said commutant is precisely $p_{K_0} L(G)p_{K_0}$.
        It follows that \eqref{eqn:lambda0-corner-endo-isos} actually gives rise to an isomorphism
        \[ 
            (WW^*)\mathcal{M}(T\subset S)(WW^*)\cong p_{K_0} L(G)p_{K_0} \,.
        \]
        Moreover, \eqref{eqn:lambda0-corner-traces} shows that the trace $\Tr$ on $\mathcal{M}(T\subset S)$ corresponds exactly to the Plancherel trace on $L(G)$ when both are restricted to these corners.

        To see why \eqref{eqn:lambda0-central-gen-dense} is true, first observe that vectors of the form $\I\otimes\delta_{\Lambda y_1,\ldots,\Lambda y_m,\Lambda}$ generate $\mathcal{H}_{\mathrm{reg}}$ as an $S$-$S$-bimodule.
        We will show that all these generators are contained in $z_{\Lambda_0}\mathcal{H}_{\mathrm{reg}}$ for some $\Lambda_0\in\mathcal{F}_{\Lambda<\Gamma}$.
        Given any tuple $\Lambda y_1,\ldots,\Lambda y_m\in (\Lambda\backslash\Gamma)^m$, put $\Lambda_0=\bigcap_{i=1}^m y_i^{-1} \Lambda y_i$ and $\xi = \I\otimes\delta_{\Lambda y_1,\ldots,\Lambda y_m,\Lambda}$.
        As in the first part of the proof, \thref{thm:i-bimod-structure} provides an $S$-$S$-bimodular isometry $W:\mathcal{K}(\Gamma; \Lambda_0\backslash\Gamma)\to \mathcal{H}_{\mathrm{reg}}$ that maps $\I\otimes\delta_{\Lambda_0}$ to $\xi$.
        Composing this with the adjoint of the isometry $W_{\Lambda_0}:\mathcal{K}(\Gamma; \Lambda_0\backslash\Gamma)\to \mathcal{H}_{\mathrm{reg}}$ that we fixed a priori, we get that
        \[
            \xi = \left(WW_{\Lambda_0}^*\right) \left(W_{\Lambda_0}(\I\otimes\delta_{\Lambda_0})\right)\qquad \in\mathcal{M}(T\subset S) (W_{\Lambda_0}W_{\Lambda_0}^*)\mathcal{H}_{\mathrm{reg}} 
            \qquad \subset z_{\Lambda_0}\mathcal{H}_{\mathrm{reg}}\,,
        \]
        as required.
    \end{proof}
\end{proposition} 

\subsection{A cohomological correspondence}
Let $\Lambda$ be an almost-normal subgroup of $\Gamma$.
Put $T=P\rtimes\Lambda$ and $\mathcal{S}=\mathrm{QN}_S(T)$.
For now, we do not assume $\Lambda<\Gamma$ to be unimodular.
In this section, we study the cohomology of the Hilbert $S$-$S$-bimodule $L^2(S)\htimes\mathcal{H}$ defined in \eqref{eqn:g-rep-bimodule} for any Hilbert space $\mathcal{H}$ equipped with a continuous unitary right action of $G$.
Combined with the results of the previous two subsections, we can then prove \thref{thm:hecke-inclusion-computation}.

First, note that
\begin{equation}
\label{eqn:qn-easy-description}
    \mathcal{S} = \mathrm{span}\{t u_g \mid t\in T, g\in\Gamma\}\,.
\end{equation}
Indeed, \cite[Lemma~2.5]{psv-cohom} tells us that $\mathcal{S}$ precisely consists of those elements of $S$ that are $T$-bounded in some irreducible bifinite $T$-$T$-subbimodule of $L^2(S)$.
However, these subbimodules are all of the form $\mathcal{K}=\bigoplus_{k\in \Lambda\cap g^{-1}\Lambda g\backslash\Lambda} L^2(T)u_{gk}\subset L^2(S)$ for $g\in\Gamma$.
The $T$-bounded vectors in $\mathcal{K}$ are spanned by $\{t u_{gk}\mid t\in T, k\in\Lambda\}$, which readily leads to \eqref{eqn:qn-easy-description}.
More generally, we have that the left $T$-module $\mathcal{S}^{\otimes_T^n}$ is free with basis
\[
    \{u_{g_1}\otimes_T \cdots \otimes_T u_{g_n}\mid g_i\in\Lambda\backslash\Gamma, i=1,\ldots,n\}\,.
\]
Indeed, the elements of this set are free in the left Hilbert $T$-module $L^2(S)^{\cfpow{T}{n}}$ by an orthogonality argument, so they are certainly free in $\mathcal{S}^{\otimes_T^n}$ \cite[see also][Lemma~6.3]{psv-cohom}.
The fact that they generate $\mathcal{S}^{\otimes_T^n}$ is clear from \eqref{eqn:qn-easy-description}.

Let $\mathcal{H}$ be a Hilbert space on which $\Gamma$ acts unitarily from the right, and consider $L^2(S)\htimes\mathcal{H}$ with the Hilbert $S$-$S$-bimodule structure defined in \eqref{eqn:g-rep-bimodule}.
There are natural isomorphisms
\begin{equation}
\label{eqn:bar-resolution-identification}
    \Hom_{T-T}(\mathcal{S}^{\otimes_T^n}, L^2(S)\htimes\mathcal{H})
    \cong \Fun((\Lambda\backslash\Gamma)^{n+1}, \mathcal{H})^\Gamma
\end{equation}
for all $n\geq 0$, defined by
\begin{align*}
    &\Phi :
    \Hom_{T-T}(\mathcal{S}^{\otimes_T^n}, L^2(S)\htimes\mathcal{H})\to \Fun((\Lambda\backslash\Gamma)^{n+1}, \mathcal{H})^\Gamma:\\
    &\qquad \langle\Phi(\phi)(g_0,\ldots,g_n),\eta\rangle=\langle u_{g_0^{-1}}\phi(u_{g_0g_1^{-1}}\otimes_T u_{g_1g_2^{-1}}\otimes_T\cdots \otimes_T u_{g_{n-1}g_n^{-1}})u_{g_n}, \I\otimes \eta\rangle\, ,\\[1ex]
    &\Psi :
    \Fun((\Lambda\backslash\Gamma)^{n+1}, \mathcal{H})^\Gamma\to \Hom_{T-T}(\mathcal{S}^{\otimes_T^n}, L^2(S)\htimes\mathcal{H}):\\
    &\qquad\Psi(f)(tu_{s_1}\otimes_T\cdots \otimes_T u_{s_n})=tu_{s_1\cdots s_n}\otimes f(s_1\cdots s_n,s_2\cdots s_n,\ldots,s_n,e)\, ,
\end{align*}
where $g_0,\ldots,g_n,s_1,\ldots,s_n\in\Gamma$, $t\in T$ and $\eta\in\mathcal{H}$.
One checks that these are well-defined, and that $\Phi\circ\Psi=\id$.
To prove that $\Psi\circ\Phi=\id$, one uses the irreducibility of $P\subset S$ to get that $\phi(s_1\otimes_T \cdots \otimes_T s_n)$ lies in $u_{s_1\cdots s_n}\otimes\mathcal{H}\subset L^2(S)\htimes\mathcal{H}$.

Following \cite[Definition~4.1]{psv-cohom}, the cohomology spaces $H^n(T\subset\mathcal{S},L^2(S)\htimes\mathcal{H})$ are defined as the $n$-th cohomology of the complex with terms
\begin{equation}
\label{eqn:ts-hecke-cohom}
    C^n = \Hom_{T-T}(\mathcal{S}^{\otimes_T^n}, L^2(S)\htimes\mathcal{H})
\end{equation}
and the usual differentials for inhomogeneous cochains; see \autoref{sec:quasireg-prelim}.

Let $G$ be the Schlichting completion of $\Gamma$, and $K$ the image of $\Lambda$ in $G$.
Let $\mathcal{H}$ be a Hilbert space equipped with a continuous unitary right representation of $G$.
The natural map $\Gamma\to G$ induces a unitary right representation of $\Gamma$.
Under the identification \eqref{eqn:bar-resolution-identification}, the cochain complex $C^n$ in \eqref{eqn:ts-hecke-cohom} is identified with
\[
    C^n\cong\Fun((\Lambda\backslash\Gamma)^{n+1}, \mathcal{H})^\Gamma=\Fun((K\backslash G)^{n+1}, \mathcal{H})^G =: \tilde{C}^n\,,
\]
where the differentials on $\tilde{C}^n$ are given by dropping coordinates, as explained in \autoref{sec:locally-profinite-cohom}.
Since $K$ is compact and open in $G$, the discussion from \cite{guichardet} reviewed in \autoref{sec:locally-profinite-cohom} shows that the cohomology of $\tilde{C}^\bullet$ is precisely $H^\bullet_c(G;\mathcal{H})$.
This leads to the following result.
\begin{proposition}
\label{thm:schlichting-completion-cohomology}
Let $\Lambda<\Gamma$ be a Hecke pair with Schlichting completion $K<G$.
Put $S=P\rtimes\Gamma$ and $\mathcal{S}=\mathrm{QN}_S(T)$.
For any continuous unitary right representation of $G\curvearrowright\mathcal{H}$, there are canonical isomorphisms in cohomology
\[
    H^n(T\subset\mathcal{S};L^2(S)\htimes\mathcal{H}) \cong H^n_c(G;\mathcal{H})\,.
\]
\end{proposition}

\begin{remark}
\label{rem:schlichting-completion-cohomology-red}
In some sense, \thref{thm:schlichting-completion-cohomology} reduces the cohomology theory of $T\subset S$ to that of the Schlichting completion $G$.
To justify this statement, fix a Hilbert $S$-$S$-bimodule $\mathcal{K}$.
Define $\mathcal{K}'$ as the largest $T$-$T$-subbimodule of $\mathcal{K}$ that can be written as a direct sum of $T$-$T$-bimodules embedding inside $L^2(S)^{\cfpow{T}{n}}$ for some $n\geq 1$.
By maximality, $\mathcal{K}'$ is in fact an $S$-$S$-subbimodule of $\mathcal{K}$.
Now, since $\mathcal{K}'$ embeds into a direct sum of copies of $\mathcal{H}_{\mathrm{reg}}$, $\mathcal{K}'$ is isomorphic to an $S$-$S$-subbimodule of $\mathcal{K}(\Gamma; I)$, where $I$ is a disjoint union of sufficiently many copies of $\bigsqcup_{n\geq 1} (\Lambda\backslash\Gamma)^n$.
By \thref{thm:i-bimod-intertw}, this means that there exists a closed $\Gamma$-subspace $\mathcal{H}$ of $\ell^2(I)$ such that $\mathcal{K}'\cong L^2(S)\htimes\mathcal{H}$ as $S$-$S$-bimodules.
Moreover, since $\mathcal{H}$ is closed, we can upgrade the action of $\Gamma$ on $\mathcal{H}$ to a continuous unitary right representation of the Schlichting completion $G$ by restricting the canonical permutation representation of $G$ on $\ell^2(I)$.

As explained in \cite[Remark~4.2(4)]{psv-cohom}, the cohomology $H^n(T\subset\mathcal{S},\mathcal{K})$ does not change if we replace $\mathcal{K}$ with the Hilbert $S$-$S$-subbimodule $\mathcal{K}'$.
This leads to canonical identifications in cohomology
\[
    H^n(T\subset\mathcal{S},\mathcal{K})
    \cong H^n(T\subset\mathcal{S},L^2(S)\htimes\mathcal{H})
    \cong H^n_c(G;\mathcal{H})\,.
\]
In \autoref{sec:cohom-prop-equiv}, we apply this principle to give cohomological proofs of the equivalence of a number of representation-theoretic properties for $G$ and $T\subset S$.
\end{remark}

\subsubsection{Application to \texorpdfstring{$L^2$}{L2}-Betti numbers}
\label{sec:hecke-l2betti-comp}
In this section, we assume $\Lambda<\Gamma$ to be unimodular.
We continue denoting the Schlichting completion by $K<G$.
Consider a subgroup $\Lambda_0\in\mathcal{F}_{\Lambda<\Gamma}$, and denote the closure of the image of $\Lambda_0$ in $G$ by $K_0$.
Take $\mathcal{H}=\ell^2(\Lambda_0\backslash\Gamma)=\ell^2(K_0\backslash G)$.
Observe that $C^n$ and $\tilde{C}^n$ both come with canonical actions of $p_{K_0}L(G)p_{K_0}$, and that the isomorphism \eqref{eqn:bar-resolution-identification} intertwines these actions.
This allows us to make a slight improvement to \thref{thm:schlichting-completion-cohomology}.
\begin{corollary}
\label{thm:cohomology-identification}
Let $\Lambda<\Gamma$ be a unimodular Hecke pair.
With the same notation as \thref{thm:schlichting-completion-cohomology}, there are isomorphisms of left $p_{K_0}L(G)p_{K_0}$-modules
\[
    H^n(T\subset\mathcal{S}; \mathcal{K}(\Gamma; \Lambda_0\backslash\Gamma)) \cong H^n_c(G; \ell^2(K_0\backslash G))
\]
for all $\Lambda_0\in\mathcal{F}_{\Lambda<\Gamma}$ and $n\geq 0$.
\end{corollary}

We now proceed to prove \thref{thm:hecke-inclusion-computation}. Throughout, we preserve the notational conventions established in this section.
\begin{proof}[Proof of \thref{thm:hecke-inclusion-computation}]  
    For all $\Lambda_0\in\mathcal{F}_{\Lambda<\Gamma}$, choose an isometry $W_{\Lambda_0}:\mathcal{K}(\Gamma; \Lambda_0\backslash\Gamma)\to \mathcal{H}_{\mathrm{reg}}$ with the properties described in \thref{thm:lambda0-central-gen-vectors}.
    Denote the central support of $W_{\Lambda_0}W_{\Lambda_0}^*$ by $z_{\Lambda_0}$.

    With some slight abuse of notation, we denote the completion of any $\Lambda_0\in\mathcal{F}_{\Lambda<\Gamma}$ inside $G$ by $K_0$, and write $p_{K_0}$ for the projection in $L(G)$ defined by \eqref{eqn:subgroup-invar-projection}.
    \thref{thm:lambda0-central-gen-vectors}(i) tells us that
    \begin{equation}
    \label{eqn:lambda0-central-gen-range-proj-trace}
        \Tr(W_{\Lambda_0}W_{\Lambda_0}^*) = [\Lambda:\Lambda_0] = \tau(p_{K_0})\,.
    \end{equation}
    Moreover, combining \thref{thm:lambda0-central-gen-vectors}(i) with \thref{thm:cohomology-identification} we get the identity
    \begin{align*}
        \dim_{(W_{\Lambda_0}W_{\Lambda_0}^*)\mathcal{M}(T\subset S)(W_{\Lambda_0}W_{\Lambda_0}^*)} &H^n(T\subset\mathcal{S}; (W_{\Lambda_0}W_{\Lambda_0}^*)\mathcal{H}_{\mathrm{reg}})\\
        &=\dim_{p_{K_0}L(G)p_{K_0}} H^n_c(G;\ell^2(K_0\backslash G))\,.
    \end{align*}
    Combining \eqref{eqn:lambda0-central-gen-range-proj-trace} with an appeal to \cite[Lemma~A.16]{kpv13} we then get that
    \begin{align*}
        \dim_{\mathcal{M}(T\subset S)} &H^n(T\subset\mathcal{S}; z_{\Lambda_0}\mathcal{H}_{\mathrm{reg}}) \\
        &= \Tr(W_{\Lambda_0}W_{\Lambda_0}^*)\dim_{(W_{\Lambda_0}W_{\Lambda_0}^*)\mathcal{M}(T\subset S)(W_{\Lambda_0}W_{\Lambda_0}^*)} H^n(T\subset\mathcal{S}; (W_{\Lambda_0}W_{\Lambda_0}^*)\mathcal{H}_{\mathrm{reg}})
\\
        &=[\Lambda:\Lambda_0]\dim_{p_{K_0}L(G)p_{K_0}} H^n_c(G;\ell^2(K_0\backslash G))\,.
    \end{align*}
    Finally, we apply \eqref{eqn:lambda0-central-gen-dense} from \thref{thm:lambda0-central-gen-vectors}(ii), take the supremum over all $\Lambda_0\in\mathcal{F}_{\Lambda<\Gamma}$ and then appeal to \eqref{eqn:subgroup-shrinking-formula} and \thref{rem:hecke-neighbourhood-basis} to find that
    \begin{align*}
        \beta^{(2)}_n(T\subset\mathcal{S}) &= \dim_{\mathcal{M}(T\subset S)} H^n(T\subset\mathcal{S}; \mathcal{H}_{\mathrm{reg}}) 
        =\sup_{\Lambda_0\in\mathcal{F}_{\Lambda<\Gamma}}\dim_{\mathcal{M}(T\subset S)} H^n(T\subset\mathcal{S}; z_{\Lambda_0}\mathcal{H}_{\mathrm{reg}}) \\
        &=\sup_{\Lambda_0\in\mathcal{F}_{\Lambda<\Gamma}}\,[\Lambda:\Lambda_0]\dim_{p_{K_0}L(G)p_{K_0}} H^n_c(G;\ell^2(K_0\backslash G))
        =\beta^{(2)}_n(G)\,,
    \end{align*}
    which concludes the proof.
\end{proof}

\subsubsection{Correspondences in one-cohomology}
\label{sec:cohom-prop-equiv}

Besides its usefulness in studying $L^2$-Betti numbers, the cohomological identification of \thref{thm:schlichting-completion-cohomology} also allows us to establish a direct link between various properties of $G$ and their relative counterparts for $T\subset S$.
As stated in \thref{rem:schlichting-completion-cohomology-red}, we now proceed to explain this in further depth.
We preserve the notation used in the previous subsection, but we no longer assume $G$ to be unimodular.

Suppose that we are given a continuous unitary right action of $G$ on a separable Hilbert space $\mathcal{H}$.
Recall that a 1-cocycle for $G$ with coefficients in $\mathcal{H}$ is a continuous function $c:G\to\mathcal{H}$ satisfying the relation $c(xy)=c(x)y+c(y)$.
Similarly, \cite{psv-cohom} defines a 1-cocycle for $T\subset S$ with coefficients in a Hilbert $S$-$S$-bimodule $\mathcal{K}$ as a $T$-bimodular derivation $c:\mathcal{S}\to\mathcal{K}$, where $\mathcal{S}=\mathrm{QN}_S(T)$.
It is immediate from the definition that the following properties are equivalent for any given 1-cocycle $c:G\to\mathcal{H}$:
\begin{itemize}
    \item $c$ is left $K$-invariant, i.e.\@ $c(kx)=c(x)$ for all $k\in K$ and $x\in G$;
    \item $c$ is right $K$-equivariant, i.e.\@ $c(xk)=c(x)k$ for all $k\in K$ and $x\in G$;
    \item $c$ vanishes on $K$.
\end{itemize}
We will call such 1-cocycles $K$\textit{-trivial}.
Observe that any 1-cocycle $c:G\to\mathcal{H}$ is cohomologous to a $K$-trivial one.
Indeed, if we put $\xi=\mu(K)^{-1}\int_K c(k)\diff \mu(k)$ w.r.t.\@ a Haar measure $\mu$, then the coboundary $d(x)=\xi-\xi x$ agrees with $c$ on $K$, so $c-d$ is $K$-trivial.
In particular, a coboundary is $K$-trivial if and only if it is given by a $K$-invariant vector.
Under the isomorphism \eqref{eqn:bar-resolution-identification} and after further identifying $\Fun(\Lambda\backslash\Gamma\times\Lambda\backslash\Gamma,\mathcal{H})^\Gamma\cong \Fun(\Lambda\backslash\Gamma, \mathcal{H})^\Lambda\cong\Fun(K\backslash G, \mathcal{H})^K$, the 1-cocycles of $T\subset S$ with coefficients in $L^2(S)\htimes\mathcal{H}$ correspond exactly to $K$-trivial 1-cocycles of $G$ with coefficients in $\mathcal{H}$.
Given a $K$-trivial 1-cocycle $c:G\to\mathcal{H}$, the corresponding $T$-$T$-bimodular derivation $\tilde{c}:\mathcal{S}\to L^2(S)\htimes\mathcal{H}$ is given by
\begin{equation}
\label{eqn:cocycle-correspondence}
    \tilde{c}(tu_g)=tu_g\otimes c(u_g)\,,
\end{equation}
where $t\in T$ and $g\in \Gamma$.
Conversely, if we start with a $T$-$T$-bimodular derivation $\tilde{c}:\mathcal{S}\to L^2(S)\htimes\mathcal{H}$, the irreducibility of $P\subset S$ ensures that substituting $t=\I$ in \eqref{eqn:cocycle-correspondence} defines a left $\Lambda$-invariant 1-cocycle $c:\Gamma\to \mathcal{H}$.
The left $\Lambda$-invariance then allows us to extend $c$ uniquely to a $K$-trivial $G$-cocycle.

This correspondence preserves all properties of 1-cocycles defined in \cite[Definition~9.15]{psv-cohom}: 
\begin{definition}[\cite{psv-cohom}]
\label{thm:one-cocycle-defs-qr}
    Consider an $S$-$S$-bimodule $\mathcal{K}$ and a 1-cocycle $c:\mathcal{S}\to \mathcal{K}$.
    Then
    \begin{itemize}
        \item $c$ is \textit{inner} if there exists a $T$-central vector $\xi$ in $\mathcal{K}$ such that $c(x)=x\xi-\xi x$ for all $x\in \mathcal{S}$\,;
        \item $c$ is \textit{approximately inner} if there exists a sequence $(\xi_n)_{n\in\NN}$ of $T$-central vectors in $\mathcal{K}$ such that $c(x)-(x\xi_n-\xi_n x)\to 0$ for all $x\in\mathcal{S}$ \,;
        \item $c$ is \textit{proper} if for any $\kappa>0$ the set of $x\in\mathcal{S}$ such that $\|c(x)\|\leq\kappa\|x\|_2$ is contained in a bifinite $T$-$T$-subbimodule of $L^2(S)$.
    \end{itemize}
\end{definition}
For $G$-cocycles, the analogous definitions are as follows.
\begin{definition}
\label{thm:one-cocycle-defs-groups}
    Consider a Hilbert space $\mathcal{H}$ equipped with a continuous right action of $G$ by unitaries, and a 1-cocycle $c:G\to\mathcal{H}$.
    Then
    \begin{itemize}
        \item $c$ is \textit{inner} if there exists $\xi$ in $\mathcal{H}$ such that $c(g)=\xi-\xi g$ for all $g\in G$\,;
        \item $c$ is \textit{approximately inner} if there exists a sequence $(\xi_n)_{n\in\NN}$ in $\mathcal{H}$ such that $c(g)-(\xi_n-\xi_n g)\to 0$ for $g\in G$ uniformly on compact subsets of $G$ \cite[\S~III.2.3]{guichardet}\,;
        \item $c$ is \textit{proper} if $g\mapsto \|c(g)\|$ is a proper map, i.e. if for any $\kappa>0$ the set of $g\in G$ such that $\|c(g)\|\leq\kappa$ is compact.
    \end{itemize}
\end{definition}
\begin{proposition}
\label{thm:cohom-correspondence}
    Let $\mathcal{H}$ be a Hilbert space equipped with a continuous right action of $G$ by unitaries.
    A $K$-trivial 1-cocycle $c:G\to\mathcal{H}$ is inner (resp.\@ approximately inner, proper) if and only if the corresponding $T$-$T$-bimodular derivation $\tilde{c}:\mathcal{S}\to L^2(S)\htimes\mathcal{H}$ specified in \eqref{eqn:cocycle-correspondence} is.
\begin{proof}
    First, observe that the $T$-central vectors in $L^2(S)\htimes\mathcal{H}$ are precisely given by $\I\otimes\xi$ where $\xi\in\mathcal{H}$ is $K$-invariant.
    For such vectors, we have that 
    \begin{equation}
    \label{eqn:inner-diff-equation}
        \tilde{c}(u_g)-(u_g(\I\otimes\xi)-(\I\otimes\xi) u_g) = u_g\otimes (c(g) - (\xi-\xi g))\,.
    \end{equation}
    If $\tilde{c}$ is approximately inner and $(\xi_n)_n$ a sequence of $K$-invariant vectors in $\mathcal{H}$ such that $\tilde{c}(x)-(x(\I\otimes\xi_n)-(\I\otimes\xi_n) x)\to 0$, the above identity tells us that $c(g)-(\xi_n-\xi_n g)\to 0$ uniformly on compacts.
    Indeed, for any fixed compact set $Q\subset G$, we can find a finite set $\{g_1,\ldots,g_m\}\subset \Gamma$ such that $Q\subset \bigcup_i Kg_i$.
    For any $\epsilon>0$, fix $N$ such that $\|\tilde{c}(u_{g_i})-(u_{g_i}(\I\otimes\xi_n)-(\I\otimes\xi_n)u_{g_i})\|<\epsilon$ for all $n\geq N$.
    Then \eqref{eqn:inner-diff-equation} yields that $\|c(g)-(\xi_n-\xi_ng)\|<\epsilon$ for all $g\in Q$, proving the claim.
    To prove that $\tilde{c}$ is approximately inner exactly when $c$ is, it is now sufficient to justify that the approximate innerness of $c$ can be witnessed by a sequence of $K$-invariant vectors.
    This is indeed the case: let $p_K:\mathcal{H}\to\mathcal{H}$ be the projection onto the space of $K$-invariant vectors, and fix a sequence $(\xi_n)_n$ in $\mathcal{H}$ such that $c(g)-(\xi_n - \xi_n g)\to 0$ uniformly on compact subsets of $G$.
    Then
    \begin{align*}
        c(g)-(p_K(\xi_n) - p_K(\xi_n)g)
        &= \frac{1}{\mu(K)} \int_K c(g) - (\xi_n k - \xi_n k g) \diff \mu(k)\\
        &= \frac{1}{\mu(K)} \int_K \left(c(kgk^{-1})-(\xi_n - \xi_n kgk^{-1})\right)k\diff\mu(k)\,,
    \end{align*}
    which converges to zero because $c(h)-(\xi_n-\xi_n h)\to 0$ uniformly on $KgK$.
    This leaves the equivalence of properness to be explained.
    Since $c$ is $K$-trivial, the function $g\mapsto\|c(g)\|$ is constant on double cosets in $K\backslash G/K$.
    First assume $c$ to be proper and fix $\kappa>0$.
    The set $\{g\in G\mid \|c(g)\|\leq\kappa\}$ is then compact by hypothesis.
    Our earlier observation now tells us that there is a finite set $\mathcal{F}\subset K\backslash G/K$ such that $\|c(g)\|>\kappa$ whenever $KgK\notin\mathcal{F}$. 
    Equivalently, $\|\tilde{c}(tu_g)\|=\|t\|_2\|c(g)\|>\kappa\|t\|_2$ for $t\in T, g\in\Gamma$ whenever $KgK\notin\mathcal{F}$.
    But since double cosets correspond to irreducible bifinite $T$-$T$-bimodules in the Jones tower of $T\subset S$, this precisely means that $\tilde{c}$ is proper.
    Conversely, if $\tilde{c}$ is proper and $\kappa>0$, the set $\{g\in\Gamma\mid \|\tilde{c}(u_g)\|\leq\kappa\}$ covers only finitely many double cosets $K\backslash G/K$.
    It follows that $\{g\in G\mid \|c(g)\|\leq \kappa\}$ is covered by the same double cosets, and hence compact, as required.
\end{proof}
\end{proposition}

This dictionary of 1-cohomology allows us to easily transfer various properties between $T\subset S$ and $G$, which we do in the proposition formulated below.
All of these can be obtained through other means, see e.g.\@ \cite[Remark~8.11]{psv-cohom} for the statement about amenability, and \cite{eymard,larsen-palma,ad-hal,popa-annals-01} for related results concerning amenability, property~(T) and the Haagerup property in this setting.
\begin{proposition}
\label{thm:cohom-prop-equiv}
Let $\Lambda<\Gamma$ be a Hecke pair with Schlichting completion $K<G$, where $\Gamma$ acts on a type $\mathrm{II}_1$ factor $P$ by outer automorphisms.
Assume that $[\Gamma:\Lambda]=\infty$.
Then $G$ has $\mathcal{P}$ whenever $S=P\rtimes\Gamma$ has $\mathcal{P}$ relative to $T=P\rtimes\Lambda$, where $\mathcal{P}$ is either property (T), the Haagerup property or amenability.
\begin{proof}
    As before, put $\mathcal{S}=\mathrm{QN}_S(T)$.
    By the Delorme--Guichardet theorem, property (T) for $G$ is equivalent to $H^1_c(G;\mathcal{H})\cong H^1(T\subset\mathcal{S}; L^2(S)\htimes\mathcal{H})$ vanishing for all continuous unitary right representations $G\curvearrowright\mathcal{H}$.
    But if $H^1(T\subset\mathcal{S};-)$ vanishes on all $S$-$S$-bimodules of the form $L^2(S)\htimes\mathcal{H}$, it must vanish on any $S$-$S$-bimodule, as explained in \thref{rem:schlichting-completion-cohomology-red}.
    This in turn is equivalent to the relative property (T) for $T\subset S$ as a consequence of \cite[Theorem~9.16]{psv-cohom}.

    As for the Haagerup property, we know from \cite[Theorem~9.16]{psv-cohom} that $T\subset S$ has the Haagerup property if and only if there exists a proper 1-cocycle $c:\mathcal{S}\to\mathcal{K}$ for some $S$-$S$-bimodule $\mathcal{K}$.
    Again appealing to \thref{rem:schlichting-completion-cohomology-red}, we can take $\mathcal{K}$ to be of the form $L^2(S)\htimes\mathcal{H}$ for some Hilbert space $\mathcal{H}$ carrying a unitary continuous right $G$-representation.
    By \thref{thm:cohom-correspondence}, $c$ induces a proper cocycle on $G$ with values in $\mathcal{H}$.
    Conversely, if $c:G\to\mathcal{H}$ is a proper 1-cocycle, replacing it by a cohomologous $K$-trivial 1-cocycle preserves the properness, so we can again appeal to \thref{thm:cohom-correspondence} to obtain a proper 1-cocycle taking values in $L^2(S)\htimes\mathcal{H}$.

    Finally, \cite[Theorem~9.16]{psv-cohom} tells us that $T\subset S$ is amenable if and only if there is an approximately inner but non-inner cocycle with values in $L^2(S)\cftimes{T}L^2(S)\cong L^2(S)\htimes\ell^2(K\backslash G)$.
    Invoking \thref{thm:cohom-correspondence} once more, this is the same as $G$ admitting an approximately inner but non-inner $K$-trivial cocycle with values in $\ell^2(K\backslash G)\subset L^2(G)$.
    By \cite[Corollaire~III.2.4]{guichardet}, this is equivalent to the amenability of $G$.
    Conversely, if $G$ is amenable, the $G$-representation $\ell^2(K\backslash G)$ weakly contains the trivial representation by \cite[Proposition~3.4]{ad-hal}.
    Appealing to \cite[Corollaire~III.2.3]{guichardet}, this means that $\ell^2(K\backslash G)$ admits an approximately inner but non-inner 1-cocycle $c:G\to \ell^2(K\backslash G)$.
    Since this cocycle is cohomologous to a $K$-trivial 1-cocycle with the same properties, the converse implication follows after one final appeal to \thref{thm:cohom-correspondence}.
\end{proof}
\end{proposition}

\begin{example}
    The above equivalences highlight another way in which the behaviour of $P\rtimes\Lambda\subset P\rtimes\Gamma$ and $G$ is subtly different from that of $\mathcal{C}_f(K<G)$.
    Indeed, as shown in \cite[\S~4]{arano-vaes}, the conditions for $\mathcal{C}_f(K<G)$ to have property (T) or the Haagerup property are stronger than the corresponding property for $G$.
    In particular, the Hecke pair $\mathrm{SL}_2(\ZZ)<\mathrm{SL}_2(\ZZ[1/p])$ leads to a Schlichting completion $G=\mathrm{PSL}_2(\QQ_p)$ with the Haagerup property, but $\mathcal{C}_f(K<G)$ does \textit{not} have the Haagerup property \cite[see][Proposition~4.2]{arano-vaes}.

    Additionally, there are locally compact totally disconnected groups $G$ with property (T) such that the associated rigid C*-tensor categories $\mathcal{C}_f(K<G)$ do not have property (T) \cite[see e.g.][Example~4.3]{arano-vaes}.
\end{example}

\section{Concluding remarks}
Putting our vanishing result in the more general framework of quasi-regular inclusions through the lens of \thref{rem:almost-normal-qr-version}, we get the following very appealing generalisation question.
\begin{question}
    Suppose that we are given $N\in\NN$ and a tower $T\subset S\subset R$ of type $\mathrm{II}_1$ factors such that $T\subset S$, $T\subset R$ and $S\subset R$ are all quasi-regular, with $T\subset S$ and $T\subset R$ additionally unimodular.
    If $\beta^{(2)}_n(T\subset S)=0$ for all $n\leq N$, does it follow that $\beta^{(2)}_n(T\subset R)=0$ for all $n\leq N$?
\end{question}
In the case $P\subset P\rtimes\Lambda\subset P\rtimes\Gamma$ with $\Lambda<\Gamma$ almost normal and $\Gamma$ acting outerly on $P$, this recovers the vanishing theorem in the discrete group setting.

In the same setting, supposing that $\Lambda<\Gamma$ is unimodular, one a priori has three sets of $L^2$-Betti numbers associated with the Hecke pair $\Lambda<\Gamma$:
\begin{itemize}
    \item the $L^2$-Betti numbers of the Schlichting completion $G$\,;
    \item the $L^2$-Betti numbers of $P\rtimes\Lambda\subset P\rtimes\Gamma$\,;
    \item the $L^2$-Betti numbers of $\mathcal{C}_f(K<G)$, where $K<G$ is the Schlichting inclusion of $\Lambda<\Gamma$.
\end{itemize}
Moreover, by \cite[Proposition~2.7]{arano-vaes} the category of bifinite $(P\rtimes\Lambda)$-$(P\rtimes\Lambda)$-bimodules generated by the $(P\rtimes\Lambda)$-$(P\rtimes\Lambda)$-bifinite subbimodules of $L^2(P\rtimes\Gamma)$ is unitarily monoidally equivalent to $\mathcal{C}_f(K<G)$.
Our results imply that the former two sequences always agree, while the latter always vanishes (even irrespective of any unimodularity hypothesis on $\Lambda<\Gamma$).
This naturally leads one to ask what the fundamental difference is between the respective tube algebras of $P\rtimes\Lambda\subset P\rtimes\Gamma$ and $\mathcal{C}_f(K<G)$.
As shown in \cite{arano-vaes}, the tube algebra of $\mathcal{C}_f(K<G)$ is given by an explicit subalgebra of $L^\infty(G)\rtimes G$, where $G$ acts on itself by conjugation.
It is however unclear what one should expect for $P\rtimes\Lambda\subset P\rtimes\Gamma$.

\providecommand{\bysame}{\leavevmode\hbox to3em{\hrulefill}\thinspace}
\providecommand{\href}[2]{#2}
\setlength{\bibsep}{0.0pt}

\end{document}